%% file: root.tex
\documentclass[11pt]{amsart}

\usepackage{geometry}                
\usepackage{graphicx}

\usepackage{natbib} 		

\usepackage{amsmath,amssymb,amsthm}

\usepackage{pdflscape} 

\usepackage{epstopdf}
\usepackage{hyperref}

\hypersetup{
  colorlinks,
  bookmarksnumbered,
  bookmarkstype={toc},
  bookmarksopen={true},
  bookmarksopenlevel={1},
  pdfstartview={FitH},
  citecolor={blue},
  linkcolor={blue},
  urlcolor={blue},
  pdfpagemode={UseOutlines},
  breaklinks=true
 }

\newtheorem{theorem}{Theorem}[section]
\newtheorem{lemma}[theorem]{Lemma}
\newtheorem{proposition}[theorem]{Proposition}
\newtheorem{corollary}[theorem]{Corollary}

\newtheorem{assumption}{Assumption}[section]

\newtheorem{example}{Example}
\newtheorem{claim}{Claim}

\newcommand{\N}{\ensuremath{\mathbb{N}}}

\newcommand{\Z}{\ensuremath{\mathbb{Z}}}

\newcommand{\R}{\ensuremath{\mathbb{R}}}

\newcommand{\order}{\ensuremath{\mathcal{O}}}

\newcommand{\abs}[1]{\ensuremath{\left\lvert #1 \right\rvert}}
\newcommand{\norm}[1]{\ensuremath{\lvert\lvert #1\rvert\rvert}}

\newcommand{\bsym}[1]{\ensuremath{\boldsymbol{#1}}}
\renewcommand{\vec}[1]{\ensuremath{\mathbf{#1}}}
\newcommand{\set}[1]{\ensuremath{\mathcal{#1}}}

\newcommand{\inv}{\ensuremath{^{-1}}}
\newcommand{\Prob}{\ensuremath{\mathbb{P}}}
\newcommand{\Expect}{\ensuremath{\mathbb{E}}}

\newcommand{\diag}{\ensuremath{\mathrm{diag}}}
\newcommand{\sign}{\ensuremath{\mathrm{sign}}}

\newcommand{\esssup}{\ensuremath{\mathrm{ess}\sup}}

\newcommand{\mxm}{\ensuremath{\text{maximize}}}

\newcommand{\wrt}{\ensuremath{\mathrm{with \; respect \; to}}}

\newcommand{\ith}{\ensuremath{^\text{th}}}

\title[Fixed-Point Approaches to Computing Bertrand-Nash Equilibrium Prices]{Fixed-Point Approaches to Computing Bertrand-Nash Equilibrium Prices Under Mixed-Logit Demand: A Technical Framework for Analysis and Efficient Computational Methods.}

\author{W. Ross Morrow}
\address{Departments of Mechanical Engineering and Economics, Iowa State University, Ames IA 50011}

\author{Steven J. Skerlos}
\address{Department of Mechanical Engineering, University of Michigan, Ann Arbor MI 48104}

\date{\today}     

\begin{document}


\maketitle

	\tableofcontents

\input{intro}

\input{framework}

%
%
	\input{extendednummethods}

\input{gnhs}

\input{othermethods}

	\section{Acknowledgements}
	This research was supported by the National Science Foundation, the University of Michigan Transportation Research Institute's Doctoral Studies Program, and a research fellowship at the Belfer Center for Science and International Affairs at the Harvard Kennedy School. Both authors wish to thank Walter McManus, Brock Palen, Divakar Viswanath, Erin MacDonald, the editors, and the three anonymous reviewers for their contributions to this research. We would also like to acknowledge helpful suggestions offered by Fred Feinberg, John Hauser, Meredith Fowlie, Kenneth Judd, and the support of Kelly Sims-Gallagher and Henry Lee at the Belfer Center. Three anonymous reviewers, an associate editor, and Duncan Simester at {\em Operations Research} provided extremely helpful suggestions. 
	
	\bibliographystyle{chicago}
	
	\bibliography{root}

\end{document}

%% file: intro.tex
\section{Introduction}
\label{SEC:Intro}
	
	Bertrand competiton has been a prominent paradigm for the empirical study of differentiated product markets for at least twenty years. Firms engaged in Bertrand competition maximize profits by choosing prices for portfolios of differentiated products, and Bertrand-Nash equilibrium prices simultaneously maximize profits for all firms. Models combining Bertrand competition with the Mixed Logit discrete choice model of consumer demand have been used to study the automotive industry, electronics, entertainment, and food products and services; see \cite{Dube02}. 
	
	Many applications of Bertrand competition rely on counterfactual experiments: exercises in which hypothetical market conditions are simulated with an estimated model. Such experiments have been used to study corporate mergers \citep{Nevo00a}, novel products and  services \citep{Petrin02, Goolsbee04, Beresteanu08}, store locations \citep{Thomadsen05}, and regulatory policy changes \citep{Goldberg95, Goldberg98, Beresteanu08}. By definition, simulating market outcomes in counterfactual experiments requires computing equilibrium prices after changing the values of exogenous variables such as the number of firms or the products offered. Numerical methods for computing equilibrium prices have not yet received a thorough treatment in the literature, which currently focuses on model specification and estimation; see \cite{Knittel08, Dube08, Su08} for recent developments in estimation. \cite{Morrow10} fills this gap with a detailed investigation of four approaches for computing Bertrand-Nash equilibrium prices in single-period, multi-firm models with Mixed Logit demand. This working paper provides most of the technical background for that investigation. 
	
	
	Applying Newton's method to some form of the first-order or ``simultaneous stationarity'' condition is currently the {\em de facto} approach for computing equilibrium prices; see, for example, \cite{Nevo97, Nevo00a, Petrin02, Smith04, Doraszelski06, Jacobsen06}. Newton's method applied directly to the first-order condition may converge when started at observed prices if changes in exogenous variables have a marginal impact on equilibrium prices. However, when the changes to exogenous variables imply significant changes in product prices Newton's method applied directly to the first-order conditions may fail to compute equilibrium prices. Furthermore analyses that do not have observed prices to use as an initial guess will require methods with greater reliability. 
	
	\cite{Morrow10} demonstrate that solving fixed-point equations equivalent to the first-order condition for equilibrium is more reliable and efficient than solving the first-order condition itself. One fixed-point equation equivalent to the first-order conditions is the BLP-markup equation popularized by \cite{Berry95}. A second fixed-point equation, here termed the $\bsym{\zeta}$-markup equation, is a novel way to write the same condition on markups. Both markup equations lead to more robust numerical methods than found with a simple application of Newton's method to the first-order condition. Using the fixed-point expressions in this way can be considered ``nonlinearly'' or ``analytically'' pre-conditioning the first-order condition satisfied by equilibrium prices, a technique well-known in applied mathematics \citep{Brown90, Cai02}. 
	
	The existence of fixed-point equations for equilibrium suggests applying fixed-point iteration \citep{Judd98} to compute equilibrium prices, instead of Newton's method. The BLP-markup equation does not appear to be well-suited to fixed-point iteration. Example \ref{EX:ConvExam} in \cite{Morrow10} provides a case in which iterating on the BLP-markup equation is not necessarily locally convergent, while iterating on the $\bsym{\zeta}$-markup equation is superlinearly locally convergent. Iterating on the $\bsym{\zeta}$-markup equation also eliminates the need to solve linear systems, required to implement Newton's method and to iterate on the BLP-markup equation. This property makes fixed-point steps based on the $\zeta$-markup equation very inexpensive relative to Newton steps, an essential property to obtaining fast computations from generally linearly convergent fixed-point iterations. 
	
	
	Besides Newton's method and fixed-point iteration, few other practical approaches to the computation of equilibrium prices exist. Variational formulations, widely applied in economic and engineering problems \citep{Ferris97}, contain many solutions that need not be equilibria of the original problem. Explicit least-square minimization or Gauss-Newton methods can also be implemented, but are computational disadvantages relative to applications of standard Newton-type methods for nonlinear systems. Some authors apply tattonement $-$ iterating on a game's best response correspondence $-$ to compute equilibrium in prices or other strategic variables including product mix \citep{Choi90}, product characteristics \citep{CBO03, Austin05, Bento05}, and engineering variables \citep{Michalek04}. Tattonement, however, has three issues: it requires the iterative computation of profit-optimal prices (a special case of the problem discussed in this article), should be inefficient relative to direct methods whenever optimal strategies are coupled, and lacks the global convergence guarantees of contemporary Newton solvers. Section \ref{ECSEC:OtherMethods} reviews these conclusions in more detail. 
	
	This article should be viewed as a companion to \cite{Morrow10}; some of our notation and text may seem out of place without first reviewing that article. In several places, text from \cite{Morrow10} is repeated. 
	

%% file: framework.tex
	
\section{A Technical Framework}
\label{ECSEC:Framework}

	This section describes the mathematical framework employed in \cite{Morrow10}. Several key assumptions are introduced and summarized in Table \ref{TAB:KeyAssumptions}. 
	
	\begin{table}
		\caption{List of important assumptions used in this section.}
		\label{TAB:KeyAssumptions}
		\begin{tabular}{cl}
		Assumption & Purpose \\ \hline
		\ref{MixedLogitUtilityAssump}
			& To provide a general form for utility functions \\
		\ref{BoundedProfitsAssump} 
			& To ensure profits are bounded and vanish as prices increase without bound \\
		\ref{LeibnizRuleCondition}
			& To ensure the Leibniz Rule holds, validating Eqn. (\ref{Decomposition}) \\
		\ref{BoundedAssumption}
			& To ensure that $\bsym{\eta}$ is bounded. Implies the coercivity of $\vec{F}_\eta,\vec{F}_\zeta$ \\
			& and the existence of simultaneously stationary prices. \\
		\ref{BoundedOmegaProductsAssumption}
			& To ensure that $\bsym{\zeta}$ is bounded. Implies the coercivity of $\vec{F}_\zeta$ \\
			& and the existence of simultaneously stationary prices. \\
		\ref{LVA}
			& To ensure that the derivatives of profit vanish as prices increase without bound \\
		\ref{BoundedAssumption2}
			& To ensure the coercivity of $\vec{F}_\eta,\vec{F}_\zeta$ under weaker conditions than \\
			& Assumption \ref{BoundedAssumption}. \\
		\hline
		\end{tabular}
	\end{table}
	

	\subsection{Mathematical Notation} 
	\label{SUBSECNotation}
	
	\subsubsection{Sets.} Table \ref{TAB:Sets} lists some important sets and the symbols used for them. $\N$ denotes the natural numbers $\{1,2,\dotsc\}$, and $\N(N)$ denotes the natural numbers up to $N$, that is, $\N(N) = \{1,\dotsc,N\}$. $\R$ denotes the set of real numbers $(-\infty,\infty)$, $[0,\infty)$ denotes the non-negative real numbers, and $[0,\infty]$ denotes the extended non-negative half-line. We denote the $(J-1)$-dimensional simplex $\{ (x_1,\dotsc,x_N) \in [0,1]^N : \sum_{n=1}^N x_n = 1 \}$ by $\mathbb{S}(N)$, and the $J$-dimensional ``pyramid'' $\{ (x_1,\dotsc,x_N) \in [0,1]^N : \sum_{n=1}^N x_n \leq 1 \}$ by $\triangle(J)$. Hyper-rectangles in $\R^N$, i.e. sets of the form $[a_1,b_1] \times \dotsb \times [a_N,b_N]$ for some $a_n,b_n \in \R$ with $a_n < b_n$ for all $n \in \N(N)$, are denoted by $[\vec{a},\vec{b}]$ where $\vec{a} = (a_1,\dotsc,a_N)$ and $\vec{b} = (b_1,\dotsc,b_N)$. $\set{P}$ always denotes the non-negative numbers: $\set{P} = [0,\infty)$. For other sets, we typically use calligraphic upper case letters such as ``$\set{A}$''. For any set $\set{A}$, $\abs{\set{A}}$ denotes its cardinality. For any $\set{B} \subset \set{A}$, $\set{A} \setminus \set{B}$ denotes the set $\{ b \in \set{A} : b \notin \set{B} \}$. 
	
	\begin{table}
		\caption{Important sets.}
		\label{TAB:Sets}
		\begin{tabular}{rcll}
			\multicolumn{3}{c}{Symbol} & Description \\ \hline
			$\N$ & $=$ & $\{1,2,\dotsc\}$ & Natural numbers \\
			$\R$ & $=$ & $(-\infty,\infty)$ & Real numbers \\
			$\set{P}$ & $=$ & $[0,\infty)$ & Non-negative real numbers \\
			$\set{J}$ & $=$ & $\{1,\dotsc,J\}$ & Set of product indices \\
			$\set{X}$ & $\subset$ & $\R^K$ & Set of product characteristics \\
			$\set{T}$ & $\subset$ & $\R^L$ & Set of individual characteristics \\ 
			\hline
		\end{tabular}
	\end{table}
	
	\subsubsection{Symbols.} Table \ref{TAB:Symbols} itemizes specific symbols used in the text.

	\begin{table}
		\caption{Summary of important symbols.}
		\label{TAB:Symbols}
		\begin{tabular}{rllll}
		\multicolumn{3}{c}{Symbol} & Description & Defined in \\ \hline \\
		\multicolumn{5}{l}{Products (see Section \ref{SUBSEC:PPCP})} \\ \hline
		$J$ 
			& $\in$ & $\N$ 
			& number of products \\
		$K$ 
			& $\in$ & $\N$ 
			& number of non-price product characteristics \\
		$\vec{x}_j$ 
			& $\in$ & $\set{X}$ 
			& non-price characteristics of product $j$ \\
		$p_j$ 
			& $\in$ & $\set{P}$ 
			& price of product $j$ \\
		$\vec{p}$ 
			& $\in$ & $\set{P}^J$ 
			& vector of all product prices \\
		\\ \multicolumn{5}{l}{Individual Characteristics (see Section \ref{SUBSEC:PPCP})} \\ \hline
		$\bsym{\theta}$ 
			& $\in$ & $\set{T}$ 
			& individual characteristics, including observed \\
			&&& demographics and ``random coefficients'' \\
		$\mu$ 
			& $-$ & $-$ 
			& distribution of individual characteristics \\
		\\ \multicolumn{5}{l}{Choice Probabilities (see Section \ref{SUBSEC:PPCP})} \\ \hline
		$u_j(\bsym{\theta},p_j)$ 
			& $\in$ & $[-\infty,\infty)$ 
			& utility of product $j$ \\
		$\vartheta(\bsym{\theta})$
			& $\in$ & $[-\infty,\infty)$ 
			& utility of the outside good \\
		$P_j^L(\bsym{\theta},\vec{p})$
			& $\in$ & [0,1] 
			& Logit choice probability for product $j$
			& Eqn. (\ref{LogitChoiceProb}) \\
		$P_j(\vec{p})$
			& $\in$ & [0,1] 
			& Mixed Logit choice probability for product $j$ \\
		$\vec{P}(\vec{p})$
			& $\in$ & $[0,1]^J$ 
			& vector of Mixed Logit choice probabilities for all \\
			&&& products \\
		\\ \multicolumn{5}{l}{Firms, Costs, Profits, and Stationarity (see Section \ref{SUBSEC:Profits}, \ref{SUBSEC:EquilibriumAndTheSSCs})} \\ \hline
		$F$ 
			& $\in$ & $\N$ 
			& number of firms \\
		$\set{J}_f$ 
			& $\subset$ & $\set{J}$ 
			& indices of the products offered by firm $f$ \\
		$c_j$ 
			& $\in$ & $\set{P}$ 
			& (fixed) unit cost of product $j$ \\
		$\vec{c}$ 
			& $\in$ & $\set{P}^J$ 
			& vector of all (fixed) unit costs \\
		$\hat{\pi}_f(\vec{p})$
			& $\in$ & $\R$
			& expected profits for firm $f$
			& Eqn. (\ref{ExpectedProfits}) \\
		$(D_k\hat{\pi}_f)(\vec{p})$
			& $\in$ & $\R$
			& derivative of firm $f$'s profits, with respect to the
			& Eqn. (\ref{EQN:ProfitGradientFormula}) \\
			&&& price of product $k$ \\
		$(\tilde{\nabla}\hat{\pi})(\vec{p})$
			& $\in$ & $\R^J$
			& Combined Gradient of profits
			& Prop. \ref{SSCProp}, Eqn. (\ref{SSCs}) \\
		\\ \multicolumn{5}{l}{Choice Probability Derivatives (see Sections \ref{SUBSEC:EquilibriumAndTheSSCs}, \ref{SUBSEC:ZetaMap})} \\ \hline
		$(D_kP_j)(\vec{p})$
			& $\in$ & $\R$
			& derivative of product $j$'s choice probability 
			& \\
			&&& with respect to the price of product $k$ \\
		$(\tilde{D}\vec{P})(\vec{p})$
			& $\in$ & $\R^{J \times J}$
			& ``intra-firm'' Jacobian matrix of the choice
			& Eqn. (\ref{IntraFirmJacobian}) \\
			&&& probability vector \\
		$\bsym{\Lambda}(\vec{p})$, $\tilde{\bsym{\Gamma}}(\vec{p})$
			& $\in$ & $\R^{J \times J}$
			& matrices appearing in our decomposition of $(\tilde{D}\vec{P})(\vec{p})$
			& Eqn. (\ref{Decomposition}), \\
			&&&& 
		\\ \multicolumn{5}{l}{Fixed-Point Equations (see Sections \ref{SUBSEC:MarkupEqn}, \ref{SUBSEC:ZetaMap})} \\ \hline
		$\bsym{\eta}(\vec{p})$
			& $\in$ & $\R^J$
			& the BLP-markup function \citep{Berry95}
			& Eqn. (\ref{MarkupEquation}) \\
		$\bsym{\zeta}(\vec{p})$
			& $\in$ & $\R^J$
			& our $\bsym{\zeta}$-markup function
			& Eqn. (\ref{ZetaMap}) \\
		\hline
		\end{tabular}
	\end{table}
	
	Bold, un-italicized symbols (e.g., ``$\vec{x}$'') denote vectors and matrices; typically we reserve lower case letters to refer to vectors and use upper case letters to refer to matrices; the vector of choice probabilities ``$\vec{P}$'' is an exception. Throughout we use $\vec{1}$ to denote a vector of ones of the appropriate size for the context in which it appears. $\vec{I}$ always denotes the identity matrix of a size appropriate for the context. For any $\vec{x} \in \R^N$, $\diag(\vec{x})$ denotes the $N \times N$ diagonal matrix whose diagonal is $\vec{x}$. Any vector inequalities between vectors are to be taken componentwise: for example, $\vec{x} < \vec{y}$ means $x_n < y_n$ for all $n$. 
	
	Random variables are denoted with capital letters ``$X$'', with random vectors being denoted with bold capital letters (e.g., ``$\vec{Q}$''). While this overlaps with our notation for matrices, it should not cause any confusion. $\Prob$ denotes a probability and $\Expect$ denotes an expectation. $\esssup_\mu f$ denotes the essential supremum of the (measurable) function $f$ over $\set{T}$, with respect to the measure $\mu$; see, e.g., \cite{Bartle66}. 
	
	$\log$ always denotes the natural (base $e$) logarithm. We use the ``Big-O'' notation $\order(g)$ as follows: If there exists some $M < \infty$ such that $\lim_{p \to q} [ f(p)/g(p) ] \leq M$, we say $f \in \order(g)$; the point $q$ is left implicit. 
	
	\subsubsection{Differentiation.} Our conventions for denoting differentiation follow \cite{Munkres91}. We use the symbol ``$D$'' to denote differentiation using subscripts to invoke additional specificity. Letting $\vec{f} : \R^M \to \R^N$, $(D_m f_n)(\vec{x})$ denotes the derivative of the $n\ith$ component function with respect to the $m\ith$ variable and $(D\vec{f})(\vec{x})$ is the $N \times M$ derivative matrix of $\vec{f}$ at $\vec{x}$ with components $( (D\vec{f})(\vec{x}) )_{n,m} = (D_m f_n)(\vec{x})$. Thus for $f:\R^M \to \R$, $(Df)(\vec{x})$ is a row vector. If $f : \R^M \to \R$, we define the gradient $(\nabla f)(\vec{x}) \in \R^M$ as the transposed derivative: $(\nabla f) (\vec{x}) = (Df)(\vec{x})^\top$. 
	
	
	\subsection{Consumers, Products, and Choice Probabilities}
	\label{SUBSEC:PPCP}
	
	A collection of $F \in \N$ firms offer a total of $J \in \N$ products to a population of individuals (or households). Each product $j \in \set{J} = \{1,\dotsc,J\}$ is defined by a price, $p_j \in \set{P} = [0,\infty)$, and a vector of $K \in \N$ product ``characteristics'' $\vec{x}_j \in \set{X} \subset \R^K$. Individuals are identified by a vector of characteristics $\bsym{\theta}$ from some set $\set{T}$. These individual characteristics can include both observed demographics and ``random coefficients''  \citep{Berry95,Nevo00b,Train03} that characterize unobserved individual-specific heterogeneity with respect to preference for product characteristics. The relative density of individual characteristic vectors in the population is described by a probability distribution $\mu$ over $\set{T}$. 
	
	An individual identified by $\bsym{\theta} \in \set{T}$ receives the (random) utility 
	\begin{equation*}
		U_j(\bsym{\theta},\vec{x}_j, p_j) = u(\bsym{\theta},\vec{x}_j, p_j) + \set{E}_j
	\end{equation*}
	from purchasing product $j \in \set{J}$, and 
	\begin{equation*}
		U_0(\bsym{\theta}) = \vartheta(\bsym{\theta}) + \set{E}_0
	\end{equation*}
	for forgoing purchase of any of these products; i.e. ``purchasing the outside good.'' Individuals choose the ``product'' $j \in \{0,\dotsc,J\}$ with maximum utility. Here $u : \set{T} \times \set{X} \times \set{P} \to [-\infty,\infty)$ is a systematic utility function, $\vartheta : \set{T} \to \R$ is a valuation of the no-purchase option or ``outside good,'' and $\bsym{\set{E}} = \{ \set{E}_j \}_{j = 0}^J$ is a random vector of i.i.d. standard extreme value variables. Section \ref{SUBSEC:Utilities} below gives a general specification of utility functions appropriate for equilibrium pricing. The basic requirements are that $u$ is continuously differentiable and strictly decreasing in price, and without lower bound as prices increase. 
	
	Demand for each product $j$ is characterized by choice probabilities $P_j : \set{P}^J \to [0,1]$ derived from (random) utility maximization. Given the distributional assumption on $\bsym{\set{E}}$, the choice probabilities for an individual characterized by $\bsym{\theta} \in \set{T}$ are those of the Logit model \citep[Chapter 3]{Train03}:
	\begin{equation}
		\label{LogitChoiceProb}
		P_j^L(\bsym{\theta},\vec{p})
			= \frac{ e^{u_j(\bsym{\theta},p_j)} }
				{ e^{\vartheta(\bsym{\theta})} + \sum_{k=1}^J e^{u_k(\bsym{\theta},p_k)} }. 
	\end{equation}
	The vector $\vec{p} \in \set{P}^J$ denotes the vector of all product prices. Product-specific utility functions $u_j : \set{T} \times \set{P} \to [-\infty,\infty)$ for all $j$, defined by $u_j(\bsym{\theta},p) = u(\bsym{\theta},\vec{x}_j,p)$ for all $(\bsym{\theta},p) \in \set{T} \times \set{P}$, are used in Eqn. (\ref{LogitChoiceProb}) and in the following sections. The Mixed Logit choice probabilities $P_j(\vec{p}) = \int P_j^L(\bsym{\theta},\vec{p}) d\mu(\bsym{\theta})$ follow from integrating over the distribution of individual characteristics \citep[Chapter 6]{Train03}. The vector of Mixed Logit choice probabilities for all products is denoted by $\vec{P}(\vec{p}) \in [0,1]^J$. 
	
	The examples below review several instances of this choice model. Examples \ref{EX:BM80Example} and \ref{EX:BLPExample} are used in \cite{Morrow10}. Example \ref{EX:NevoExample} illustrates the type general specifications used in estimation. Example \ref{EX:Simulation} describes one kind of ``simulation'' of a Mixed Logit model \citep{Train03}. 
	
	\begin{example} \citep{Boyd80}
		\label{EX:BM80Example}
		Take $\set{T} = \set{P} \times \R^K$, denoting $\bsym{\theta} = (\alpha,\bsym{\beta})$ for $\alpha \in \set{P}$ and $\bsym{\beta} \in \R^K$. Set $u(\alpha,\bsym{\beta},\vec{x},p) = - \alpha p + \bsym{\beta}^\top\vec{x}$ and $\vartheta(\alpha,\bsym{\beta}) = -\infty$ for all $(\alpha,\bsym{\beta}) \in \set{P} \times \R^K$. $\mu$ is defined by specifying that $\alpha$ and $\bsym{\beta}$ are independently lognormally distributed (with appropriately chosen signs, means, and variances).
	\end{example}
	
	\begin{example} \citep{Berry95}
		\label{EX:BLPExample}
		Take $\set{T} = \set{P} \times \R^K \times \R$, denoting $\bsym{\theta} = (\phi,\bsym{\beta},\beta_0)$ for $\phi \in \set{P}$, $\bsym{\beta} \in \R^K$, and $\beta_0 \in \R$. Set 
		\begin{equation*}
			u(\phi,\bsym{\beta},\vec{x},p) 
				= \left\{ \begin{aligned}
					&\alpha \log( \phi - p ) + \bsym{\beta}^\top\vec{x} 
						&&\quad\text{if } p < \phi \\
					&-\infty &&\quad\text{otherwise}
				\end{aligned} \right .
			\quad\quad\text{and}\quad\quad
			\vartheta(\phi,\beta_0)
				= \alpha \log \phi + \beta_0
		\end{equation*}
		for some fixed coefficient $\alpha > 0$. $\phi$ represents income and is given a lognormal distribution, while the random coefficients $\bsym{\beta},\beta_0$ are independently normally distributed with some mean and variance. Note that income ($\phi$) serves as an upper bound on the price an individual can pay for {\em any} product.
		
	\end{example}
	
	\begin{example} \citep{Nevo00b}
		\label{EX:NevoExample}
		Take $\set{T} = \set{P} \times \R^D \times \R^{K+2}$, denoting $\bsym{\theta} = (\phi,\vec{d},\bsym{\nu})$ for $\phi \in \set{P}$, $\vec{d} \in \R^D$, and $\bsym{\nu} \in \R^{K+2}$. Again, $\phi$ represents income; $\vec{d} \in \R^D$ represents a vector of $D$ observed demographic variables (which may include income); $\bsym{\nu} \in \R^{K+2}$ represents a vector of $K+2$ random coefficients: one for each product characteristic, one for price, and one for the outside good. Set 
		\begin{align*}
			u(\phi,\vec{d},\bsym{\nu},\vec{x},p) 
				&= ( \alpha + \bsym{\pi}_p^\top\vec{d} + \bsym{\sigma}_p^\top \bsym{\nu} ) ( \phi - p ) 
					+ \left( \bsym{\beta} + \bsym{\Pi}\vec{d} + \bsym{\Sigma}\bsym{\nu} \right)^\top\vec{x} \\
			\vartheta(\phi,\vec{d},\bsym{\nu})
				&= ( \alpha + \bsym{\pi}_p^\top\vec{d} + \bsym{\sigma}_p^\top \bsym{\nu} ) \phi 
						+ \bsym{\pi}_0^\top\vec{d} + \bsym{\sigma}_0^\top\bsym{\nu}
		\end{align*}
		where $\alpha \in \R$, $\bsym{\beta} \in \R^K$, $\bsym{\pi}_p,\bsym{\pi}_0 \in \R^{D}$, $\bsym{\Pi} \in \R^{K \times D}$, $\bsym{\sigma}_p,\bsym{\sigma}_0 \in \R^{K+2}$, and $\bsym{\Sigma} \in \R^{K \times (K+2)}$ are coefficients. The distribution of $\vec{d}$ is estimated from available data (e.g., Census data) and $\bsym{\nu}$ is assumed to be {\em standard} independent multivariate normal. When $\alpha + \bsym{\pi}_p^\top\vec{d} + \bsym{\sigma}_p^\top \bsym{\nu}$, the coefficient on price, is positive, an individual prefers higher prices. 
		
		\cite{Petrin02} and \cite{Berry04a} adopt similar specifications that eliminate this counterintuitive property. \cite{Petrin02} takes the price component of utility to be $\alpha(\phi) \log( \phi - p )$, where $\alpha : \set{P} \to \set{P}$ is a step function. \cite{Berry04a} take the price component of utility to be $\alpha p$, but define $\alpha = - e^{-( \alpha + \bsym{\pi}_p^\top\vec{d} + \bsym{\sigma}_p^\top \bsym{\nu} )}$.
		
	\end{example}
	
	
	\begin{example} 
		\label{EX:Simulation}
		{\em (Simulation)}. Take any of the examples above, and draw $S \in \N$ vectors $\bsym{\theta}_s \in \set{T}$ according to the distribution $\mu$. Let $\set{T}^\prime = \{ \bsym{\theta}_s \}_{s=1}^S$ and define a probability measure $\mu^\prime$ over $\set{T}^\prime$ by $\mu^\prime(\bsym{\theta}_s) = 1/S$ for all $s$. Then $(u,\vartheta,\set{T}^\prime,\mu^\prime)$ defines a {\em simulator} of the ``full'' Mixed Logit model with $(u,\vartheta,\set{T},\mu)$; see \cite{Train03}. These approximations are essential in estimation of Mixed Logit models and in computations of equilibrium prices.
	\end{example}
	
	
	
	\subsection{Utility Specification}
	\label{SUBSEC:Utilities}
	
	This section presents a generalization of the systematic utility functions used in the examples given in the text, a specification closely related to the one introduced by \cite{Caplin91}. \cite{Morrow08, Morrow08a} use a similar specification to analyze equilibrium prices in simple Logit models. 
	\begin{assumption}
		\label{MixedLogitUtilityAssump}
		For all $j$, there exist functions $w_j : \set{T} \times \set{P} \to [-\infty,\infty)$ and $v_j : \set{T} \to (-\infty,\infty)$ such that the systematic utility function $u_j : \set{T} \times \set{P} \to [-\infty,\infty)$ can be written $u_j(\bsym{\theta},p) = w_j(\bsym{\theta},p) + v_j(\bsym{\theta})$. Furthermore there exists $\varsigma : \set{T} \to (0,\infty]$ such that $w_j : \set{T} \times [0,\infty] \to [-\infty,\infty)$ satisfies, for all $j$ and $\mu$-almost every (a.e.) $\bsym{\theta} \in \set{T}$, 
		\begin{itemize}
			\item[(a)] $w_j(\bsym{\theta},\cdot) : (0,\varsigma(\bsym{\theta})) \to [-\infty,\infty)$ is continuously differentiable, strictly decreasing, and finite
			\item[(b)] $w_j(\bsym{\theta},p) = -\infty$ for all $p \geq \varsigma(\bsym{\theta})$, and 
			\item[(c)] $w_j(\bsym{\theta},p) \downarrow -\infty$ as $p \uparrow \varsigma(\bsym{\theta})$. 
		\end{itemize}
		$v_j : \set{T} \to (-\infty,\infty)$ is arbitrary. 
	\end{assumption}
	
	Note that we have not restricted $\mu$, the distribution of individual characteristics, with Assumption \ref{MixedLogitUtilityAssump}. Important examples of $\mu$ from the econometrics and marketing literature include finitely supported distributions (often empirical frequency distributions for integral observed demographic variables), standard continuous distributions (e.g. normal, lognormal and $\chi^2$), truncated standard continuous distributions, finite mixtures of standard continuous distributions, and independent products of any of these types of distributions. This generality allows us to address a wide variety of otherwise disparate examples with a single notation. In particular, this generality allows us to use a single framework to treat both ``full'' Mixed Logit models defined by some $\mu$ with uncountable support and simulation-based approximations to such models. 
	
	Some existing empirical specifications violate Assumption \ref{MixedLogitUtilityAssump} by admiting positive price coefficients for $\bsym{\theta} \in \set{T}^\prime \subset \set{T}$, where $\set{T}^\prime \subset \set{T}$ has nonzero $\mu$-measure. See, for example, \cite{Nevo00a} (Example \ref{EX:NevoExample}) or \cite{Brownstone00}. This implies that $w(\bsym{\theta},\cdot)$ is increasing on $\set{T}^\prime$. If $w(\bsym{\theta},\cdot)$ is {\em not} decreasing for $\mu$-a.e. $\bsym{\theta}$, or at least {\em eventually} decreasing for $\mu$-a.e. $\bsym{\theta}$ in the sense that there are always prices large enough to ensure that $w(\bsym{\theta},\cdot)$ is decreasing for $\mu$-a.e. $\bsym{\theta}$, then profit-optimal pricing is not a well-posed problem and finite equilibrium prices will not exist. 
	
	The variable $\Sigma = \varsigma(\bsym{\theta})$ represents an individual-specific reservation price. As in the \cite{Berry95} model of Example \ref{EX:BLPExample}, this reservation price is most often derived from household or individual income. Correspondingly, $\Sigma$ is often given a lognormal distribution to (roughly) fit empirical income data. In principle, this reservation price could be related to purchasing power derived from observed demographic variables other than income, or unobserved demographic variables such as family wealth. Thus we allow this reservation price to be specified as a function of all ``demographic'' characteristics, $\bsym{\theta}$. Conditions (b) and (c) in Assumption \ref{MixedLogitUtilityAssump} imply that the probability an individual characterized by $\bsym{\theta}$ will purchase a product is zero for any price above $\varsigma(\bsym{\theta})$ {\em and} vanishes as the price approaches $\varsigma(\bsym{\theta})$. We set $\varsigma_* = \esssup \varsigma$ and allow, but do not require, $\varsigma_* = \infty$. For example, simulation-based approximations to the \citeauthor{Berry95} demand model have $\varsigma_* < \infty$, as can be easily checked. 
	
	Note also that Condition (c) in Assumption \ref{MixedLogitUtilityAssump} ensures the continuity of $P_j^L(\bsym{\theta},\vec{p})$ at any vector of prices with some component equal to $\varsigma(\bsym{\theta})$. We must require this of the {\em Logit} choice probabilities to obtain {\em Mixed Logit} choice probabilities that are continuous on $(0,\varsigma_*)^J$ for the important class of simulation-based approximations with finitely supported $\mu$. Continuous Logit choice probabilities also imply continuous Mixed Logit choice probabilities, by the Dominated Convergence Theorem. 
	
	
	\subsection{Profits}
	\label{SUBSEC:Profits}
	
	To describe the optimal pricing problems faced by each firm we use the following notation. Let $F \in \N$ denote the number of firms. For each $f \in \{ 1 , \dotsc , F\}$, there exists a set $\set{J}_f \subset \set{J}$ of indices that corresponds to the $J_f = \abs{\set{J}_f}$ products offered by firm $f$. The collection of all these sets, $\{\set{J}_f\}_{f=1}^F$, forms a partition of $\set{J}$. Subsequently, in writing ``$f(j)$'' for some $j \in \set{J}$, we mean the unique $f \in \{1,\dotsc,F\}$ such that $j \in \set{J}_f$. The vector $\vec{p}_f \in \R^{J_f}$ refers to the vector of prices of the products offered by firm $f$. Negative subscripts denote competitor's variables as in, for instance, $\vec{p}_{-f} \in \R^{J_{-f}}$, where $J_{-f} = \sum_{g \neq f} J_g$, is the vector of prices for products offered by all of firm $f$'s competitors. Firm-specific choice probability functions are denoted by $\vec{P}_f(\vec{p}) \in \R^{J_f}$. 
	
	Two additional assumptions are required to complete the definition of firms' profits in a manner consistent with empirical applications of Bertrand competition. First, we must specify unit and fixed costs: for each product $j$ there exists a unit cost $c_j \in \set{P}$ and for each firm there exists a fixed cost $c_f^F \in \set{P}$. Both $c_j$ and $c_f^F$ depend only on the collection of product characteristics chosen by the firm, and not on the quantity sold by the firm during the purchasing period for the reasons discussed below. We let $\vec{c}_f \in \set{P}^{J_f}$ denote the vector of unit costs for the products offered by firm $f$, and $\vec{c} \in \set{P}^J$ denote the vector of unit costs for all products. 
	
	Second, Bertrand competition entails the following ``comittment'' assumption on the quantities produced \citep{Baye08}. Let $Q_j(\vec{p})$ denote the (random) quantity of product $j$ that the population will demand during the purchasing period, given prices for all products $\vec{p}$. These random demands are derived from random utility maximization. We assume each firm commits to producing exactly $Q_j(\vec{p})$ units of each product $j \in \set{J}_f$ during the purchasing period. This implies either that there are no production capacity constraints that limit a firm's ability to meet {\em any} demands that arise during the purchase period, or that production backlogs do not affect demand. 
	
	With the commitment and constant costs assumptions, the total cost firm $f$ incurs in producing (and selling) $Q_j(\vec{p})$ units of product $j$ during the purchasing period are given by the random variable 
	\begin{equation*}
		\sum_{j \in \set{J}_f} c_j Q_j(\vec{p}) + c_f^F. 
	\end{equation*}
	Random revenues are, of course, given by $\sum_{j \in \set{J}_f} Q_j(\vec{p})p_j$. The random variable 
	\begin{equation*}
		\Pi_f(\vec{p}) 
			= \sum_{j \in \set{J}_f} Q_j(\vec{p}) ( p_j - c_j ) - c_f^F
	\end{equation*}
	then gives firm $f$'s (random) profits for the purchasing period as a function of all product prices. Following most of the theoretical and empirical literature in both marketing and economics, we assume that firms take expected profits, 
	\begin{equation}
		\label{ExpectedProfits}
		\Expect[\Pi_f(\vec{p})] 
			= I\hat{\pi}_f(\vec{p}) - c_f^F
		\quad\text{where}\quad
		\hat{\pi}_f(\vec{p}) 
			= \sum_{j\in\set{J}_f} P_j(\vec{p})(p_j - c_j)
	\end{equation}
	as the metric by which they optimize their pricing decisions in this stochastic optimization problem. Here $I \in \N$ denotes the number of individuals in the population. 
	
	Eqn. (\ref{ExpectedProfits}) demonstrates that neither the total firm fixed costs $c_f^F$ nor the population size $I$ play a role in determining the prices that maximize expected profits under the assumptions described above. Henceforth we focus on the ``population-normalized gross expected profits'' $\hat{\pi}_f(\vec{p})$, referred to in the text and below as simply ``profits''. Firms thus solve
	\begin{equation}
		\label{EQN:Opt}
		\begin{aligned}
			\mxm &\quad \hat{\pi}_f(\vec{p}) = \sum_{j \in \set{J}_f} P_j(\vec{p})(p_j - c_j) \\
			\wrt &\quad \vec{p}_f \in \set{P}^{J_f}
		\end{aligned}
	\end{equation}
	
	Before continuing with our framework, we discuss quantity-dependent costs and clarify when profits are bounded. 
	
	
	\subsubsection{Quantity-Dependent Costs}
	
	Including costs that depend on quantities produced is certainly possible, though this should introduce extra terms into the first-order equations presented below (Eqn. (\ref{SSCs})). Generally speaking, unit costs that depend on the quantity produced would be expressed as $c_j : \Z_+ \to \set{P}$, and unit costs that depended on the {\em expected} quantity produced would be expressed as $c_j : \set{P} \to \set{P}$. If unit costs depend on the quantity produced, then product $j$'s unit costs for the purchasing period (i) are random and (ii) depend on prices. To see this, simply note that product $j$'s unit costs for the purchasing period are $c_j(Q_j(\vec{p}))$. Assuming quantity-dependent costs also obscures expected profits, since there are now nonlinear terms $Q_j(\vec{p})c_j(Q_j(\vec{p}))$ in the formula for random profits. If unit costs depend only on the expected quantity produced, then unit costs are not random but still depend on prices: $c_j(\Expect[Q_j(\vec{p})]) = c_j(IP_j(\vec{p}))$. In either case the derivatives of unit costs with respect to prices should appear in the first-order conditions. This is acknowledged in the theoretical literature. As these terms have not yet been included in the empirical literature, even when costs are assumed to depend on quantities produced \citep{Berry95,Petrin02}, we focus on costs that are independent of the quantity produced. 
	
	
	\subsubsection{Bounded and Vanishing Profits}

	Here we present a technical assumption that ensures that profits are not only bounded, but vanish as all prices approach $\varsigma_*$. 
	
	\begin{assumption}
		\label{BoundedProfitsAssump}
		For all $j$ there exists some $r_j : \set{T} \to (1,\infty)$ and some $\bar{p}_j : \set{T} \to \set{P}$ satisfying
		\begin{equation}
			\label{PCond}
			\sup \Big\{ \; p \mu(\{ \bsym{\theta} : \bar{p}(\bsym{\theta}) > p \}) \; : \; p \in (0,\varsigma_*) \; \Big\} < \infty. 
		\end{equation} 
		such that 
		\begin{equation}
			\label{UtilityBound}
			u_j(\bsym{\theta},p) \leq -r_j(\bsym{\theta}) \log p + \vartheta(\bsym{\theta})
		\end{equation}
		for all $p \geq \bar{p}_j(\bsym{\theta})$, $\mu$-a.e. 
	\end{assumption}
	
	\begin{lemma}
		If Assumption \ref{BoundedProfitsAssump} holds, then $\hat{\pi}_f(\vec{p})$ is bounded in $\vec{p}$ and vanishes as $\vec{p}_f \to \varsigma_* \vec{1} \in \R^{J_f}$. 
	\end{lemma}
	\proof{}
		We use the Dominated Convergence Theorem. Eqn. (\ref{UtilityBound}) ensures that $p_j P_j^L(\bsym{\theta},\vec{p})$ vanishes $\mu$-a.e. as $p_j \uparrow \varsigma_*$; see also \cite{Morrow08a}. Eqn. (\ref{PCond}) ensures that
		\begin{equation*}
			\hat{\pi}_f(\vec{p})
				= \sum_{j \in \set{J}_f} \int p_j P_j^L(\bsym{\theta},\vec{p}) d\mu(\bsym{\theta})
					- \sum_{j \in \set{J}_f} c_j \int P_j(\vec{p}) 
		\end{equation*}
		is bounded as prices approach $\varsigma_*$, as we now show. 
		
		The key quantities in this integral are 
		\begin{equation*}
			\int p_j P_j^L(\bsym{\theta},\vec{p}) d\mu(\bsym{\theta})
				\leq \int p_j \left( 
						\frac{ e^{u_j(\bsym{\theta},p_j) - \vartheta(\bsym{\theta})} }
							{ 1 + e^{u_j(\bsym{\theta},p_j) - \vartheta(\bsym{\theta})} } \right) 
						d\mu(\bsym{\theta}); 
		\end{equation*}
		the $c_j P_j(\vec{p})$ terms vanish if $p_j \uparrow \varsigma_*$ since $P_j(\vec{p})$ vanishes. We must show that these terms are bounded as $p_j \uparrow \varsigma_*$. By assumption, 
		\begin{equation*}
			p_j \left( \frac{ e^{u_j(\bsym{\theta},p_j) - \vartheta(\bsym{\theta})} }
							{ 1 + e^{u_j(\bsym{\theta},p_j) - \vartheta(\bsym{\theta})} } \right) 
				\leq \left( \frac{ 1 }{ p_j } \right)^{r_j(\bsym{\theta})-1}
					\left( \frac{ 1 }{ 1 + e^{u_j(\bsym{\theta},p_j) - \vartheta(\bsym{\theta})} } \right) 
		\end{equation*}
		for all $p_j \geq \bar{p}_j(\bsym{\theta})$. Thus we write
		\begin{align*}
			\int p_j \left( 
						\frac{ e^{u_j(\bsym{\theta},p_j) - \vartheta(\bsym{\theta})} }
							{ 1 + e^{u_j(\bsym{\theta},p_j) - \vartheta(\bsym{\theta})} } \right) 
						d\mu(\bsym{\theta})
				&= \int_{ \{ \bsym{\theta} : p_j < \bar{p}_j(\bsym{\theta}) \} } p_j \left( 
						\frac{ e^{u_j(\bsym{\theta},p_j) - \vartheta(\bsym{\theta})} }
							{ 1 + e^{u_j(\bsym{\theta},p_j) - \vartheta(\bsym{\theta})} } \right) 
						d\mu(\bsym{\theta}) \\
				&\quad\quad \quad\quad \quad\quad
					+ \int_{ \{ \bsym{\theta} : p_j \geq \bar{p}_j(\bsym{\theta}) \} } p_j \left( 
						\frac{ e^{u_j(\bsym{\theta},p_j) - \vartheta(\bsym{\theta})} }
							{ 1 + e^{u_j(\bsym{\theta},p_j) - \vartheta(\bsym{\theta})} } \right) 
						d\mu(\bsym{\theta}) \\
				&\leq p_j \mu{ \{ \bsym{\theta} : p_j < \bar{p}_j(\bsym{\theta}) \} } \\
				&\quad\quad \quad\quad \quad\quad
					+ \int_{ \{\bsym{\theta} : p_j \geq \bar{p}_j(\bsym{\theta}) \} } 
						\left( \frac{ 1 }{ p_j } \right)^{r_j(\bsym{\theta})-1}
						d\mu(\bsym{\theta}). 
		\end{align*}
		By Eqn. (\ref{PCond}), the first term is bounded. We take $p_j > 1$, without loss of generality, so that $1/p^{r_j(\bsym{\theta})-1} \leq 1$ for $\mu$-a.e. $\bsym{\theta}$ and the second term is bounded.
	\endproof
	
	We now make some remarks regarding Assumption \ref{BoundedProfitsAssump}. 
	
	Note that if $\varsigma(\bsym{\theta}) < \infty$ then Eqn. (\ref{UtilityBound}) holds for {\em any} $r(\bsym{\theta}) > 1$ by taking $\bar{p}(\bsym{\theta}) = \varsigma(\bsym{\theta})$. If $\varsigma(\bsym{\theta}) = \infty$, Eqn. (\ref{UtilityBound}) admits any utility function $u(\bsym{\theta},\cdot)$ that is (eventually) concave in price. 
	
	If $\varsigma_* < \infty$, then $\varsigma(\bsym{\theta}) < \infty$ for $\mu$-a.e. $\bsym{\theta}$. Furthermore, Eqn. (\ref{PCond}) is trivial. 
	
	To further analyze Eqn. (\ref{PCond}), we assume $\varsigma_* = \infty$. We define $Z = \bar{p}(\bsym{\Theta})$, where $\bsym{\Theta}$ is the $\set{T}$-valued random variable with $\Prob( \bsym{\Theta} \in \set{A} ) = \mu(\set{A}) = \int_{\set{A}} d\mu(\bsym{\theta})$. If $\varsigma(\bsym{\theta}) < \infty$ for $\mu$-a.e. $\bsym{\theta}$, then we can take $Z = \Sigma = \varsigma(\bsym{\Theta})$. Eqn. (\ref{PCond}) can be re-written as $\sup\{ p \Prob( Z > p ) : p \in (0,\infty) \} < \infty$, or equivalently $\lim_{ p \to \infty} [ p \Prob( Z > p ) ] < \infty$. Eqn. (\ref{PCond}) admits any $Z$ with finite expectation, and even admits any $Z$ with a ``fat-tailed'' distribution satisfying $p^{1+\beta}\Prob( Z > p) \to 1$ as $p \to \infty$ for some $\beta > 0$. Eqn. (\ref{PCond}) can be written $\Prob(Z > p) = \order(1/p)$. 
	

	\subsection{Local Equilibrium and the Simultaneous Stationarity Conditions}
	\label{SUBSEC:EquilibriumAndTheSSCs}
	
	Assuming that the choice probabilities are continuously differentiable in prices, at equilibrium each firm's prices satisfy the stationarity condition 
	\begin{equation}
		\label{EQN:ProfitGradientFormula}
		(D_k \hat{\pi}_f)(\vec{p})
			= \sum_{j \in \set{J}_f} (D_kP_j)(\vec{p})(p_j - c_j) + P_k(\vec{p})
		\quad\quad\text{for all } k \in \set{J}_f. 
	\end{equation}
	
	Combining the stationarity condition for each firm we obtain the {\em Simultaneous Stationarity Condition}, a first-order (necessary) condition for local equilibrium prices. 
	
	\begin{proposition}[Simultaneous Stationarity Condition]
		\label{SSCProp}
		Suppose $\vec{P}$ is continuously differentiable. Let $(\tilde{\nabla}\hat{\pi})(\vec{p}) \in \R^J$ denote the ``combined gradient'' with components $((\tilde{\nabla}\hat{\pi})(\vec{p}))_j = (D_j\hat{\pi}_{f(j)})(\vec{p})$ where $f(j)$ denotes the index of the firm offering product $j$. If $\vec{p}$ is a local equilibrium, then 
		\begin{equation}
			\label{SSCs}
			(\tilde{\nabla}\hat{\pi})(\vec{p})
				= (\tilde{D}\vec{P})(\vec{p})^\top (\vec{p} - \vec{c}) + \vec{P}(\vec{p}) 
				= \vec{0}. 
		\end{equation}
		where $(\tilde{D}\vec{P})(\vec{p}) \in \R^{J \times J}$ is the ``intra-firm'' Jacobian matrix of price derivatives of the choice probabilities defined by
		\begin{equation}
			\label{IntraFirmJacobian}
			\big( (\tilde{D}\vec{P})(\vec{p}) \big)_{j,k} 
				= \left\{ \begin{aligned}
					&(D_kP_j)(\vec{p})	
						&&\quad\text{if products $j$ and $k$ are offered by the same firm } \\
					&\quad\quad 0		
						&&\quad\text{otherwise }
				\end{aligned} \right.
		\end{equation}
		Prices $\vec{p}$ satisfying Eqn. (\ref{SSCs}) are called ``simultaneously stationary.'' 
	\end{proposition}
	
	The matrix $- (\tilde{D}\vec{P})(\vec{p})$ has previously been denoted by ``$\triangle$'' \citep{Berry95,Petrin02,Beresteanu08}, ``$\bsym{\Omega}$'' \citep{Nevo00a}, and ``$\bsym{\Phi}$'' \citep{Dube02}. We prefer the ``$D$'' notation to emphasize the relationship of $(\tilde{D}\vec{P})(\vec{p})$ to the Jacobian matrix of the choice probabilities $\vec{P}$, while using the superscript ``$\sim$'' to denote the intra-firm sparsity structure. 
	
%
	
	A set of simultaneously stationary prices are a local equilibrium {\em only} if every firm's profits are locally maximized at those prices; this can be verified by confirming that every firm's profits are locally concave (Section \ref{SUBSECSufficiency}). Note that there is no convenient condition to verify that every firm's profits are {\em globally} maximized at a particular local equilibrium. That is, there is no convenient condition to ensure that certain prices are a proper equilibrium. 
	
	
	\subsection{Choice Probability Derivatives}
	\label{SUBSECChoiceProbabilityDerivatives}
	
	In this section we examine the price derivatives of Mixed Logit choice probabilities. In what follows, $(Dw_j)(p_j)$ denotes the derivative of the price component of utility, $w_j$, with respect to price. 
	\begin{proposition}
		\label{MixedLogitDerivatives}
		Fix $\vec{p} \in (0,\varsigma_*)^J$, let $u_j$ be given as in Assumption \ref{MixedLogitUtilityAssump} for all $j$, and suppose the Leibniz Rule holds for the Mixed Logit choice probabilities $P_j(\vec{p}) = \int P_j^L(\bsym{\theta},\vec{p}) d\mu(\bsym{\theta})$; that is, $(D_k P_j)(\vec{p}) = \int (D_jP_k^L)(\bsym{\theta},\vec{p}) d\mu(\bsym{\theta})$. Then the Jacobian matrix of $\vec{P}$ is given by 
		\begin{equation}
			\label{Decomposition}
			(D\vec{P})(\vec{p}) = \bsym{\Lambda}(\vec{p}) - \bsym{\Gamma}(\vec{p})
		\end{equation}
		where $\bsym{\Lambda}(\vec{p}) \in \R^{J \times J}$ is the diagonal matrix with diagonal entries 
		\begin{equation*}
			\lambda_j(\vec{p}) = \int_{\set{L}(p_j)} (Dw_j)(\bsym{\theta},p_j) P_j^L(\bsym{\theta},\vec{p}) d\mu(\bsym{\theta}),
			\quad\quad
			\set{L}(p) = \{ \bsym{\theta} : \varsigma(\bsym{\theta}) > p \}
		\end{equation*}
		and $\bsym{\Gamma}(\vec{p})$ is the full $J \times J$ matrix with entries
		\begin{equation*}
			\gamma_{j,k}(\vec{p}) 
				= \int_{\set{G}(p_j,p_k)} P_j^L(\bsym{\theta},\vec{p})P_k^L(\bsym{\theta},\vec{p})(Dw_k)(\bsym{\theta},p_k) d\mu(\bsym{\theta}),
			\quad\quad
			\set{G}(p,q) = \set{L}(p) \cap \set{L}(q).
		\end{equation*}
		The intra-firm price derivatives of the Mixed Logit choice probabilities are given by $(\tilde{D}\vec{P})(\vec{p}) = \bsym{\Lambda}(\vec{p}) - \tilde{\bsym{\Gamma}}(\vec{p})$ where $\big( \tilde{\bsym{\Gamma}}(\vec{p}) \big)_{j,k} = \gamma_{j,k}(\vec{p})$ if $f(j) = f(k)$ and $\big( \tilde{\bsym{\Gamma}}(\vec{p}) \big)_{j,k} = 0$ otherwise. 
	\end{proposition}
	\proof{}
		We first characterize the Logit choice probabilities. For all $j,k$ we have
		\begin{align*}
			(D_kP_j^L)(\bsym{\theta},\vec{p})
				&= P_j^L(\bsym{\theta},\vec{p}) ( \delta_{j,k} - P_k^L(\bsym{\theta},\vec{p}) )(Dw_k)(\bsym{\theta},p_k) \\
				&= \delta_{j,k} P_k^L(\bsym{\theta},\vec{p}) (Dw_k)(\bsym{\theta},p_k)
					- P_j^L(\bsym{\theta},\vec{p}) P_k^L(\bsym{\theta},\vec{p}) (Dw_k)(\bsym{\theta},p_k)
		\end{align*}
		for any $\bsym{\theta} \in \set{L}(p_k)$ and $(D_kP_j^L)(\bsym{\theta},\vec{p}) = 0$ for any $\bsym{\theta} \in  \{ \bsym{\theta}^\prime \in \set{T} : p_k > \varsigma(\bsym{\theta}^\prime)\}$ (because $P_j^L(\bsym{\theta},\cdot)$ is identically zero in a neighborhood of $\vec{p}$). Neglecting values $\bsym{\theta} \in \varsigma\inv(p_k)$ for the moment, we observe that these formulae and the Leibniz rule generate the desired expression for the Mixed Logit choice probabilities. 
		
		We complete the proof by considering $\bsym{\theta} \in \varsigma\inv(p_k)$. If $\varsigma\inv(p_k)$ has $\mu$-measure zero for any $p_k$, then we do not need to worry about Logit choice probability derivatives at $\bsym{\theta} \in \varsigma\inv(p_k)$. On the other hand if $\varsigma\inv(p_k)$ has positive $\mu$-measure for some $p_k$, we must assume continuity of the Logit choice probability derivatives: i.e. $(D_kP_j^L)(\bsym{\theta},\vec{p}) \to 0$ as $p_k \uparrow \varsigma(\bsym{\theta})$. Otherwise, the Logit choice probability derivative is not defined on a set of demographics with positive measure.
	\endproof
	
	$\bsym{\lambda}$ is closely related to a familiar economic quantity. Recall that the ``inclusive value,'' or expected maximum utility, conditional on demographics is given by \citep{Small81,Train03}
	\begin{equation*}
		\iota^L(\bsym{\theta},\vec{p}) 
			= \log \left( e^{\vartheta(\bsym{\theta})} + \sum_{j=1}^J e^{u_j(\bsym{\theta},p_j)} \right)
	\end{equation*}
	It is easy to see that $\lambda_k$ is the derivative of the ``aggregate inclusive value'' $\iota(\vec{p}) = \int \iota^L(\bsym{\theta},\vec{p}) d\mu(\bsym{\theta})$ with respect to the $k\ith$ price: $\lambda_k(\vec{p}) = (D_k\iota)(\vec{p}) = \int (D_k\iota^L)(\bsym{\theta},\vec{p}) d\mu(\bsym{\theta})$. 

	Note that $\bsym{\Gamma}(\vec{p})$ and $\tilde{\bsym{\Gamma}}(\vec{p})$ are not necessarily symmetric for all $\vec{p}$. If $(Dw_k)(\bsym{\theta},p)$ is independent of both $k$ and $p$, as in the case of the \cite{Boyd80} model presented in Example \ref{EX:BM80Example} above, then $\bsym{\Gamma}(\vec{p})$ (and thus $\tilde{\bsym{\Gamma}}(\vec{p})$) is symmetric for all $\vec{p}$. On the other hand if $(Dw_k)(\bsym{\theta},\cdot)$ is independent of $k$ and strictly monotone in $p$, as is the case of the strictly concave in price utility from \cite{Berry95}, then $\gamma_{j,k}(\vec{p}) = \gamma_{k,j}(\vec{p})$ if and only if $p_j = p_k$. 
	
	The following assumption gives a simple, abstract condition on $(\vec{u},\vartheta,\mu)$ that guarantees the Leibniz Rule holds and defines continuously differentiable choice probabilities. 
	\begin{assumption}
		\label{LeibnizRuleCondition}
		Let $k$ be arbitrary and define $\psi_k : \set{T} \times \set{P} \to \set{P}$ by 
		\begin{equation*}
			\psi_k(\bsym{\theta},p) 
				= \left\{ \begin{aligned}
					&\abs{(Dw_k)(\bsym{\theta},p)}
						\left( \frac{ e^{u_k(\bsym{\theta},p)} }
							{ e^{\vartheta(\bsym{\theta})} + e^{u_k(\bsym{\theta},p)} } \right)
					&&\quad\text{if } p < \varsigma(\bsym{\theta})\\
					&\quad\quad\quad\quad\quad
						0 &&\quad\text{if } p \geq \varsigma(\bsym{\theta})\\ 
				\end{aligned} \right .
		\end{equation*}
		Assume (i) $\psi_k(\bsym{\theta},\cdot) : (0,\varsigma_*) \to \set{P}$ is continuous for $\mu$-a.e. $\bsym{\theta} \in \set{T}$; that is, $\psi_k(\bsym{\theta},q) \to \psi_k(\bsym{\theta},p)$ as $q \to p$ for any $p \in (0,\varsigma_*)$. (ii) $\psi_k(\cdot,p^\prime) : \set{T} \to \set{P}$ is uniformly $\mu$-integrable for all $p^\prime$ in some neighborhood of any $p \in (0,\varsigma_*)$; that is, there exists some $\varphi : \set{T} \to [0,\infty)$ with $\int \varphi (\bsym{\theta}) d\mu(\bsym{\theta}) < \infty$ (that may depend on $k$ and $p$), such that $\psi_k(\bsym{\theta},p^\prime) \leq \varphi (\bsym{\theta})$ for all $p^\prime$ in some neighborhood of $p$. 
	\end{assumption}
	Note that under Assumption \ref{MixedLogitUtilityAssump}, (i) requires only that $\psi_k(\bsym{\theta},p) \to 0$ as $p \uparrow \varsigma(\bsym{\theta})$ for $\mu$-a.e. $\bsym{\theta}$. 
	
	\begin{proposition}
		\label{LeibnizRule}
		If Assumption \ref{LeibnizRuleCondition} holds, then the Leibniz Rule holds for the Mixed Logit choice probabilities which are also continuously differentiable on $(0,\varsigma_*)^J$. 
	\end{proposition}
	
	\proof{}
		Taking for granted that $(D_kP_j^L)(\bsym{\theta},\cdot)$ is continuous at $\vec{p}$ and the differences
		\begin{equation}
			\label{Differences}
			h\inv \big( P_j^L(\bsym{\theta},\vec{p} + h\vec{e}_k) - P_j^L(\bsym{\theta},\vec{p}) \big)
		\end{equation}
		are uniformly $\mu$-integrable for small enough $h$, the Dominated Convergence Theorem implies that
		\begin{align*} 
			&\lim_{h \to 0} h\inv \Big( \int P_j^L(\bsym{\theta},\vec{p} + h\vec{e}_k) d\mu(\bsym{\theta})
									- \int P_j^L(\bsym{\theta},\vec{p}) d\mu(\bsym{\theta}) \Big) \\
			&\quad\quad\quad\quad
				= \lim_{h \to 0} \int h\inv \big( P_j^L(\bsym{\theta},\vec{p} + h\vec{e}_k) - P_j^L(\bsym{\theta},\vec{p}) \big) d\mu(\bsym{\theta}) \\
			&\quad\quad\quad\quad
				= \int \lim_{h \to 0} h\inv \big( P_j^L(\bsym{\theta},\vec{p} + h\vec{e}_k) - P_j^L(\bsym{\theta},\vec{p}) \big) d\mu(\bsym{\theta}) \\
			&\quad\quad\quad\quad
				= \int (D_kP_j^L)(\bsym{\theta},\vec{p}) d\mu(\bsym{\theta}). 
		\end{align*}
		This validates the Leibniz Rule. This proof is essentially that given in a general setting by \cite[Chapter 5, pg. 46]{Bartle66}.
		
		To complete the proof we must validate that $(D_kP_j^L)(\bsym{\theta},\cdot)$ is continuous in $p_k$ and the differences in Eqn. (\ref{Differences}) are uniformly $\mu$-integrable in a neighborhood of $p_k$. It is easy to see that the desired continuity follows from Assumption \ref{MixedLogitUtilityAssump} and Assumption \ref{LeibnizRuleCondition}, Condition (i). Specifically, note that $(D_kP_j^L)(\bsym{\theta},\vec{p}) = 0 = \psi_k(\bsym{\theta},p_k)$ for $\bsym{\theta} \in \{ \bsym{\theta}^\prime \in \set{T} : p_k > \varsigma(\bsym{\theta}^\prime)\}$ and
		\begin{align*}
			(D_kP_j^L)(\bsym{\theta},\vec{p})
				&= \big( \delta_{j,k} - P_j^L(\bsym{\theta},\vec{p}) \big) P_k^L(\bsym{\theta},\vec{p}) (Dw_k)(\bsym{\theta},p_k) \\
				&= \big( \delta_{j,k} - P_j^L(\bsym{\theta},\vec{p}) \big) 
						\left( \frac{ e^{\vartheta(\bsym{\theta})} + e^{u_k(\bsym{\theta},p_k)} }
								{ e^{\vartheta(\bsym{\theta})} + \sum_{i=1}^J e^{u_i(\bsym{\theta},p_i)} } \right)
						\left( \frac{ e^{u_k(\bsym{\theta},p_k)} }
								{ e^{\vartheta(\bsym{\theta})} +e^{u_k(\bsym{\theta},p_k)} } \right) (Dw_k)(\bsym{\theta},p_k) \\
				&= \big( \delta_{j,k} - P_j^L(\bsym{\theta},\vec{p}) \big) 
						\left( \frac{ e^{\vartheta(\bsym{\theta})} + e^{u_k(\bsym{\theta},p_k)} }
								{ e^{\vartheta(\bsym{\theta})} + \sum_{i=1}^J e^{u_i(\bsym{\theta},p_i)} } \right)
						\psi_k(\bsym{\theta},p_k)
		\end{align*}
		for $\bsym{\theta} \in \set{L}(p_k)$. Suppose $p_k = \varsigma(\bsym{\theta})$. By Assumption \ref{MixedLogitUtilityAssump} (a) and (b), the first two terms are continuous. By Assumption \ref{MixedLogitUtilityAssump} (c), 
		\begin{align*}
			\lim_{\vec{q} \to \vec{p}, q_k < \varsigma(\bsym{\theta})} (D_kP_j^L)(\bsym{\theta},\vec{p})
				&= \left( \delta_{j,k} - \frac{ e^{u_j(\bsym{\theta},p_j)} }
								{ e^{\vartheta(\bsym{\theta})} + \sum_{i \neq k} e^{u_i(\bsym{\theta},p_i)} } \right ) 
					\left( \frac{ e^{\vartheta(\bsym{\theta})}}
								{ e^{\vartheta(\bsym{\theta})} + \sum_{i \neq k} e^{u_i(\bsym{\theta},p_i)} } \right)
						\lim_{q_k \uparrow \varsigma(\bsym{\theta})} \psi_k(\bsym{\theta},q_k)
		\end{align*}
		Assumption \ref{LeibnizRuleCondition}, Condition (i) is then necessary for the continuity of $(D_kP_j^L)(\bsym{\theta},\vec{p})$ for all $j,k$ and $\vec{p} \in (\vec{0},\varsigma_*\vec{1})$. Specifically if $\psi_k(\bsym{\theta},\cdot)$ is discontinuous at $\varsigma(\bsym{\theta})$, then 
		\begin{align*}
			\lim_{\vec{q} \to \vec{p}, q_k < \varsigma(\bsym{\theta})} (D_kP_k^L)(\bsym{\theta},\vec{p})
				&= \left( \frac{ e^{\vartheta(\bsym{\theta})}}
								{ e^{\vartheta(\bsym{\theta})} + \sum_{i \neq k} e^{u_i(\bsym{\theta},p_i)} } \right)
						\lim_{q_k \uparrow \varsigma(\bsym{\theta})} \psi_k(\bsym{\theta},q_k)
		\end{align*}
		
		To prove the integrability, we first note that for all $j,k$ and $\vec{p}$ we have $\abs{ (D_kP_j^L)(\bsym{\theta},\vec{p}) } \leq \psi_k(\bsym{\theta},p_k)$. This bound is a consequence of the formula above, and is tight as $\vec{p}_{-k}$ varies. The mean value theorem for functions of a single real variable states that
		\begin{align*}
			h\inv( P_j^L(\bsym{\theta},\vec{p} + h\vec{e}_k) - P_j^L(\bsym{\theta},\vec{p}) )
				= (D_kP_j^L)(\bsym{\theta},\vec{p} + \eta\vec{e}_k)
		\end{align*}
		for some $\eta$ such that $\abs{\eta} < \abs{h}$, and thus
		\begin{align*}
			\abs{h}\inv\abs{ P_j^L(\bsym{\theta},\vec{p} + h\vec{e}_k) - P_j^L(\bsym{\theta},\vec{p}) }
				\leq \psi_k(\bsym{\theta},p_k+\eta)
				\leq \varphi(\bsym{\theta})
		\end{align*}
		for $\mu$-a.e. $\bsym{\theta} \in \set{T}$ and small enough $h$. Thus, the desired uniform $\mu$-integrability follows from Assumption \ref{LeibnizRuleCondition}, Condition (ii).
	\endproof
	
	An ``easier'' bound is simply $\abs{(D_kP_j^L)(\bsym{\theta},\vec{p})} \leq \abs{(Dw_k)(\bsym{\theta},p_k) }$, and thus we might consider changing the statement of Proposition \ref{LeibnizRule} to hypothesize only the uniform $\mu$-integrability of the utility price derivatives. In fact, this bound can be used to validate the Leibniz Rule for the \citeauthor{Boyd80} model of Example \ref{EX:BM80Example} that lacks an outside good. However, this bound fails to be useful for the \citeauthor{Berry95} model of Example \ref{EX:BLPExample}, since $w(p) = \alpha \log( \varsigma(\bsym{\theta}) - p )$ and $\abs{(Dw_k)(\bsym{\theta},p_k) } = \alpha / ( \varsigma(\bsym{\theta}) - p )$ is singular on $\varsigma\inv(p)$. In empirical applications, $\varsigma$ is onto, generating a singularity somewhere in $\set{T}$ for all $p$; this singularity cannot be ``controlled'' for all $p$ by choosing the measure $\mu$. In this case, a hypothesis only about the price derivatives of utility is not useful. 
	
	We close this section by stating some basic results concerning $(\tilde{D}\vec{P})(\vec{p})$ that are used below. 
	
	\begin{lemma}
		\label{LambdaNonsingular}
		Under Assumption \ref{MixedLogitUtilityAssump}, $P_j(\vec{p})$ and $\lambda_j(\vec{p})$ are never zero on $(0,\varsigma_*)^J$. Thus $\bsym{\Lambda}(\vec{p})$ is nonsingular for all $\vec{p} \in (0,\varsigma_*)^J$. 
	\end{lemma}
	\proof{}
		Note that $\set{L}(p_j)$ is nonempty and has positive $\mu$-measure, $P_j^L(\cdot,\vec{p})$ is strictly positive on $\set{L}(p_j)$, and $(Dw_j)(\cdot,p_j)P_j^L(\cdot,\vec{p})$ is strictly negative on $\set{L}(p_j)$. It follows that $P_j(\vec{p})$ and $\lambda_j(\vec{p})$ are nonzero.
	\endproof
	
%
%
%
	
	\begin{lemma}
		\label{LamGamLemma}
		Let $\vec{p} \in (0,\varsigma_*)^J$,  suppose $\vartheta : \set{T} \to (-\infty,\infty)$, and define
		\begin{align}
			\label{OmegaDef}
			\bsym{\Omega}_f(\vec{p}) 
				&= \bsym{\Lambda}_f(\vec{p})\inv \bsym{\Gamma}_f(\vec{p})^\top
				\in \R^{J_f \times J_f}
				\text{ for all } f \\
			\tilde{\bsym{\Omega}}(\vec{p}) 
				&= \bsym{\Lambda}(\vec{p})\inv \tilde{\bsym{\Gamma}}(\vec{p})^\top
				\in \R^{J \times J}. 
		\end{align}
		These matrices are well-defined by Lemma \ref{LambdaNonsingular}, and have the following properties:
		\begin{itemize}
		
			\item[(i)] $(D_f\vec{P}_f)(\vec{p})^\top = \bsym{\Lambda}_f(\vec{p}) ( \vec{I} - \bsym{\Omega}_f(\vec{p}) )$ and $(\tilde{D}\vec{P})(\vec{p})^\top = \bsym{\Lambda}(\vec{p}) ( \vec{I} - \tilde{\bsym{\Omega}}(\vec{p}) )$. 
			
			\item[(ii)] $\norm{ \bsym{\Omega}_f(\vec{p}) }_\infty < 1$ and $\norm{ \tilde{\bsym{\Omega}}(\vec{p}) }_\infty < 1$. 
			
			\item[(iii)] $\vec{I} - \bsym{\Omega}_f(\vec{p}) \in \R^{J_f \times J_f}$ and $\vec{I} - \tilde{\bsym{\Omega}}(\vec{p}) \in \R^{J \times J}$ are strictly diagonally dominant and nonsingular.  
			
			\item[(iv)] $(\vec{I} - \bsym{\Omega}_f(\vec{p}))\inv \in \R^{J_f \times J_f}$ and $(\vec{I} - \tilde{\bsym{\Omega}}(\vec{p}))\inv \in \R^{J \times J}$ map positive vectors to positive vectors. 
		\end{itemize}
	\end{lemma}
	
	\proof{}
	\begin{itemize}
		\item[(i)] This follows immediately from Prop. \ref{MixedLogitDerivatives}. 
		\item[(ii)] We note that
		\begin{align*}
			\omega_{k,l}(\vec{p})
				= \frac{ \gamma_{l,k}(\vec{p}) }{ \lambda_k(\vec{p}) }
				= \int_{\set{L}(p_k)} P_l^L(\bsym{\theta},\vec{p}) d\mu_{k,\vec{p}}(\bsym{\theta})
		\end{align*}
		where $\mu_{k,\vec{p}}$ is the probability distribution with density, with respect to $\mu$, given by
		\begin{equation*}
			d\mu_{k,\vec{p}}(\bsym{\theta})
				= \frac{ P_k^L(\bsym{\theta},\vec{p}) \abs{(Dw_k)(\bsym{\theta},p_k)} d\mu(\bsym{\theta}) }
									{ \int_{\set{L}(p_k)} P_k^L(\bsym{\phi},\vec{p}) \abs{(Dw_k)(\bsym{\phi},p_k)} d\mu(\bsym{\phi}) }. 
		\end{equation*}
		Thus $\bsym{\Lambda}_f(\vec{p})\inv \bsym{\Gamma}_f(\vec{p})^\top$ has row sums
		\begin{align*}
			\int \left( \sum_{j \in \set{J}_f} P_j^L(\bsym{\theta},\vec{p}) \right) d\mu_{k,\vec{p}}(\bsym{\theta}) < 1. 
		\end{align*}
		The additional assumption that $\vartheta : \set{T} \to (-\infty,\infty)$ plays a role in establishing this inequality because then there is always a set $\set{T}_k^\prime \subset \set{T}$ with $\mu_{k,\vec{p}}(\set{T}_k^\prime) > 0$ on which $\sum_{j \in \set{J}_f} P_j^L(\bsym{\theta},\vec{p}) < 1$. 
		
		\item[(iii)] The inequality
		\begin{align*}
			1 > \int \left( \sum_{j \in \set{J}_f} P_j^L(\bsym{\theta},\vec{p}) \right) d\mu_{k,\vec{p}}(\bsym{\theta})
		\end{align*}
		is equivalent to
		\begin{align*}
			\abs{ 1 - \omega_{k,k}(\vec{p}) }
				&= 1 - \int P_k^L(\bsym{\theta},\vec{p}) d\mu_{k,\vec{p}}(\bsym{\theta}) \\
				&> \int \left( \sum_{j \in \set{J}_f \setminus k } P_j^L(\bsym{\theta},\vec{p}) \right) d\mu_{k,\vec{p}}(\bsym{\theta})
				= \sum_{l \neq k} \omega_{k,l}(\vec{p}).  
		\end{align*}
		The claim follows. 
		
		\item[(iii)] Because $\bsym{\Omega}_f(\vec{p})$ maps positive vectors to positive vectors, so does its power series
		\begin{equation*}
			\sum_{n=1}^\infty \bsym{\Omega}_f(\vec{p})^n
				= (\vec{I} - \bsym{\Omega}_f(\vec{p})) \inv.
		\end{equation*}
		\end{itemize}
	\endproof
	
	\begin{corollary}
		\label{DPNonSingular}
		$(D_f\vec{P}_f)(\vec{p})^\top$ and $(\tilde{D}\vec{P})(\vec{p})^\top$ are strictly diagonally dominant and nonsingular for $\vec{p} \in (\vec{0},\varsigma_*\vec{1})$. 
	\end{corollary}
	\proof{}
		This follows directly from Lemma \ref{LamGamLemma}, claims (i) and (iii).
	\endproof
		
	
	\subsection{The BLP-Markup Equation}
	\label{SUBSEC:MarkupEqn}
	
	A prominent form of the first-order conditions Eqn. (\ref{SSCs}) is the {\em BLP-markup equation}:
	\begin{equation}
		\label{MarkupEquation}
		\vec{p} = \vec{c} + \bsym{\eta}(\vec{p})
		\quad\text{where}\quad
		\bsym{\eta}(\vec{p}) 
			= - (\tilde{D}\vec{P})(\vec{p})^{-\top} \vec{P}(\vec{p}). 
	\end{equation}
	Note that $\bsym{\eta}$ is defined for {\em any} continuously differentiable choice probabilities with nonsingular $(\tilde{D}\vec{P})(\vec{p})^\top$. We have shown above that this applies to certain Mixed Logit models (Section \ref{SUBSECChoiceProbabilityDerivatives}). Eqn. (\ref{MarkupEquation}) and Corollary \ref{DPNonSingular} show that $\bsym{\eta}$ is well-defined and continuous, at least for $\vec{p} \in (0,\varsigma_*)^J$. 
	
	Traditionally, the BLP-markup equation (\ref{MarkupEquation}) has been used to estimate {\em costs} assuming {\em observed} prices are in equilibrium via the formula $\vec{c} = \vec{p} - \bsym{\eta}(\vec{p})$; see, e.g., \cite{Berry95} or \cite{Nevo00a}. These cost estimates form the basis of counterfactual experiments with an estimated demand model. \cite{Beresteanu08} have recently suggested that the BLP-markup equation is also useful for computing equilibrium prices, a suggestion we explore further below. Note that the BLP-markup equation must be interpreted as a nonlinear fixed-point equation when applied to compute equilibrium prices. 
	
	
	We now derive several important properties of $\bsym{\eta}$ from an alternative form of Eqn. (\ref{MarkupEquation}) based on Lemma \ref{LamGamLemma}, valid when $\vec{p} \in (0,\varsigma_*)^J$:
	\begin{equation}
		\label{EtaEqn}
		\Big( \vec{I} - \tilde{\bsym{\Omega}}(\vec{p}) \Big) 
				\bsym{\eta}(\vec{p})
			= - \bsym{\Lambda}(\vec{p})\inv\vec{P}(\vec{p}).
	\end{equation}
	
	First, Eqn. (\ref{EtaEqn}) proves that profit-optimal markups are positive for the class of Mixed Logit models we consider, thanks to Lemma \ref{LamGamLemma}, claim (iv).  
	\begin{corollary}
		\label{PositiveMixedLogitMarkups}
		Suppose Assumptions \ref{MixedLogitUtilityAssump}-\ref{LeibnizRuleCondition} hold. Then $\bsym{\eta}(\vec{p}) > \vec{0}$ for all $\vec{p} \in (0,\varsigma_*)^J$. Hence if $\vec{p} \in (0,\varsigma_*)^J$ is a local equilibrium, then $\vec{p} > \vec{c}$. 
	\end{corollary}
	
	Second, Eqn. (\ref{EtaEqn}), rather than Eqn. (\ref{MarkupEquation}), should be used to actually compute $\bsym{\eta}$. Recall that $\kappa_2(\vec{A})$ denotes the 2-norm condition number of the matrix $\vec{A}$ \citep{Trefethen97}. 
	\begin{lemma}
		\label{EtaCondition}
		Suppose Assumptions \ref{MixedLogitUtilityAssump}-\ref{LeibnizRuleCondition} hold. Eqn. (\ref{EtaEqn}) is better conditioned than Eqn. (\ref{MarkupEquation}), in the sense that
		\begin{equation*}
			\kappa_2 \big( \vec{I} - \tilde{\bsym{\Omega}}(\vec{p}) \big)
				\leq \left( \frac{ \min_j \abs{ \lambda_j(\vec{p}) }  }
						{ \max_j \abs{ \lambda_j(\vec{p}) } } \right)
					\kappa_2 \big( (\tilde{D}\vec{P})(\vec{p})^\top \big)
		\end{equation*}
		for all $\vec{p} \in (\vec{0},\varsigma_*\vec{1})$. 
	\end{lemma}
	\proof{}
		This follows from Lemma \ref{LamGamLemma}, claim (i), the inequality $\kappa_2(\vec{AB}) \leq \kappa_2(\vec{A})\kappa_2(\vec{B})$ valid for any matrices $\vec{A}$ and $\vec{B}$, and the formula
		\begin{equation*}
			\kappa_2(\bsym{\Lambda}(\vec{p})\inv ) 
				= \frac{ \min_j \abs{ \lambda_j(\vec{p}) }  }
						{ \max_j \abs{ \lambda_j(\vec{p}) } }.
		\end{equation*}
	\endproof
	Lemma \ref{EtaCondition} states that the greater the variation in aggregate absolute rate of change in inclusive values, the more poorly conditioned $(\tilde{D}\vec{P})(\vec{p})^\top$ is relative to $\vec{I} - \tilde{\bsym{\Omega}}(\vec{p})$. Because $\bsym{\Lambda}(\vec{p})$ is diagonal, $\norm{ \bsym{\Lambda}(\vec{p}) }_1 = \norm{ \bsym{\Lambda}(\vec{p}) }_2 = \norm{ \bsym{\Lambda}(\vec{p}) }_\infty$ and thus the same bound applies for condition numbers in norms other than the 2-norm.
	
	Third, Eqn. (\ref{EtaEqn}) also provides bounds on the magnitude of values taken by $\bsym{\eta}$:
	\begin{lemma} 
		Suppose Assumptions \ref{MixedLogitUtilityAssump}-\ref{LeibnizRuleCondition} hold. For all $\vec{p} \in (0,\varsigma_*)^J$, $\bsym{\eta}$ satisfies
		\begin{equation}
			\label{EtaBounds}
			\frac{ \norm{ \bsym{\Lambda}(\vec{p})\inv\vec{P}(\vec{p}) }_\infty }
							{ 1 + \norm{ \tilde{\bsym{\Omega}}(\vec{p}) }_\infty }
					\leq \norm{ \bsym{\eta}(\vec{p}) }_\infty
					\leq \frac{\norm{\bsym{\Lambda}(\vec{p})\inv\vec{P}(\vec{p})}_\infty}
							{1 - \norm{\tilde{\bsym{\Omega}}(\vec{p})}_\infty}. 
		\end{equation}
	\end{lemma}
	\proof{}
		This follows immediately from Eqn. (\ref{EtaEqn}), using the triangle inequality.
	\endproof
	
	The upper bound suggests the following assumptions to ensure that $\bsym{\eta}$ itself is bounded:
	\begin{assumption}
		\label{BoundedAssumption}
		Suppose there exist $M \in (0,\infty)$ and $\varepsilon \in (0,1)$ such that
		\begin{align}
			\label{LPBound}
			\sup \big\{ \norm{ \bsym{\Lambda}(\vec{p})\inv\vec{P}(\vec{p}) }_\infty \; : \; \vec{p} \in (0,\varsigma_*)^J \big\} &= M < \infty \\
			\label{LGBound}
			\sup \big\{ \norm{\tilde{\bsym{\Omega}}(\vec{p})}_\infty \; : \; \vec{p} \in (0,\varsigma_*)^J \big\}
				&= 1 - \varepsilon < 1. 
		\end{align}
	\end{assumption}
	
	Under simple Logit, $P_k^L(\vec{p})/\abs{\lambda_k(\vec{p})} = \abs{(Dw_k)(p_k)}\inv$ and $\bsym{\Omega}_f(\vec{p}) = \vec{1}\vec{P}_f^L(\vec{p})^\top$. Thus Eqn. (\ref{LPBound}) is akin to concavity of $w_k$ for all sufficiently large $p_k$, and Eqn. (\ref{LGBound}) is implied by $\vartheta > -\infty$, i.e. the existence of an outside good with positive purchase probability. 
	
	\begin{lemma}
		\label{BoundedEta}
		Suppose Assumptions \ref{MixedLogitUtilityAssump}-\ref{LeibnizRuleCondition} hold.
		\begin{itemize}
			\item[(i)] If Assumption \ref{BoundedAssumption} also holds, $N = \sup \{ \norm{ \bsym{\eta}(\vec{p}) }_\infty : \vec{p} \in (0,\varsigma_*)^J \} < \infty$. 
			\item[(ii)] Moreover Eqn. (\ref{LPBound}) in Assumption \ref{BoundedAssumption} is necessary for $N$ to be finite. 
		\end{itemize}
	\end{lemma}
	
	Unfortunately some simple models do not satisfy Assumption \ref{BoundedAssumption}. A simple Logit model with $w(p) = - \alpha \log p$ for some $\alpha > 0$ violates Eqn. (\ref{LPBound}). More generally, the \cite{Boyd80} model of Example \ref{EX:BM80Example} does not satisfy Eqn. (\ref{LPBound}). This is most easily seen by noting that finite-sample approximations to this model have 
	\begin{equation*}
		\lim_{p_k \to \infty} \norm{ \bsym{\Lambda}(\vec{p})\inv\vec{P}(\vec{p}) }_\infty
			= \max_{s=1,\dotsc,S} \left\{ \frac{1}{ \alpha_s } \right\} 
	\end{equation*}
	where $\{ \alpha_s \}_{s=1}^S$ are the sampled price coefficients. Of course, as $S \to \infty$, $\min_{s=1,\dotsc,S} \{ \alpha_s \} \to 0$, and thus $\norm{ \bsym{\Lambda}(\vec{p})\inv\vec{P}(\vec{p}) }_\infty \to \infty$.

	\subsection{The $\zeta$-Markup function}
	\label{SUBSEC:ZetaMap}
	
	Substituting Eqn. (\ref{Decomposition}) into Eqn. (\ref{SSCs}) yields the $\bsym{\zeta}$-{\em markup equation} introduced in \cite{Morrow10}:
	\begin{equation}
		\label{ZetaMap}
		\vec{p} = \vec{c} + \bsym{\zeta}(\vec{p})
		\quad\text{where}\quad
		\bsym{\zeta}(\vec{p}) 
			= \bsym{\Lambda}(\vec{p})\inv\tilde{\bsym{\Gamma}}(\vec{p})^\top(\vec{p} - \vec{c}) 
				- \bsym{\Lambda}(\vec{p})\inv\vec{P}(\vec{p})
	\end{equation}
	when $\bsym{\Lambda}(\vec{p})$ is nonsingular (Section \ref{SUBSECChoiceProbabilityDerivatives}, \ref{SUBSEC:ZetaMap}). Thus the $\bsym{\zeta}$-markup equation is specific to Mixed Logit models, unlike the BLP-markup equation. 
	
	We observe a relationship between the maps $\bsym{\zeta}$ and $\bsym{\eta}$. 
	\begin{proposition}
		\label{ZetaEtaDistinct}
		Suppose Assumption \ref{MixedLogitUtilityAssump}-\ref{LeibnizRuleCondition} hold. For any $\vec{p} \in (0,\varsigma_*)^J$, $\bsym{\zeta}(\vec{p}) = \tilde{\bsym{\Omega}}(\vec{p})(\vec{p} - \vec{c}) + ( \vec{I} - \tilde{\bsym{\Omega}}(\vec{p}) ) \bsym{\eta}(\vec{p})$. 
	\end{proposition}
	\proof{}
		This follows directly from Eqns. (\ref{EtaEqn}) and (\ref{ZetaMap}).
	\endproof
	In so far as $\bsym{\eta}$ and $\bsym{\zeta}$ are distinct maps, they can generate numerical methods for the computation of equilibrium prices with entirely different properties. The equation above implies that $\bsym{\zeta}(\vec{p}) = \bsym{\eta}(\vec{p})$ if, and only if, $\vec{p} - \vec{c} - \bsym{\eta}(\vec{p}) = \vec{p} - \vec{c} - \bsym{\zeta}(\vec{p})$ lies in the null space of $\tilde{\bsym{\Omega}}(\vec{p})$. Thus if $\tilde{\bsym{\Omega}}(\vec{p})$ is full-rank, $\bsym{\zeta}$ and $\bsym{\eta}$ coincide {\em only} at simultaneously stationary prices. Simple examples of Mixed Logit models can be constructed that always have $\mathrm{rank}(\tilde{\bsym{\Omega}}(\vec{p})) = J$. For Logit, $\bsym{\Omega}_f(\vec{p}) = \vec{1}\vec{P}_f^L(\vec{p})^\top$ for all $f$ and $\tilde{\bsym{\Omega}}(\vec{p})$ always has rank $F \leq J$. However the analysis in \cite{Morrow08a} can be used to show that $\bsym{\zeta}$ and $\bsym{\eta}$ coincide only at simultaneously stationary prices. 
	
	We now explore $\bsym{\zeta}$'s asymptotic properties. 
	\begin{lemma}
		\label{ZetaGrowth}
		Under Assumption \ref{BoundedAssumption} $\norm{\bsym{\zeta}(\vec{p})}_\infty < \norm{\vec{p}-\vec{c}}_\infty$ whenever $\norm{\vec{p}-\vec{c}}_\infty > M / \varepsilon$. Moreover $\norm{\vec{p}-\vec{c}}_\infty - \norm{\bsym{\zeta}(\vec{p})}_\infty \to \infty$ as $\norm{ \vec{p} - \vec{c} }_\infty \to \infty$. 
	\end{lemma}

	\proof{}
		We simply note that
		\begin{align*}
			\norm{\bsym{\zeta}(\vec{p})}_\infty 
				&\leq \norm{ \tilde{\bsym{\Omega}}(\vec{p}) }_\infty \norm{ \vec{p} - \vec{c} }_\infty
						+ \norm{ \bsym{\Lambda}(\vec{p})\inv\vec{P}(\vec{p}) }_\infty \\
				&\leq ( 1 - \varepsilon ) \norm{ \vec{p} - \vec{c} }_\infty + M \\
				&\leq \norm{ \vec{p} - \vec{c} }_\infty - \big( \varepsilon \norm{ \vec{p} - \vec{c} }_\infty - M \big) \\
				&\leq \left[ 1 - \varepsilon + \frac{M}{\norm{ \vec{p} - \vec{c} }_\infty} \right] \norm{ \vec{p} - \vec{c} }_\infty
		\end{align*}
		Now if $\norm{\vec{p}-\vec{c}}_\infty > M / \varepsilon$, then $M / \norm{\vec{p}-\vec{c}}_\infty < \varepsilon$. Thus
		\begin{align*}
			\norm{\bsym{\zeta}(\vec{p})}_\infty 
				< \left[ 1 - \varepsilon + \varepsilon \right] \norm{ \vec{p} - \vec{c} }_\infty
				= \norm{ \vec{p} - \vec{c} }_\infty. 
		\end{align*}
		
		
		To prove that $\norm{ \vec{p} - \vec{c} }_\infty - \norm{\bsym{\zeta}(\vec{p})}_\infty \to \infty$, note that
		\begin{align*}
			\norm{ \vec{p} - \vec{c} }_\infty - \norm{\bsym{\zeta}(\vec{p})}_\infty 
				\geq \left( \varepsilon - \frac{M}{\norm{ \vec{p} - \vec{c} }_\infty} \right) \norm{ \vec{p} - \vec{c} }_\infty
		\end{align*}
		For all $\norm{ \vec{p} - \vec{c} }_\infty > M/\varepsilon$, the term in parentheses is positive. Furthermore, this term approaches $\varepsilon$ as $\norm{ \vec{p} - \vec{c} }_\infty \to \infty$. Thus $\norm{ \vec{p} - \vec{c} }_\infty - \norm{\bsym{\zeta}(\vec{p})}_\infty \to \infty$ as $\norm{ \vec{p} - \vec{c} }_\infty \to \infty$.
	\endproof
	
	A slightly different assumption concerning $\tilde{\bsym{\Omega}}(\vec{p})$ is useful when analyzing the $\bsym{\zeta}$ map. 
	\begin{assumption}
		\label{BoundedOmegaProductsAssumption}
		Suppose
		\begin{equation}
			\label{LGBound2}
			\sup \Big\{ \norm{ \tilde{\bsym{\Omega}}(\vec{p})(\vec{p} - \vec{c}) }_\infty : \vec{p} \in (0,\varsigma_*)^J \Big\} < \infty, 
		\end{equation}
	\end{assumption}
	
	\begin{lemma}
		Suppose Assumption \ref{MixedLogitUtilityAssump}-\ref{LeibnizRuleCondition} holds and $\varsigma_* = \infty$. Then $\bsym{\zeta}$ is bounded if, and only if, Eqn. (\ref{LPBound}) and Eqn. (\ref{LGBound2}) hold. 
	\end{lemma}
	\proof{}
		This follows directly from the triangle inequality and the non-negativity of $\tilde{\bsym{\Omega}}(\vec{p})(\vec{p} - \vec{c})$ and $- \bsym{\Lambda}(\vec{p})\inv\vec{P}(\vec{p})$ for all $\vec{p} \geq \vec{c}$.
	\endproof
	
	For future reference, we prove that Eqn. (\ref{LGBound2}) {\em strengthens} Eqn. (\ref{LGBound}). 
	\begin{lemma}
		\label{LGBound2ImpliesLGBound}
		If Eqn. (\ref{LGBound2}) holds, then Eqn. (\ref{LGBound}) holds. 
	\end{lemma}
	\proof{}
		Note that Eqn. (\ref{LGBound2}) implies that for any $k$,
		\begin{equation*}
			\lim_{p_j \to \infty} \big( \omega_{k,j}(\vec{p})(p_j-c_j) \big) < \infty
			\quad\text{for all}\quad
			j \in \set{J}_{f(k)}. 
		\end{equation*}
		This, in turn, implies that $\omega_{k,j}(\vec{p}) \to 0$ as $p_j \to \infty$. 
		
		Now Eqn. (\ref{LGBound}) {\em fails} only if $\lim_{\vec{p} \to \vec{q}} \norm{ \tilde{\bsym{\Omega}}(\vec{p}) }_\infty = 1$ where $\vec{q}$ has some $q_j = \infty$. But the row sums of $\tilde{\bsym{\Omega}}(\vec{p})$ satisfy
		\begin{equation*}
			\lim_{\vec{p} \to \vec{q}} \left[ \sum_{j \in \set{J}_{f(k)}} \omega_{k,j}(\vec{p}) \right]
				= \sum_{j \in \set{J}_{f(k)} \cup \{ j : q_j < \infty \} } \omega_{k,j}(\vec{q})
				< 1. 
		\end{equation*}
		Thus if Eqn. (\ref{LGBound2}) holds, Eqn. (\ref{LGBound}) cannot fail.
	\endproof

%
	
	
	\subsection{Existence of Simultaneously Stationary Prices}
	\label{SUBSECExistence}
	
	This section provides two existence results. Neither establish the existence of a local {\em equilibrium}, or the {\em uniqueness} of simultaneously stationary points. To address the existence of local equilibrium will require additional conditions to ensure that each firm's profits are locally concave at the simultaneously stationary prices whose existence can be ensured \citep{Morrow08}. Little is known about how to address the uniqueness of simultaneously stationary points. Indeed, \cite{Morrow10} provide an example of a Mixed Logit model with 9 simultaneously stationary prices, 4 of which are local equilibria and 2 of which are proper equilibria. 
	
	Assumption \ref{BoundedAssumption} ensures the existence of finite simultaneously stationary prices when $\varsigma_* = \infty$. 
	\begin{corollary}
		\label{Existence}
		Suppose $\varsigma_* = \infty$ and Assumptions \ref{MixedLogitUtilityAssump}-\ref{LeibnizRuleCondition}, \ref{BoundedAssumption} hold. Then there exists at least one vector of finite simultaneously stationary prices. 
	\end{corollary}
	\proof{}
		This is a direct consequence of Brouwer's fixed-point theorem. $\vec{c} + \bsym{\eta}(\cdot)$ is a continuous map that takes the compact, convex set $[0,M/\varepsilon]^J$ into itself, and thus there is at least one fixed-point $\vec{p} = \vec{c} + \bsym{\eta}(\vec{p}) \in [0,M/\varepsilon]^J$.
	\endproof
	
	To apply Corollary \ref{Existence} to cases when $\varsigma_* < \infty$, $\bsym{\eta}$ must be extended from $(0,\varsigma_*)$ to all of $\set{P}^J$ preserving the bounds (\ref{EtaBounds}). This is easy for many of the simulation-based approximations encountered in practice, but difficult for the general case. 
	
	We can extend this existence result using Eqn. (\ref{LPBound2}) and the $\bsym{\zeta}$ map. 
	\begin{lemma}
		\label{EventuallyPointsOutward}
		Suppose $\varsigma_* = \infty$, Assumptions \ref{MixedLogitUtilityAssump}-\ref{LeibnizRuleCondition}, Eqn. (\ref{LGBound2}), and Eqn. (\ref{LPBound2}) hold. Then there exists some $\bar{q}_k > c_k$ such that $p_k - c_k - \zeta_k(\vec{p}) > 0$ for all $\vec{p}$ with $p_k \geq \bar{q}_k$. 
	\end{lemma}
	\proof{}
		The assumed bound implies
		\begin{equation*}
			\abs{ \sum_{j\in\set{J}_f} \omega_{k,j}(\vec{p}) (p_j-c_j) } \leq L < \infty
		\end{equation*}
		for any $k$ and any $\vec{p} \in (0,\infty)^J$. Consider
		\begin{align*}
			p_k - c_k - \zeta_k(\vec{p})
				&= ( p_k - c_k ) 
					- \sum_{j\in\set{J}_f} \omega_{k,j}(\vec{p}) (p_j-c_j)
					+ \frac{P_k(\vec{p})}{\lambda_k(\vec{p})} \\
				&= \left( 1 - \frac{P_k(\vec{p})}{\abs{\lambda_k(\vec{p})}(p_k-c_k)} \right) 
						( p_k - c_k ) 
					- \sum_{j\in\set{J}_f} \omega_{k,j}(\vec{p}) (p_j-c_j)
		\end{align*}
		If Eqn. (\ref{LPBound2}) holds, then
		\begin{equation*}
			0 \leq \lim_{p_k \to \infty} \left[ \frac{P_k(\vec{p})}{\abs{\lambda_k(\vec{p})}(p_k-c_k)} \right] 
				\leq \delta < 1. 
		\end{equation*}
		Thus for any $\epsilon > 0$, there exists some $\bar{p}_k > 0$ and $\triangle(\vec{p})$ with $\abs{ \triangle(\vec{p}) } < \epsilon$ such that
		\begin{equation*}
			\frac{P_k(\vec{p})}{\lambda_k(\vec{p})(p_k-c_k)} \leq \delta + \triangle(\vec{p})
			\quad\quad\text{for all}\quad\quad
			p_k > \bar{p}_k. 
		\end{equation*}
		Thus
		\begin{align*}
			p_k - c_k - \zeta_k(\vec{p})
				\geq ( 1 - \delta + \triangle(\vec{p}) ) ( p_k - c_k ) - L
			\quad\quad\text{for all}\quad\quad
			p_k > \bar{p}_k. 
		\end{align*}
		In particular, if we choose $\epsilon \leq (1-\delta)/2$ we have 
		\begin{align*}
			p_k - c_k - \zeta_k(\vec{p})
				\geq \left( \frac{ 1 - \delta }{ 2 } \right) ( p_k - c_k ) - L
				= \left( \frac{ 1 - \delta }{ 2 } \right) 
					\left( p_k - c_k - \frac{ 2L }{ 1 - \delta } \right)
				> 0
		\end{align*}
		for all $p_k \geq \bar{q}_k = \max\{ c_k + 2L / (1-\delta) , \bar{p}_k \}$.
	\endproof
	
	One consequence of this lemma is that infinite prices are never an equilibrium. 
	\begin{corollary}
		Under the assumptions of Lemma \ref{EventuallyPointsOutward}, any profit derivative is eventually negative. 
	\end{corollary}
	\proof{}
		Note that
		\begin{equation*}
			(D_k\hat{\pi}_{f(k)})(\vec{p})
				= - \abs{ \lambda_k(\vec{p}) }
					( p_k - c_k - \zeta_k(\vec{p}) ). 
		\end{equation*}
		Since $p_k - c_k - \zeta_k(\vec{p})$ is positive for all large enough $p_k$, $(D_k\hat{\pi}_{f(k)})(\vec{p})$ is negative for all large enough $p_k$, regardless of $\vec{p}_{-k}$.
	\endproof
	
	Another consequence of Lemma \ref{EventuallyPointsOutward} is an alternative existence result. 
	\begin{corollary}
		Under the assumptions of Lemma \ref{EventuallyPointsOutward} there exists at least one simultaneously stationary point. 
	\end{corollary}
	\proof{}
		Following \cite{Morrow08a}, we prove this proposition using the Poincare-Hopf Theorem \citep{Milnor65}. The logic is simple: We will consider the vector field $\vec{p} - \vec{c} - \bsym{\zeta}(\vec{p})$ on a hyper-rectangle $[\vec{c},\bar{\vec{q}}]$ whose critical points are simultaneously stationary; $\bar{\vec{q}}$ has components $\bar{q}_k$ defined in Lemma \ref{EventuallyPointsOutward}. The Poincare-Hopf Theorem then states that the sum of the indices of all the critical points of this vector field equals one, the Euler characteristic of $[\vec{c},\bar{\vec{q}}]$. Thus it is not possible that the vector field have {\em no} critical points, for then the sum of indices would be zero. 
		
		We must only prove one property of $\vec{p} - \vec{c} - \bsym{\zeta}(\vec{p})$: that this vector field points outward on the boundary of the chosen hyper-rectangle. Half of this proof is Lemma \ref{EventuallyPointsOutward}, in which we prove that $p_k - c_k - \zeta_k(\vec{p}) > 0$ if $\vec{p} \in [\vec{c},\bar{\vec{q}}]$ with $p_k = \bar{q}_k$. We must also show that $p_k - c_k - \zeta_k(\vec{p}) < 0$ if $\vec{p} \in [\vec{c},\bar{\vec{q}}]$ with $p_k = c_k$. But 
		\begin{align*}
			c_k - c_k - \zeta_k(\vec{p})
				= - \left( \sum_{j \in \set{J}_f} \omega_{k,j}(\vec{p})(p_j - c_j) 
						+ \frac{P_k(\vec{p})}{\abs{\lambda_k(\vec{p})}} \right)
				< 0.
		\end{align*}
	\endproof
	
	This proof does not need to make any claims about the number of critical points, or of their indices. If it can be shown that any critical point of $\vec{p} - \vec{c} - \bsym{\zeta}(\vec{p})$ cannot have a zero or negative index, then the simultaneously stationary point is unique.

%% file: extendednummethods.tex

\section{Computational Methods}
\label{ECSEC:ExtendedMethods}
	
	This section provides details for the four approaches to computing equilibrium prices described in \cite{Morrow10}; see Table \ref{TAB:MethodSummary}. Section \ref{SUBSEC:NewtonsMethod} briefly reviews Newton's method, followed by application of Newton's method to solve Eqn. (\ref{SSCs}) in Section \ref{SUBSEC:CG-NM}. Newton's method applied directly to Eqn. (\ref{SSCs}) may compute ``spurious'' solutions with infinite prices because the combined gradient vanishes as prices increase without bound. Section \ref{SUBSEC:FP-NM} avoids this difficulty by applying Newton's method to the two markup equations instead of Eqn. (\ref{SSCs}) itself. Section \ref{SUBSEC:FPIs} discusses fixed-point iterations based on the BLP- and $\bsym{\zeta}$-markup equations, and Section \ref{SUBSEC:PreliminaryConsiderations} reviews a number of practical considerations. 
	
	\begin{table}
		\caption{Summary of the numerical methods examined in this article.}
		\label{TAB:MethodSummary}
		\begin{tabular}{clcll}
			\\
			\multicolumn{5}{l}{{\bfseries Newton Methods} (NM)} \\ \hline
			\\ Abbr. & Method & Section & Advantage & Our Experience$^{\text{(a)}}$ \\ \hline
			CG-NM 
				& Solve $\vec{F}_\pi(\vec{p}) = (\tilde{\nabla}\hat{\pi})(\vec{p}) = \vec{0}$ 
				& \ref{SUBSEC:CG-NM}
				& $-$
				& Unreliable, slow \\
			$\bsym{\eta}$-NM
				& Solve $\vec{F}_\eta(\vec{p}) = \vec{p} - \vec{c} - \bsym{\eta}(\vec{p}) = \vec{0}$ 
				& \ref{SUBSEC:FP-NM}
				& Coercive
				& Reliable, slow \\
			$\bsym{\zeta}$-NM
				& Solve $\vec{F}_\zeta(\vec{p}) = \vec{p} - \vec{c} - \bsym{\zeta}(\vec{p}) = \vec{0}$
				& \ref{SUBSEC:FP-NM}
				& Coercive
				& Reliable, slow \\ 
			\\
			\multicolumn{5}{l}{{\bfseries Fixed-Point Iterations} (FPI)} \\ \hline
			\\ Abbr. & Method & Section & Advantage & Our Experience \\ \hline
			$\bsym{\zeta}$-FPI
				& Iterate $\vec{p} \leftarrow \vec{c} + \bsym{\zeta}(\vec{p})$ 
				& \ref{SUBSEC:FPIs}
				& Easy to evaluate
				& Reliable, fast \\
			$\bsym{\eta}$-FPI
				& Iterate $\vec{p} \leftarrow \vec{c} + \bsym{\eta}(\vec{p})$ 
				& \ref{SUBSEC:FPIs}
				& $-$
				& Not convergent \\ 
			\\ \hline
			\multicolumn{5}{l}{(a) Conclusions on behavior of these methods is based on the numerical experiments} \\
			\multicolumn{5}{l}{described in \cite{Morrow10}, using a novel \texttt{GMRES}-Newton method with} \\
			\multicolumn{5}{l}{Levenberg-Marquardt style trust-region global convergence strategy.}
		\end{tabular}
	\end{table}

	
	\subsection{Newton's Method}
	\label{SUBSEC:NewtonsMethod}
	
	Newton's method, a classical technique to compute a zero of an arbitrary function $\vec{F} : \R^J \to \R^J$, is now a portfolio of related approaches to solve nonlinear systems \citep{Ortega70, Kelley95, Dennis96, Judd98, Kelley03}. Generally speaking, Newton-type methods are differentiated in two relatively independent directions: (i) the technique used to approximate the Jacobian matrices $(D\vec{F})$ and solve for the Newton step and (ii) the technique used to enforce convergence from arbitrary initial conditions. See \cite{Dennis96}, \cite{Judd98}, or \cite{Kelley03} for good treatments of these issues. Choosing the right variant of Newton's method determines the reliability and efficiency of equilibrium price computations. 
		
	Problem formulation also determines the reliability and efficiency of equilibrium price computations using Newton's method. Scalings of the variables and function values are one prominent example of a problem transformation that improves the performance of Newton's method \citep{Dennis96}. Nonlinear problem preconditioning can also be important \citep{Cai02}, as the following example demonstrates. 
	
	\begin{example}
		\label{EX:SimpleExample}
		Let $\vec{F} : \R^N \to \R^N$ be defined by $\vec{F}(\vec{x}) = \vec{x} / ( 1 + \norm{ \vec{x} }_2^2 )$. Iterating Newton steps converges to the unique (finite) zero $\vec{x}_* = \vec{0}$ only from initial conditions $\vec{x}_0$ with $\norm{ \vec{x}_0 }_2 < 1/\sqrt{3}$. Newton's method diverges or fails for all other starting points. Standard global convergence strategies for Newton's method (line search, trust region methods) cannot improve this poor global convergence behavior because $\norm{\vec{F}(\vec{x})}_2$ has unbounded level sets; see \cite{Morrow10} for details. 
		
		A simple nonlinear transformation overcomes this poor global convergence behavior. Note that $\vec{F}(\vec{x}) = \vec{A}(\vec{x}) \vec{f}(\vec{x})$ where $\vec{A}(\vec{x}) = (1 + \norm{\vec{x}}_2^2)\inv\vec{I}$ and $\vec{f}(\vec{x}) = \vec{x}$. Because $\vec{A}(\vec{x})$ is nonsingular for all $\vec{x}$, the problems $\vec{F}(\vec{x}) = \vec{0}$ and $\vec{f}(\vec{x}) = \vec{0}$ have identical solution sets. However applying Newton's method to the problem $\vec{f}(\vec{x}) = \vec{0}$ trivially converges in a single step from any initial condition without a global convergence strategy.
	\end{example}
	
	Example \ref{EX:SimpleExample} illustrates why computing equilibrium prices based on the markup equations is more reliable and efficient than using Eqn. (\ref{SSCs}) directly. The following two sections echo the pattern of this example to provide the details. 
	
	
	\subsection{Newton's Method on the Combined Gradient}
	\label{SUBSEC:CG-NM}
	
	The most direct approach to compute equilibrium prices using Newton's method is to solve $\vec{F}_\pi(\vec{p}) = (\tilde{\nabla}\hat{\pi})(\vec{p}) = \vec{0}$, abbreviated CG-NM in Table \ref{TAB:MethodSummary}. This approach works well when the initial condition is near an equilibrium, as required by theory \citep{Ortega70,Kelley95,Dennis96}. In practice, computing counterfactual equilibrium prices starting with the observed prices may exploit this local convergence if changes to exogeneous variables have a relatively small impact on equilibrium prices. On the other hand, CG-NM can be unreliable when started ``far'' from equilibrium. 
	
	The challenge is the tendency for the derivatives of profits to vanish as prices become large \cite{Morrow10}, as demonstrated in Example \ref{EX:LogitProfitsVanish} below. 
	
	\begin{example}
		\label{EX:LogitProfitsVanish}
		Consider a simple Logit model with linear in price utility and an outside good: $u(p) = - \alpha p + v$ for some $\alpha > 0$ and any $v \in \R$, and $\vartheta > -\infty$. The derivative of firm $f$'s profit function with respect to the price of product $k \in \set{J}_f$ is
		\begin{equation*}
			(D_k\hat{\pi}_f)(\vec{p})
				= - \alpha P_k^L(\vec{p}) (p_k-c_k)
					+ \alpha P_k^L(\vec{p}) \hat{\pi}_f(\vec{p})
					+ P_k^L(\vec{p}). 
		\end{equation*}
		Since $P_k^L(\vec{p})$ and $P_k^L(\vec{p})(p_k-c_k)$ both vanish as $p_k \to \infty$ (as is easily checked), $\hat{\pi}_f(\vec{p})$ is bounded in $\vec{p}$. Thus $(D_k\hat{\pi}_f)(\vec{p}) \to 0$ as $p_k \to \infty$. 
	\end{example}
	
	We now provide a general assumption under which $(D_k\hat{\pi}_{f(k)})(\vec{p}) \to 0$ as $p_k \uparrow \varsigma_*$. 
	\begin{assumption}
		\label{LVA}
		Let $\psi_k$ be defined as in Assumption \ref{LeibnizRuleCondition}. Assume: (i) $p_k \psi_k(\bsym{\theta},p_k) \to 0$ as $p_k \uparrow \varsigma(\bsym{\theta})$ for $\mu$-a.e. $\bsym{\theta}$. (ii) There exists $M < \infty$ and $\bar{p}_k \in [0,\varsigma_*)$ such that $p_k \int \psi_k(\bsym{\theta},p_k) d\mu(\bsym{\theta}) \leq M$ for all $p_k \in (\bar{p}_k,\varsigma_*)$. 
	\end{assumption}
	As with Assumption \ref{LeibnizRuleCondition} above, (i) and (ii) are essentially conditions for the Dominated Convergence Theorem. 
	
	Assumption \ref{LVA} (i) extends Assumption \ref{LeibnizRuleCondition} (i) to include a neighborhood of $\varsigma_*$. Note that if $\varsigma(\bsym{\theta}) < \infty$ then (i) holds if, and only if, $\psi_k(\bsym{\theta},p_k) \to 0$ as $p_k \uparrow \varsigma(\bsym{\theta})$; i.e. $\psi_k(\bsym{\theta},\cdot)$ is continuous at $\varsigma(\bsym{\theta})$. Thus if $\varsigma(\bsym{\theta}) < \infty$ Assumption \ref{LVA} (i) and Assumption \ref{LeibnizRuleCondition} (i) are the same. If $\varsigma(\bsym{\theta}) = \infty$ and $p_k\psi_k(\bsym{\theta},p_k) \to 0$ as $p_k \uparrow \infty$, then necessarily $\psi_k(\bsym{\theta},p_k) \to 0$ as $p_k \uparrow \infty$. The converse, however, need not hold. 
	
	If $\varsigma_* < \infty$, Assumption \ref{LVA} (ii) simply says that $\int \psi_k(\bsym{\theta},p_k) d\mu(\bsym{\theta})$ is bounded as $p_k \uparrow \varsigma_*$. This is not implied by Assumption \ref{LeibnizRuleCondition} (ii), but is a natural extension of it. 
	
	\begin{lemma}
		\label{LambdaVanishes}
		Suppose Assumptions \ref{MixedLogitUtilityAssump}-\ref{LVA} hold. Then $p_k \abs{ \lambda_k(\vec{p}) } \to 0$ as $p_k \uparrow \varsigma_*$ for all $k$. Additionally, $p_k \abs{ \gamma_{j,k}(\vec{p}) } \to 0$ as $p_k \uparrow \varsigma_*$ for all $j$. Subsequently, $(D_k\hat{\pi}_{f(k)})(\vec{p}) \to 0$ as $p_k \uparrow \varsigma_*$. 
	\end{lemma}
	\proof{}
		Let $\{ p_k^{(n)} \}_{n \in \N} \subset (0,\varsigma_*)$ be any sequence converging to $\varsigma_*$. Define $\Psi_k^{(n)}: \set{T} \to \set{P}$ by $\Psi_k^{(n)}(\bsym{\theta}) = p_k^{(n)} \psi_k(\bsym{\theta},p_k^{(n)})$. The functions $\{ \Psi_k^{(n)} \}_{n \in \N}$ converge pointwise to zero and have integrals uniformly bounded by the constant $M$. By the Dominated Convergence Theorem
		\begin{equation*}
			\lim_{n \to \infty} \int \Psi_k^{(n)}(\bsym{\theta}) d \mu(\bsym{\theta})
				= \int \Big( \lim_{n \to \infty} \Psi_k^{(n)}(\bsym{\theta}) \Big) d \mu(\bsym{\theta})
				= 0.
		\end{equation*}
	\endproof
	
	In other words, under Assumption \ref{LVA} the components of $\vec{F}_\pi$ vanish as the corresponding price tends to $\varsigma_*$ even though this may not mean that $\varsigma_*$ maximizes profits. Because of this, CG-NM may converge to a zero of $\vec{F}_\pi$ at $\varsigma_*\vec{1}$, or with some components equal to $\varsigma_*$, that is {\em not} an equilibrium. 
	
	Note that even though the price derivatives vanish at infinity, this does not mean that infinite prices maximize profits. Nonetheless, CG-NM may converge to a zero of $\vec{F}_\pi$ with some components equal to infinity that is {\em not} an equilibrium. Moreover, because the components of $\vec{F}_\pi(\vec{p})$ can vanish over some divergent sequences, standard global convergence strategies based on minimizing $\norm{\vec{F}_\pi(\vec{p})}_2$ will not be effective ways of avoiding this behavior. As in Example \ref{EX:SimpleExample}, we must reformulate the problem to obtain reliable and efficient approaches for computing equilibrium prices.

	
	\subsection{Newton's Method and the Markup Equations}
	\label{SUBSEC:FP-NM}
	
	Reliable and efficient implementations of Newton's method are found by observing that the combined gradient, $\vec{F}_\pi$, can be written as follows:
	\begin{align}
		\label{Eta-CGRelation}
		\vec{F}_\pi(\vec{p})
			&= (\tilde{D}\vec{P})(\vec{p})^{\top}\vec{F}_\eta(\vec{p})
			&\quad\text{where}\quad&&
			\vec{F}_\eta(\vec{p})
				&= \vec{p} - \vec{c} - \bsym{\eta}(\vec{p}) \\
		\label{Zeta-CGRelation}
		\vec{F}_\pi(\vec{p})
			&= \bsym{\Lambda}(\vec{p})\vec{F}_\zeta(\vec{p})
			&\quad\text{where}\quad&&
			\vec{F}_\zeta(\vec{p})
				&= \vec{p} - \vec{c} - \bsym{\zeta}(\vec{p}). 
	\end{align}
	Either $\vec{F}_\eta$ or $\vec{F}_\zeta$ can be used to compute simultaneously stationary prices when $(\tilde{D}\vec{P})(\vec{p})^{\top}$ and $\bsym{\Lambda}(\vec{p})$, respectively, are nonsingular \citep{Morrow10}. Of course, $\vec{F}_\eta$ and $\vec{F}_\zeta$ recast the first-order condition as a fixed-point problem: $\vec{F}_\eta$ is zero if and only if the BLP-markup equation holds, and $\vec{F}_\zeta$ is zero if and only if the $\bsym{\zeta}$-markup equation holds.  
	
	Solving $\vec{F}_\eta(\vec{p}) = \vec{0}$ or $\vec{F}_\zeta(\vec{p}) = \vec{0}$, abbreviated $\bsym{\eta}$-NM and $\bsym{\zeta}$-NM respectively in Table \ref{TAB:MethodSummary}, requires the solution of nontrivial nonlinear systems with Newton's method. $\bsym{\eta}$-NM and $\bsym{\zeta}$-NM, however, are less likely to have the computational problems that CG-NM exhibits because they exploit {\em norm-coercivity} of the maps $\vec{F}_\eta$ and $\vec{F}_\zeta$ \citep{Morrow10}. A norm-coercive map has a norm that tends to infinity with the norm of its argument \citep{Ortega70, Harker90}. Globally convergent implementations of Newton's method that decrease the value of $\norm{ \vec{F}(\vec{p}) }_2$ in each step produce bounded sequences of iterates when $\vec{F}$ is norm-coercive. Thus, solving the BLP- or $\bsym{\zeta}$-markup equation instead of the literal first-order condition removes the tendency for applications of Newton's method to compute ``spurious'' solutions at infinity. 
	
	We now prove that the maps $\vec{F}_\eta$ and $\vec{F}_\zeta$ are indeed coercive. 
	
	\begin{lemma}
		\label{EtaCoercivityImpliesZetaCoercivity}
		Suppose $\varsigma_* = \infty$ and Assumption \ref{MixedLogitUtilityAssump}-\ref{LeibnizRuleCondition} hold. 
		\begin{itemize}
			\item[(i)] Norm-coercivity of $\vec{F}_\zeta(\vec{p})$ implies that of $\vec{F}_\eta(\vec{p})$. 
			\item[(ii)] If Eqn. (\ref{LGBound}) holds, then norm-coercivity of $\vec{F}_\eta(\vec{p})$ implies that of $\vec{F}_\zeta(\vec{p})$. 
		\end{itemize}
	\end{lemma}
	\proof{}
		Proposition \ref{ZetaEtaDistinct} implies that
		\begin{equation*}
			\vec{p} - \vec{c} - \bsym{\zeta}(\vec{p})
				= \big( \vec{I} - \tilde{\bsym{\Omega}}(\vec{p}) \big) 
					( \vec{p} - \vec{c} - \bsym{\eta}(\vec{p}) ). 
		\end{equation*}
		To prove (i), note that
		\begin{equation*}
			\norm{ \vec{p} - \vec{c} - \bsym{\eta}(\vec{p}) }_\infty
				\geq \left( \frac{1}{ 1 + \norm{ \tilde{\bsym{\Omega}}(\vec{p}) }_\infty } \right)
						\norm{ \vec{p} - \vec{c} - \bsym{\zeta}(\vec{p}) }_\infty
				\geq \left( \frac{1}{ 2 } \right)
						\norm{ \vec{p} - \vec{c} - \bsym{\zeta}(\vec{p}) }_\infty. 
		\end{equation*}
		To prove (ii), note that if Eqn. (\ref{LGBound}) holds, 
		\begin{equation*}
			\norm{ \vec{p} - \vec{c} - \bsym{\zeta}(\vec{p}) }_\infty
				\geq \big( 1 - \norm{\bsym{\Omega}(\vec{p})}_\infty \big) \norm{ \vec{p} - \vec{c} - \bsym{\eta}(\vec{p}) }_\infty
				\geq \varepsilon \norm{ \vec{p} - \vec{c} - \bsym{\eta}(\vec{p}) }_\infty.
		\end{equation*}
	\endproof
	
	\begin{lemma}
		\label{Coercivity}
		Suppose $\varsigma_* = \infty$ and Assumption \ref{MixedLogitUtilityAssump}-\ref{LeibnizRuleCondition} and \ref{BoundedAssumption} hold. Then 
		\begin{equation*}
			\lim_{\norm{\vec{p}}_\infty \to \infty} \norm{ \vec{p} - \vec{c} - \bsym{\eta}(\vec{p}) }_\infty 
				= \infty
				= \lim_{\norm{\vec{p}}_\infty \to \infty} \norm{ \vec{p} - \vec{c} - \bsym{\zeta}(\vec{p}) }_\infty. 
		\end{equation*}
	\end{lemma}
	\proof{}
		The norm-coercivity of $\bsym{\eta}$ is a trivial consequence of the boundedness of $\bsym{\eta}$ under Assumption \ref{BoundedAssumption}. The norm-coercivity of $\bsym{\zeta}$ then follows from Lemma \ref{EtaCoercivityImpliesZetaCoercivity}.
	\endproof
	
	We now weaken Assumption \ref{BoundedAssumption}'s Eqn. (\ref{LPBound}). 
	\begin{assumption}
		\label{BoundedAssumption2}
		Suppose that $\varsigma_* = \infty$ and
		\begin{align}
			\label{LPBound2}
			\lim_{ M \uparrow \infty } 
				\sup \left\{ \frac{ \norm{ \bsym{\Lambda}(\vec{p})\inv\vec{P}(\vec{p}) }_\infty }
								{ \norm{ \vec{p} }_\infty } 
					: \vec{p} \in \set{P}^J, \norm{ \vec{p} }_\infty \geq M \right\}
				= \delta \in [0,1). 
		\end{align}
	\end{assumption}
	Note that the limit is of a non-increasing sequence of non-negative numbers, and thus exists. 
	
	\begin{lemma}
		Assuming Eqn. (\ref{LPBound2}) is equivalent to assuming that for any sequence $\vec{p}_n$ with $\norm{\vec{p}_n}_\infty \to \infty$, $\lim_{n\to\infty} \norm{ \bsym{\Lambda}(\vec{p}_n)\inv\vec{P}(\vec{p}_n) }_\infty/\norm{\vec{p}_n}_\infty \leq \delta$. 
	\end{lemma}
	\proof{}
		If Eqn. (\ref{LPBound2}) holds, then for any $\varepsilon > 0$ there exists an $M > 0$ such that
		\begin{equation*}
			\sup \left\{ \frac{ \norm{ \bsym{\Lambda}(\vec{p})\inv\vec{P}(\vec{p}) }_\infty }
								{ \norm{ \vec{p} }_\infty } 
					: \vec{p} \in \set{P}^J, \norm{ \vec{p} }_\infty \geq M \right\}
				< \delta + \varepsilon. 
		\end{equation*}
		If $\norm{ \vec{p}_n }_\infty \to \infty$, then there is also an $N_\epsilon$ such that $\norm{ \vec{p}_n }_\infty \geq M$ for all $n > N_\epsilon$. Thus
		\begin{equation*}
			 \frac{ \norm{ \bsym{\Lambda}(\vec{p}_n)\inv\vec{P}(\vec{p}_n) }_\infty }
								{ \norm{ \vec{p}_n }_\infty } 
				< \delta + \varepsilon
		\end{equation*}
		for all $n > N_\epsilon$, and thus
		\begin{equation*}
			\lim_{n\to\infty} \left[ \frac{ \norm{ \bsym{\Lambda}(\vec{p}_n)\inv\vec{P}(\vec{p}_n) }_\infty }
								{ \norm{ \vec{p}_n }_\infty } \right] \leq \delta. 
		\end{equation*}
		Conversely, if Eqn. (\ref{LPBound2}) fails, then there is a $\bar{M} > 0$ such that 
		\begin{equation*}
			S(M) = \sup \left\{ \frac{ \norm{ \bsym{\Lambda}(\vec{p})\inv\vec{P}(\vec{p}) }_\infty }
								{ \norm{ \vec{p} }_\infty } 
					: \vec{p} \in \set{P}^J, \norm{ \vec{p} }_\infty \geq M \right\}
				\geq 1
		\end{equation*}
		for all $M \geq \bar{M}$. We can thus choose $\vec{p}_M$ with $\norm{ \vec{p}_M }_\infty \geq M$ satisfying
		\begin{equation*}
			1 - \frac{ \norm{ \bsym{\Lambda}(\vec{p}_M)\inv\vec{P}(\vec{p}_M) }_\infty }
								{ \norm{ \vec{p}_M }_\infty } 
				\leq S(M) - \frac{ \norm{ \bsym{\Lambda}(\vec{p}_M)\inv\vec{P}(\vec{p}_M) }_\infty }
								{ \norm{ \vec{p}_M }_\infty }
				\leq \frac{1}{M}. 
		\end{equation*}
		In other words,
		\begin{equation*}
			 \frac{ \norm{ \bsym{\Lambda}(\vec{p}_M)\inv\vec{P}(\vec{p}_M) }_\infty }
								{ \norm{ \vec{p}_M }_\infty }
				\geq 1 - \frac{1}{M}
		\end{equation*}
		for all $M \geq \bar{M}$, and thus
		\begin{equation*}
			\lim_{M \to \infty} 
				\left[ \frac{ \norm{ \bsym{\Lambda}(\vec{p}_M)\inv\vec{P}(\vec{p}_M) }_\infty }
								{ \norm{ \vec{p}_M }_\infty } \right]
				\geq 1. 
		\end{equation*}
		Hence the ``sequence version'' of Eqn. (\ref{LPBound2}) fails, and thus by contraposition the sequence version and Eqn. (\ref{LPBound2}) are identical.
	\endproof

	
	Next we note that Eqn. (\ref{LPBound2}) {\em weakens} Eqn. (\ref{LPBound}). 
	\begin{lemma}
		If Eqn. (\ref{LPBound}) holds, then Eqn. (\ref{LPBound2}) holds.
	\end{lemma}
	\proof{}
		If $\sup \{ \norm{ \bsym{\Lambda}(\vec{p})\inv\vec{P}(\vec{p}) }_\infty : \vec{p} \in \set{P}^J \} \leq M$, then 
		\begin{equation*}
			S(L) = \sup \left\{ \frac{ \norm{ \bsym{\Lambda}(\vec{p})\inv\vec{P}(\vec{p}) }_\infty }
							{ \norm{ \vec{p} }_\infty } 
				: \vec{p} \in \set{P}^J, \norm{ \vec{p} }_\infty \geq L \right\}
			\leq \frac{M}{L}. 
		\end{equation*}
		Thus $\lim_{L \to \infty} S(L) = 0$, a special case of Eqn. (\ref{LPBound2}).
	\endproof
	
	Now we prove the alternative coercivity result. 
	\begin{lemma}
		Suppose $\varsigma_* = \infty$ and Assumptions \ref{MixedLogitUtilityAssump}-\ref{LeibnizRuleCondition}, \ref{BoundedOmegaProductsAssumption} and \ref{BoundedAssumption2} hold. Then 
		\begin{equation*}
			\lim_{\norm{\vec{p}}_\infty \to \infty} \norm{ \vec{p} - \vec{c} - \bsym{\eta}(\vec{p}) }_\infty 
				= \infty
				= \lim_{\norm{\vec{p}}_\infty \to \infty} \norm{ \vec{p} - \vec{c} - \bsym{\zeta}(\vec{p}) }_\infty. 
		\end{equation*}
	\end{lemma}
	\proof{}
		We prove the claim for $\bsym{\zeta}$; the result for $\bsym{\eta}$ then follows from Lemma \ref{LGBound2ImpliesLGBound}. Note that
		\begin{align*}
			&\abs{ p_k - c_k - \sum_{j \in \set{J}_{f(k)}} \omega_{k,j}(\vec{p})(p_j - c_j) 
						- \frac{P_k(\vec{p})}{\lambda_k(\vec{p})} } \\
			&\quad\quad\quad\quad\quad\quad
				= \abs{ \left( 1 - \frac{P_k(\vec{p})}{p_k\lambda_k(\vec{p})} \right) p_k
						- c_k - \sum_{j \in \set{J}_{f(k)}} \omega_{k,j}(\vec{p})(p_j - c_j) } \\
			&\quad\quad\quad\quad\quad\quad
				\geq \abs{ \; \abs{ 1 - \frac{P_k(\vec{p})}{p_k\lambda_k(\vec{p})} } p_k
						- \Bigg | c_k + \sum_{j \in \set{J}_{f(k)}} \omega_{k,j}(\vec{p})(p_j - c_j) \Bigg | \; }
		\end{align*}
		Suppose that $p_k \to \infty$. By assumption, 
		\begin{equation*}
			\lim_{n \to \infty} \left[ 1 - \frac{P_k(\vec{p})}{p_k\lambda_k(\vec{p})} \right]
				\geq 1 - \delta > 0
		\end{equation*}
		while the second term is bounded. Thus 
		\begin{align*}
			\abs{ p_k - c_k - \sum_{j \in \set{J}_{f(k)}} \omega_{k,j}(\vec{p})(p_j - c_j) 
						- \frac{P_k(\vec{p})}{\lambda_k(\vec{p})} }
				\to \infty.
		\end{align*}
	\endproof
	
	Note that since we did not require that Eqn. (\ref{LGBound}) held, $\bsym{\zeta}$ need not be bounded for $\vec{F}_\eta$ and $\vec{F}_\zeta$ to be coercive. 

	
	\subsection{Fixed-Point Iteration}
	\label{SUBSEC:FPIs}
	
	In addition to applications of Newton's method, the BLP- and $\bsym{\zeta}$-markup equations suggest applying fixed-point iteration to solve for equilibrium prices. We examine fixed-point iterations based on both equations. 
	
	
	\subsubsection{$\zeta$ Fixed-Point Iteration}
	
	The fixed-point iteration $\vec{p} \leftarrow \vec{c} + \bsym{\zeta}(\vec{p})$ based on the $\bsym{\zeta}$-markup equation, here abbreviated $\bsym{\zeta}$-FPI, can efficiently compute equilibrium prices for some problems. $\bsym{\zeta}$-FPI has relatively efficient steps because no linear systems need to be solved, unlike every other method listed in Table \ref{TAB:MethodSummary}. While we are not aware of a general convergence proof for $\bsym{\zeta}$-FPI, this iteration has converged reliably on test problems including the examples in \cite{Morrow10}. 
	
	The first observation we make is that the $\bsym{\zeta}$-FPI steps always point in directions of ``myopic gradient ascent.''  
	\begin{lemma}
	 	Let $\vec{p} \in (0,\varsigma_*)^J$, and let $\bsym{\delta}\vec{p} = \vec{c} + \bsym{\zeta}(\vec{p}) - \vec{p}$ denote the $\bsym{\zeta}$-FPI step. Then
		\begin{equation*}
			\frac{1}{ \max_j \abs{ \lambda_j(\vec{p}) } }
				\leq \frac{(\tilde{\nabla}\hat{\pi})(\vec{p})^\top\bsym{\delta}\vec{p}}
					{(\tilde{\nabla}\hat{\pi})(\vec{p})^\top(\tilde{\nabla}\hat{\pi})(\vec{p})}
				\leq \frac{1}{ \min_j \abs{ \lambda_j(\vec{p}) } }. 
		\end{equation*}
		Similarly, let $\theta(\vec{p})$ denote the angle between $\bsym{\delta}\vec{p}$ and $(\tilde{\nabla}\hat{\pi})(\vec{p})$, and suppose $\vec{p}$ is not simultaneously stationary. Then
		\begin{equation*}
			\cos \theta(\vec{p}) 
				\geq \frac{ \min_j \abs{\lambda_j(\vec{p})} }
							{ \max_j \abs{\lambda_j(\vec{p})} }. 
		\end{equation*}
	\end{lemma}
	\proof{}
		Both results follows directly from the equation $(\tilde{\nabla}\hat{\pi})(\vec{p}) = \abs{ \bsym{\Lambda}(\vec{p}) } \bsym{\delta}\vec{p}$ where $\abs{\bsym{\Lambda}(\vec{p})}$ denotes the absolute value of the components of $\bsym{\Lambda}(\vec{p})$. 
	\endproof
	Specifically, the $\bsym{\zeta}$-FPI steps have a positive projection onto the combined gradient, and cannot become orthogonal to the combined gradient over any sequence of non-simultaneously stationary prices that stay in $(0,\varsigma_*)^J$. 
	
	If $F = 1$, and the equilibrium problem is an optimization problem, this implies $\bsym{\zeta}$-FPI has steps that point in gradient ascent directions and, when properly scaled, converge to local maximizers of profit. More specifically, $\bsym{\zeta}$-FPI cannot converge to minimizers of profits. This may generate the properties of $\bsym{\zeta}$-FPI observed in Example 10 from \cite{Morrow10}. 
	
	\begin{corollary}
		Let Assumptions \ref{MixedLogitUtilityAssump}-\ref{BoundedAssumption} hold, and suppose $\{ \vec{p}^{(n)} \}_{n=1}^\infty$ is the $\bsym{\zeta}$-FPI sequence. Then $\{ \vec{p}^{(n)} \}_{n=1}^\infty$ is bounded. 
	\end{corollary}
	\proof{}
		By Lemma \ref{ZetaGrowth}, for any sufficiently large $M > 0$ we can find some $L > 0$ such that
		\begin{equation*}
			\norm{ \bsym{\zeta}(\vec{p}) }_\infty < \norm{ \vec{p} - \vec{c} }_\infty - M
			\quad\text{for all}\quad
			\norm{ \vec{p} - \vec{c} }_\infty > L. 
		\end{equation*}
		If the $\bsym{\zeta}$-FPI sequence diverges, then for any such $L$ there is an $N$ such that
		\begin{equation*}
			\norm{ \vec{p}^{(n)} - \vec{c} }_\infty > L
			\quad\text{for all}\quad
			n > N. 
		\end{equation*}
		But then 
		\begin{equation*}
			\norm{ \vec{p}^{(n+1)} - \vec{c} }_\infty 
				= \norm{ \bsym{\zeta}(\vec{p}^{(n)}) }_\infty 
				< \norm{ \vec{p}^{(n)} - \vec{c} }_\infty - M
				< \norm{ \vec{p}^{(n)} - \vec{c} }_\infty
			\quad\text{for all}\quad
			n > N, 
		\end{equation*}
		which states that the $\bsym{\zeta}$-FPI sequence is decreasing. This is a contradiction of the hypothesis that the $\bsym{\zeta}$-FPI sequence diverges.
	\endproof
	
	To implement $\bsym{\zeta}$-FPI, one simply needs to iterate the assignment $\vec{p} \leftarrow \vec{c} + \bsym{\zeta}(\vec{p})$ where Eqn. (\ref{ZetaMap}) defines $\bsym{\zeta}(\vec{p})$. As shown in Table \ref{TAB:FlopCounts} below, integral approximations, rather than the actual computation of the step, drive the computational burden. Given a price vector, utilities, and utility derivatives, computing $\vec{P}(\vec{p})$, $\bsym{\Lambda}(\vec{p})$, and $\tilde{\bsym{\Gamma}}(\vec{p})$ for a set of $S$ samples requires $\order( S \sum_{f=1}^F J_f^2 )$ floating point operations (\texttt{flops}), while the fixed-point step itself only requires $\order( \sum_{f=1}^F J_f^2 )$ \texttt{flops}. Note that computing the fixed-point {\em step} $\vec{c} + \bsym{\zeta}(\vec{p})$ requires an equivalent amount of work as computing the combined gradient $(\tilde{\nabla}\hat{\pi})(\vec{p})$. Furthermore, because $\bsym{\Lambda}(\vec{p})$ is a diagonal matrix, no serious obstacles to computing the fixed point step arise as $J$ becomes large. 
	
	
	\subsubsection{$\eta$ Fixed-Point Iteration}
	
	The fixed-point iteration $\vec{p} \leftarrow \vec{c} + \bsym{\eta}(\vec{p})$, abbreviated $\bsym{\eta}$-FPI, based on the BLP-markup equation need not converge. Example \ref{EX:ConvExam} below, repeated from \cite{Morrow10}, gives a case in which $\bsym{\eta}$ can fail to be even locally convergent. 
	
	\begin{example}
		\label{EX:ConvExam}
		Consider multi-product monopoly pricing with a simple Logit model having $u_j(p) = - \alpha p + v_j$ for some $\alpha > 0$, any $v_j \in \R$, and $\vartheta > -\infty$. It is well known that for a single-product firm, unique profit-maximizing prices exist \citep{Anderson88, Milgrom90, Caplin91}. \cite{Morrow08} proves that profit-optimal prices $\vec{p}_*$ are unique for the multi-product case $-$ and even so with multiple firms $-$ even though profits are not quasi-concave \citep{Hanson96}. 
		
		In this example, $\bsym{\eta}$-FPI is not always locally convergent near $\vec{p}_*$, while $\bsym{\zeta}$-FPI is always {\em superlinearly} locally convergent. For an arbitrary continuously differentiable function $\vec{F}$ and $\vec{p}_* = \vec{F}(\vec{p}_*)$, $\vec{F}$ is contractive on some neighborhood of $\vec{p}_*$ in some norm $\norm{\cdot}$ if $\rho( (D\vec{F})(\vec{p}_*) ) < 1$ where $\rho(\vec{A})$ \citep{Ortega70}. We show that $\rho( (D\bsym{\eta})(\vec{p}_*) ) > 1$ may hold while $\rho( (D\bsym{\zeta})(\vec{p}_*) ) = 0$, where $\rho(\vec{A})$ denotes the spectral radius of the matrix $\vec{A}$. 
		
		The components of the BLP-markup function $\bsym{\eta}$ are given by $\eta_k(\vec{p}) = \alpha\inv(1 - \sum_{j=1}^J P_j^L(\vec{p}) )\inv$ for all $k$. From this formula the equation
		\begin{equation*}
			\rho( (D\bsym{\eta})(\vec{p}_*) )
				= \frac{ \sum_{j=1}^J P_j^L(\vec{p}_*) }{ 1 - \sum_{j=1}^J P_j^L(\vec{p}_*) }
				= \sum_{j=1}^J e^{u_j(p_{j,*}) - \vartheta}
		\end{equation*}
		can be derived. For valuations of the outside good, $\vartheta$, sufficiently close to $-\infty$, $\rho( (D\bsym{\eta})(\vec{p}_*) ) > 1$ can hold; see \cite{Morrow10} for details. 
		
		To prove the claim regarding $\rho( (D\bsym{\zeta})(\vec{p}_*) )$, note that $\zeta_k(\vec{p}) = \hat{\pi}(\vec{p}) + 1/\alpha$, and thus $(D_l\zeta_k)(\vec{p}_*) = (D_l\hat{\pi})(\vec{p}_*) = 0$ for all $k,l$.
		
	\end{example}
	
	Even if the BLP-markup equation does generate a convergent fixed-point iteration, evaluating $\bsym{\eta}$ involves the solution of $F$ linear systems that grow in size with the number of products offered by the firms. The work required to evaluate $\bsym{\eta}$ using a direct method like \texttt{PLU} or \texttt{QR} factorization is $\order([\max_f J_f ]^3)$, given values of $\vec{P}(\vec{p})$, $\bsym{\Lambda}(\vec{p})$, and $\tilde{\bsym{\Gamma}}(\vec{p})$ as approximated using simulation. The work to evaluate $\bsym{\zeta}$ is only $\order([ \max_f J_f ]^2)$ given $\vec{P}(\vec{p})$, $\bsym{\Lambda}(\vec{p})$, and $\tilde{\bsym{\Gamma}}(\vec{p})$ (Table \ref{TAB:FlopCounts}). Generally speaking, function evaluations must be cheap for the linear convergence of fixed-point iterations to result in faster computations than the superlinearly or quadratically convergent variants of Newton's method. 

	
	\subsection{Practical Considerations}
	\label{SUBSEC:PreliminaryConsiderations}
	
	This section addresses several practical considerations. 
	
	
	\subsubsection{Simulation}
	\label{SUBSEC:FiniteSample}
	
	Any method for computing equilibrium prices under Mixed Logit models faces a common obstacle: the integrals that define the choice probabilities ($\vec{P}$) and their derivatives ($\bsym{\Lambda},\tilde{\bsym{\Gamma}}$) cannot be computed exactly. We employ finite-sample versions of the methods discussed below by drawing $S \in \N$ samples from the demographic distribution and applying the method to the finite-sample model thus generated. Particularly, these samples are used to compute approximate $\vec{P}(\vec{p})$, $\bsym{\Lambda}(\vec{p})$, and $\tilde{\bsym{\Gamma}}(\vec{p})$; see Table \ref{TAB:FlopCounts}. These samples are kept fixed for {\em all} steps of the method and, in principle, can be generated in any way. We draw directly from the demographic distribution, although importance and quasi-random sampling (e.g., see \cite{Train03}) can also be employed. The Law of Large Numbers motivates this widely-used approach to econometric analysis (e.g., see \cite{McFadden89} and \cite{Draganska04}). While all numerical approaches for computing equilibrium prices described here rely on a Law of Large Numbers for simultaneously stationary prices, we do not provide a formal convergence theorem. We do provide numerical evidence that computed equilibrium prices based on the fixed-point iteration for our examples do indeed follow such a law. 
	
	
	\subsubsection{Truncation of Low Purchase Probability Products}
	\label{SUBSEC:Truncation}
	
	All of the methods we implement can be built to ignore products with excessively low choice probabilities. That is, one can ignore price updates for all products with $P_j(\vec{p}) \leq \varepsilon_P$, where $\varepsilon_P$ is some small value (say $10^{-10}$). Products with a choice probability this small (or smaller) need not be considered a part of the market in the price equilibrium computations. For example, \cite{Wards} reports total sales of cars and light trucks during 2005 as $N = 16,947,754$. Particularly, 7,667,066 cars and 9,280,688 light trucks. Because expected demand is defined by $\Expect[ Q_j(\vec{p}) ] = N P_j(\vec{p})$, any $\varepsilon_P \leq 0.5 * N\inv \approx 3 \times 10^{-8}$ ignores any vehicle that, as priced, is not {\em expected} to have a {\em single} customer out of the millions of customers that bought or considered buying new vehicles. There are also technical reasons for this truncation. Particularly, $\bsym{\Lambda}(\vec{p})$ and $(D\tilde{\nabla}\hat{\pi})(\vec{p})$ become singular as $P_j(\vec{p}) \to 0$, for any $j$. Truncating avoids this non-singularity and hopefully helps conditioning. 
	
	
	\subsubsection{Termination Conditions}
	\label{SUBSEC:Termination}
	
	We terminate {\em all} iterations with the numerical simultaneous stationarity condition $\norm{ (\tilde{\nabla}\hat{\pi})(\vec{p}) }_\infty \leq \varepsilon_T$ where $\varepsilon_T$ is some small number (e.g., $10^{-6}$). Note that a standard application of Newton's method to solve $\vec{F}_\eta(\vec{p}) = \vec{0}$ or $\vec{F}_\zeta(\vec{p}) = \vec{0}$ would terminate when either
	\begin{equation}
		\label{TermInequalities}
		\norm{ \vec{p} - \vec{c} - \bsym{\eta}(\vec{p}) }_\infty \leq \varepsilon_T
		\quad\quad\text{or}\quad\quad
		\norm{ \vec{p} - \vec{c} - \bsym{\zeta}(\vec{p}) }_\infty \leq \varepsilon_T, 
	\end{equation}
	respectively. For example, \cite{Aguirregabiria06} use the condition $\norm{ \vec{p} - \vec{c} - \bsym{\eta}(\vec{p}) }_\infty \leq \varepsilon_T$. Ensuring that Eqn. (\ref{TermInequalities}) holds does {\em not} necessarily imply that $\norm{ (\tilde{\nabla}\hat{\pi})(\vec{p}) }_\infty \leq \varepsilon_T$, the strictly interpreted first-order condition. 
	
	Because
	\begin{equation*}
		(\tilde{D}\vec{P})(\vec{p})^\top( \vec{p} - \vec{c} - \bsym{\eta}(\vec{p}) )
			= (\tilde{\nabla}\hat{\pi})(\vec{p})
			= \bsym{\Lambda}(\vec{p})( \vec{p} - \vec{c} - \bsym{\zeta}(\vec{p}) ), 
	\end{equation*}
	it is easy to terminate all methods, CG-NM, $\bsym{\eta}$-NM, $\bsym{\zeta}$-NM, and $\bsym{\zeta}$-FPI, when $\norm{ (\tilde{\nabla}\hat{\pi})(\vec{p}) }_\infty \leq \varepsilon_T$. While this is done here to ensure consistency in our comparisons of different methods, $\norm{ (\tilde{\nabla}\hat{\pi})(\vec{p}) }_\infty \leq \varepsilon_T$ should always be the termination condition for price equilibrium computations. 
	
	Three other standard termination conditions are used \citep{Brown90, Dennis96}. We terminate the iteration if the (relative) step length becomes too small, if a maximum number of iterations is exceeded, or if an exceptional event occurs (e.g. division by zero). These three conditions are considered ``failure'' as the iteration has failed to compute a numerically simultaneously stationary point in the sense of the first termination condition. 
	
	
	\subsubsection{Second-Order Conditions.} 
	\label{SUBSEC:Sufficiency}
	
	Each method in Table \ref{TAB:MethodSummary} finds simultaneously stationary points, rather than local equilibria. Unlike in optimization, there is no {\em a priori} assurance that first-order iterative methods for equilibrium problems will converge to certain types of stationary points. Thus in computing equilibria it is vitally important to check the second-order sufficient conditions to verify that a local equilibrium has indeed been found. 
	
	In local equilibrium every firm's profit Hessian, $(D_f\nabla_f\hat{\pi}_f)(\vec{p})$, should also be negative definite. The formulas given in Proposition \ref{MixedLogitSecondDerivatives} below provide an expression for $(D_f\nabla_f\hat{\pi}_f)(\vec{p})$ that we use to check the second-order sufficient condition. Cholesky factorization, rather than direct approximation of the spectrum, is used to test the negative definiteness of $(D_f\nabla_f\hat{\pi}_f)(\vec{p})$ \citep{Golub96}. 
	
	
	\subsubsection{Computational Burden}
	\label{SUBSEC:CompBurden}
	
	Table \ref{TAB:FlopCounts} reviews the formulae and computational burden of computing $(\tilde{\nabla}\hat{\pi})$, $\bsym{\eta}$, and $\bsym{\zeta}$. 
	
	\begin{table}
		\begin{small}
		\begin{center}
		\caption[Work required to evaluate $(\tilde{\nabla}\hat{\pi})$, $\bsym{\eta}$, and $\bsym{\zeta}$. ]{Work required to evaluate $(\tilde{\nabla}\hat{\pi})$, $\bsym{\eta}$, and $\bsym{\zeta}$ given $S$ samples $\{ \bsym{\theta}_s \}_{s=1}^S \subset \set{T}$, an $S \times J$ matrix $\vec{L}(\vec{p})$ of Logit choice probabilities ($(\vec{L}(\vec{p}))_{s,j} = P_j^L(\bsym{\theta}_s,\vec{p})$), and an $S \times J$ matrix of utility derivatives $\vec{D}(\vec{p})$ ($(\vec{D}(\vec{p}))_{s,j} = (Dw_j)(\bsym{\theta}_s,p_j)$). The first section gives work required for sample-average approximations to $\vec{P}(\vec{p})$, $\bsym{\Lambda}(\vec{p})$, and $\tilde{\bsym{\Gamma}}(\vec{p})$. The second section takes $\vec{P}(\vec{p})$, $\bsym{\Lambda}(\vec{p})$, and $\tilde{\bsym{\Gamma}}(\vec{p})$ as given.}
		\label{TAB:FlopCounts}
		\begin{tabular}{lcc}
			\\
			Quantity & Formula & \texttt{flops} \\ \hline
			\\
			$\vec{P}(\vec{p})$ & $S\inv \vec{L}(\vec{p})^\top\vec{1}$ & $SJ$ \\
			$\vec{V}(\vec{p})$ & $\vec{L}(\vec{p}) \cdot \vec{D}(\vec{p})$$^{\text{(a)}}$ & $SJ$ \\
			$\bsym{\Lambda}(\vec{p})$ & $S\inv \vec{V}(\vec{p})^\top\vec{1}$ & $SJ$ \\
			$\tilde{\bsym{\Gamma}}(\vec{p})$ & $S\inv \vec{L}(\vec{p})^\top\vec{V}(\vec{p})$ & $2 S \sum_{f=1}^F J_f^2$ \\
			\\
			\multicolumn{2}{l}{Total work to compute $\vec{P}(\vec{p})$, $\bsym{\Lambda}(\vec{p})$, and $\tilde{\bsym{\Gamma}}(\vec{p})$} & $S\left( 3J + 2 \sum_{f=1}^F J_f^2 \right)$ \\
			\\ \hline
			\\
			$\bsym{\zeta}(\vec{p})$
				& $\tilde{\bsym{\Omega}}(\vec{p})(\vec{p}-\vec{c}) 
						- \bsym{\Lambda}(\vec{p})\inv\vec{P}(\vec{p})$
				& $2 \sum_{f=1}^F J_f^2 + 4J$ \\
			$\bsym{\eta}(\vec{p})$
				& $(\vec{I} - \tilde{\bsym{\Omega}}(\vec{p})) \bsym{\eta}(\vec{p}) = - \bsym{\Lambda}(\vec{p})\inv\vec{P}(\vec{p})$
				& $\left( \frac{4}{3} \right) \sum_{f=1}^F J_f^3 
					+ \left( \frac{7}{2} \right) \left( \sum_{f=1}^F J_f^2 + J \right)
					- 2$
				\\
			\\
			$(\tilde{\nabla}\hat{\pi})(\vec{p})$
				& $(\bsym{\Lambda}(\vec{p}) - \tilde{\bsym{\Gamma}}(\vec{p})^\top)(\vec{p}-\vec{c}) + \vec{P}(\vec{p})$
				& $2 \sum_{f=1}^F J_f^2 + 5J$ \\
				& $= \bsym{\Lambda}(\vec{p})( \vec{p} - \vec{c} - \bsym{\zeta}(\vec{p}) )$
				& $2 \sum_{f=1}^F J_f^2 + 6J$ \\
			\\ \hline
			\multicolumn{3}{l}{(a) ``$\cdot$'' here denotes element-by-element multiplication.}
		\end{tabular}
		\end{center}
		\end{small}
	\end{table}
	
	Computing $\bsym{\eta}$ and applying Newton's method to $\vec{F}_\eta$ requires solving linear systems. We give some more details regarding these computations here. As stated above, the linear system
	\begin{equation*}
		(\vec{I} - \tilde{\bsym{\Omega}}(\vec{p}))\bsym{\eta}(\vec{p})
			= - \bsym{\Lambda}(\vec{p})\inv\vec{P}(\vec{p})
	\end{equation*}
	should be used to solve for $\bsym{\eta}(\vec{p})$. Note also that only the systems 
	\begin{equation*}
		(\vec{I} - \bsym{\Omega}_f(\vec{p}) )\bsym{\eta}_f(\vec{p}) 
			= - \bsym{\Lambda}_f(\vec{p})\inv\vec{P}_f(\vec{p})
	\end{equation*}
	for all $f$ need be solved. Of course, our condition bound applies within firms as well: 
	\begin{equation*}
		\kappa_2 \big( (D_f\vec{P}_f)(\vec{p})^\top \big)
			\geq \left( \frac{ \max_{j \in \set{J}_f} \abs{ \lambda_j(\vec{p}) } }{ \min_{j \in \set{J}_f} \abs{ \lambda_j(\vec{p}) }  } \right)
					\kappa_2 \big( \vec{I} - \bsym{\Omega}_f(\vec{p}) \big). 
	\end{equation*}
	If Householder QR factorization is used to solve these systems, then computing $\bsym{\eta}(\vec{p})$ from $\vec{P}(\vec{p})$, $\bsym{\Lambda}(\vec{p})$, and $\tilde{\bsym{\Gamma}}(\vec{p})$ requires $\order(\sum_{f=1}^F J_f^3)$ \texttt{flops} (Table \ref{TAB:FlopCounts}). 

	This is a significant increase in computational effort relative to computing $\bsym{\zeta}(\vec{p})$ or $(\tilde{\nabla}\hat{\pi})(\vec{p})$. The diagonal dominance of $\vec{I} - \tilde{\bsym{\Omega}}(\vec{p})$, indeed of $(\tilde{D}\vec{P})(\vec{p})$ itself, suggests that Jacobi, Gauss-Seidel, and Successive Over-Relaxation (SOR) iterations \citep{Golub96} may be a relatively efficient way to compute $\bsym{\eta}$. 
	
	Additional work is required to compute $(D\bsym{\eta})(\vec{p})$, if this is to be used in Newton's method. Though it requires solving a matrix-linear system of the type $(\tilde{D}\vec{P})(\vec{p})(D\bsym{\eta})(\vec{p}) = \vec{B}(\vec{p})$, the required matrix factorizations of $\vec{I} - \bsym{\Omega}_f(\vec{p})$ need only be computed once to compute both $\bsym{\eta}$ and $(D\bsym{\eta})$, but must be updated for each vector of prices. 
	
	
	\subsection{Computing Jacobian Matrices for Newton's Method}
	\label{SUBSEC:Jacobians}
	
	Standard ``exact'' or Quasi-Newton methods to solve $\vec{F}(\vec{x}) = \vec{0}$ either always or periodically require the Jacobian matrix $(D\vec{F})(\vec{x})$. Using finite differences to approximate Jacobian matrices requires $J$ evaluations of the function $\vec{F}$, an unacceptable workload. In the 993 vehicle example from \cite{Morrow10}, approximating $(D\vec{F})(\vec{x})$ {\em once} with finite differences would take roughly 993 evaluations of $\vec{F}$, when the work of less than 50 evaluations appears to sufficient to converge to equilibrium prices using the $\bsym{\zeta}$-FPI. 
	
	We recommend directly approximating $(D\vec{F})(\vec{x})$ using integral expressions for $(D\tilde{\nabla}\hat{\pi})(\vec{p})$, $(D\bsym{\eta})(\vec{p})$, and $(D\bsym{\zeta})(\vec{p})$ provided below. An alternative is to use automatic differentiation, but we are skeptical that this would in fact be faster than the direct formulae provided here. 
	
	
	\subsubsection{Jacobian of the Combined Gradient}
	
	Assuming a second application of the Leibniz Rule holds, we can derive integral expressions for the second derivatives $(D_lD_k\hat{\pi}_{f(k)})(\vec{p})$ through
	\begin{equation*}
		\big( (D\tilde{\nabla}\hat{\pi})(\vec{p}) \big)_{k,l}
			= (D_l D_k \hat{\pi}_{f(k)})(\vec{p})
			= \int (D_l D_k \hat{\pi}_{f(k)}^L)(\bsym{\theta},\vec{p}) d\mu(\bsym{\theta}). 
	\end{equation*}
	
	\begin{proposition}
		\label{MixedLogitSecondDerivatives}
		Let $w$ be twice continuously differentiable in $p$ and suppose a second application of the Leibniz Rule holds for the Mixed Logit choice probabilities at $\vec{p}$. Set 
		\begin{align*}
			\phi_{k,l}(\vec{p}) 
				&= \int (Dw_k)(\bsym{\theta},p_k)P_k^L(\bsym{\theta},\vec{p})
						P_l^L(\bsym{\theta},\vec{p})(Dw_l)(\bsym{\theta},p_l) d\mu(\bsym{\theta}) \\
			\psi_{k,l}(\vec{p}) 
				&= \int (Dw_k)(\bsym{\theta},p_k)P_k^L(\bsym{\theta},\vec{p})
						\hat{\pi}_{f(k)}^L(\bsym{\theta},\vec{p}) P_l^L(\bsym{\theta},\vec{p})(Dw_l)(\bsym{\theta},p_l) d \mu(\bsym{\theta}) \\
			\chi_k(\vec{p})
				&= \left( \frac{1}{2} \right) \int \big( (D^2w_k)(\bsym{\theta},p_k) + (Dw_k)(\bsym{\theta},p_k)^2 \big) \\ 
				&\quad\quad\quad\quad\quad\quad
					\times P_k^L(\bsym{\theta},\vec{p}) 
						\big( (p_k - c_k) - \hat{\pi}_{f(k)}^L(\bsym{\theta},\vec{p}) \big) d \mu(\bsym{\theta})
		\end{align*}
		\begin{itemize}
			\item[(i)] {\em Component form:} Setting
			\begin{equation*}
				\xi_{k,l}(\vec{p}) = \delta_{k,l}(\lambda_k(\vec{p}) + \chi_k(\vec{p})) - \gamma_{k,l}(\vec{p}) - (p_k-c_k)\varphi_{k,l}(\vec{p})
			\end{equation*}
			we have
			\begin{align*}
				(D_lD_k\hat{\pi}_{f(k)})(\vec{p})
					= \xi_{k,l}(\vec{p}) + 2\psi_{k,l}(\vec{p}) + \delta_{f(k),f(l)} \xi_{l,k}(\vec{p})
			\end{align*}
			\item[(ii)] {\em Matrix form:} Let $\bsym{\Phi}(\vec{p})$, $\bsym{\Psi}(\vec{p})$ and $\vec{X}(\vec{p}) = \diag(\bsym{\chi}(\vec{p}))$ be the matrices of these quantities. Also set
			\begin{equation*}
				\bsym{\Xi}(\vec{p}) 
					= \bsym{\Lambda}(\vec{p})  
						- \bsym{\Gamma}(\vec{p})
						- \diag( \vec{p} - \vec{c} ) \bsym{\Phi}(\vec{p}) 
						+ \vec{X}(\vec{p}). 
			\end{equation*}
			and 
			\begin{equation*}
				( \tilde{\bsym{\Xi}}(\vec{p}) )_{k,l}
					= \left\{ \begin{aligned}
						&\xi_{k,l}(\vec{p}) &&\quad \text{if } f(k) = f(l) \\
						&\quad 0 &&\quad\text{if } f(k) \neq f(l)
					\end{aligned} \right .
			\end{equation*}
			Then
			\begin{equation}
				(D\tilde{\nabla}\hat{\pi})(\vec{p}) 
					= \bsym{\Xi}(\vec{p}) + 2 \bsym{\Psi}(\vec{p}) + \tilde{\bsym{\Xi}}(\vec{p})^\top. 
			\end{equation}
		\end{itemize}
	\end{proposition}
	
	\proof{}
		To see that this only relies on a second application of the Leibniz Rule to the {\em choice probabilities}, note that
	\begin{equation*}
		(D_lD_k\hat{\pi}_{f(k)})(\vec{p})
			= \sum_{j \in \set{J}_{f(k)}} (D_lD_kP_j)(\vec{p})(p_j - c_j) + \delta_{f(k),f(l)} (D_kP_l)(\vec{p}) + (D_lP_k)(\vec{p})
	\end{equation*}
	and thus the continuous second-order differentiability of $\hat{\pi}_f(\vec{p})$ depends only on the second-order continuous differentiability of $\vec{P}_f$. This result is then an immediate consequence of the validity of the Leibniz Rule, if a bit tedious to derive. 
	\endproof
	
	The validity of a second application of the Leibniz Rule to the choice probabilities is ensured by the following condition. 
	
	\begin{proposition}
		\label{SecondLeibnizRule}
		Let $(u,\vartheta,\mu) = (w+v,\vartheta,\mu)$ be such that 
		\begin{itemize}
			\item[(i)] $w(\bsym{\theta},\vec{y},\cdot) : (0,\varsigma_*) \to \R$ is twice continuously differentiable for all $\vec{y} \in \set{Y}$ and $\mu$-a.e. $\bsym{\theta} \in \set{T}$ 
			\item[(ii)] for all $(\vec{y},p) \in \set{Y} \times (0,\varsigma_*)$, $\abs{ (D^2w)(\cdot,\vec{y},q) + (Dw)(\cdot,\vec{y},q)^2 } e^{u(\cdot,\vec{y},q) - \vartheta(\cdot)} : \set{T} \to [0,\infty)$ is uniformly $\mu$-integrable for all $q$ in some neighborhood of $p$. 
			\item[(iii)] for all $(\vec{y},p),(\vec{y}^\prime,p^\prime) \in \set{Y} \times (0,\varsigma_*)$, 
			\begin{equation*}
				\abs{(Dw)(\cdot,\vec{y},q)}e^{u(\cdot,\vec{y},q) - \vartheta(\cdot)} 
						e^{u(\cdot,\vec{y}^\prime,q^\prime) - \vartheta(\cdot)}\abs{(Dw)(\cdot,\vec{y}^\prime,q^\prime)} : \set{T} \to [0,\infty)
			\end{equation*}
			is uniformly $\mu$-integrable for all $(q,q^\prime)$ in some neighborhood of $(p,p^\prime)$. 
		\end{itemize}
		Then a second application of the Leibniz Rule holds for the Mixed Logit choice probabilities, which are also continuously differentiable on $(\vec{0},\varsigma_*\vec{1})$. 
	\end{proposition}
	
	This is proved in the same manner as Proposition \ref{LeibnizRule}.
	
	We also observe the following. 
	\begin{proposition}
		\label{ZeroSecondDers}
		If $P_k(\vec{p}) = 0$ then $(D_l D_k \hat{\pi}_{f(k)})(\vec{p}) = (D_k D_l \hat{\pi}_{f(l)})(\vec{p}) = 0$ for all $l \in \N(J)$. 
	\end{proposition}
	The proof follows from the derivative formulae given above. Of course, if $P_k(\vec{p}) = 0$ then $(D_k \hat{\pi}_{f(k)})(\vec{p}) = 0$ as well and we have the following situation: (i) the Newton system is consistent for any $s_k^N(\vec{p}) \in \R$ and (ii) $s_l^N(\vec{p})$ does not depend on $s_k^N(\vec{p})$ for all $l \in \N(J)\setminus\{k\}$. Thus, in practice one can restrict attention to the Newton step defined by the submatrix of $(D\tilde{\nabla}\hat{\pi})(\vec{p})$ formed by rows and columns indexed by $\{ j : P_j(\vec{p}) > \varepsilon_P \}$. 
	
	The formulae above give the following expression of the profit Hessians. 
	\begin{corollary}
		\label{Hessians}
		Let $w$ be twice continuously differentiable in $p$ and suppose a second application of the Leibniz Rule holds for the Mixed Logit choice probabilities. Firm $f$'s profit Hessian is given by
		\begin{equation*}
			(D_f\nabla_f\hat{\pi}_f)(\vec{p}) 
				= \bsym{\Xi}_{f,f}(\vec{p}) + 2 \bsym{\Psi}_{f,f}(\vec{p}) + \bsym{\Xi}_{f,f}(\vec{p})^\top. 
		\end{equation*}
	\end{corollary}
	
	
	\subsubsection{The $\bsym{\eta}$ map.}
	
	For $\vec{F}_\eta(\vec{p}) = \vec{p} - \vec{c} - \bsym{\eta}(\vec{p})$, we have $(D\vec{F}_\eta)(\vec{p}) = \vec{I} - (D\bsym{\eta})(\vec{p})$ where $(D\bsym{\eta})(\vec{p})$ solves the linear matrix equation 
	\begin{equation*}
		(\tilde{D}\vec{P})(\vec{p})^\top(D\bsym{\eta})(\vec{p}) 
			= - ( \vec{A}(\vec{p}) + (D\vec{P})(\vec{p}) ). 
	\end{equation*}
	Here $( \vec{A}(\vec{p}) )_{k,l} = \sum_{j \in \set{J}_{f(k)}} (D_lD_kP_j)(\vec{p})\eta_j(\vec{p})$. This is easily derived from the defining formula $(\tilde{D}\vec{P})(\vec{p})^\top \bsym{\eta}(\vec{p}) = - \vec{P}(\vec{p})$. 
	
	
	\subsubsection{The $\bsym{\zeta}$ map.}
	
	For $\vec{F}_\zeta(\vec{p}) = \vec{p} - \vec{c} - \bsym{\zeta}(\vec{p})$, we have $(D\vec{F}_\zeta)(\vec{p}) = \vec{I} - (D\bsym{\zeta})(\vec{p})$ where $(D\bsym{\zeta})(\vec{p})$ can be computed using the following formula: 
		\begin{align*}
			(D_l\zeta_k)
				&= \lambda_k\inv \Bigg[ \delta_{k,l} \left[ \int P_k^L \big( (D^2w_k) + (Dw_k)^2 \big) \left( \hat{\pi}_{f(k)}^L - \zeta_k \right) 
						- \lambda_k \right] \\
				&\quad\quad\quad\quad\quad\quad
					+ \zeta_k \phi_{k,l} + \gamma_{k,l} 
					+ \delta_{f(k),f(l)} \phi_{k,l} (p_l-c_l) + \delta_{f(k),f(l)} \gamma_{l,k} 
					- 2 \psi_{k,l} \Bigg]. 
		\end{align*}

%% file: gnhs.tex
	\section{The \texttt{GMRES}-Newton Hookstep Method}
	
	In this section we provide some details regarding the \texttt{GMRES}-Newton Hookstep method employed in \cite{Morrow10}. For complete details, see \cite{Morrow10b}. 
	
	
	\subsection{Inexact Newton Methods}
	
	A strong theory of ``Inexact'' Newton methods exists for the solution of systems of nonlinear equations when there are ``many'' variables. Inexact Newton steps are simply ``inexact'' solutions to the Newton system; that is, an inexact Newton step $\vec{s}^{IN}$ is any vector that satisfies
	\begin{equation}
		\label{InexactNewtonCondition}
		\lvert\lvert \vec{F}(\vec{x}) + (D\vec{F})(\vec{x})\vec{s}^{IN} \rvert\rvert
			\leq \delta \lvert\lvert \vec{F}(\vec{x}) \rvert\rvert
	\end{equation}
	for some fixed $\delta \in (0,1)$ \citep{Dembo82,Brown90,Eisenstat94,Eisenstat96,Pernice98}. The name ``truncated'' Newton method has also been used for the specific case when the inexactness comes from the use of iterative linear system solvers like \texttt{GMRES} \citep{Saad86,Walker88} or \texttt{BiCGSTAB} \citep{vanDerVorst92,Sleijpen93}. We focus on \texttt{GMRES}, a particularly simple yet strong iterative method for general linear systems that has been consistently used in the context of solving nonlinear systems \citep{Brown90}. 
	
	By appropriately choosing a sequence of $\delta$'s, the local asymptotic convergence rate of an inexact Newton's method can be fully quadratic \citep{Dembo82,Eisenstat94}. Of course, taking $\delta \to 0$ to achieve the quadratic convergence rate will also require increasingly burdensome computations of inexact Newton steps that satisfy increasingly strict inexact Newton conditions. On the other hand, $\delta$ can be chosen to be a constant if a linear locally asymptotic convergence rate is suitable \citep{Pernice98}. 
	
	Generally speaking there are three reasons to adopt the inexact perspective. First, direct methods like QR factorization may not be the most effective means to solve the Newton system when this system is large, because of computational burden and accumulation of roundoff errors. Instead, iterative solution methods are often used to solve linear systems with many variables; see, e.g. \cite{Trefethen97}. Second, iterative methods like \texttt{GMRES} require only matrix-vector products $(D\vec{F})(\vec{p})\vec{s}$ that can be approximated with finite directional derivatives \citep{Brown90,Pernice98}. Thus inexact Newton's methods can be ``matrix-free''; see Section \ref{ECSUBSEC:DirectionalFiniteDerivatives} below. Third, Newton steps often point in inaccurate directions when far from a solution \citep{Pernice98}. Thus solving for exact Newton steps may involve wasted effort, especially when there are many variables. 
	
	\texttt{matlab}'s \texttt{fsolve} function implements a related approach using the (preconditioned) Conjugate Gradient (\texttt{CG}) method applied to the normal equation for the Newton system, $(D\vec{F})(\vec{p})^\top(D\vec{F})(\vec{p}) \vec{s}^{IN} = - (D\vec{F})(\vec{p})^\top\vec{F}(\vec{p})$. Use of the normal equations is required because \texttt{CG} is applicable only to symmetric systems \citep{Trefethen97}. Note that this requires that the Jacobian $(D\vec{F})$ is explicitly available. Although this holds for price equilibrium problems under Mixed Logit models, it can be a significant restriction for general problems. By requiring products $(D\vec{F})(\vec{p})^\top \vec{h}$ in each step of the iterative linear solver, this approach also increases the work by $\order(NJ^2)$ \texttt{flops} where the solver takes $N$ steps. Finally, this approach can also be less accurate: using the normal equation squares the linear problem's condition number, and thus risks serious degradation in solution quality \citep{Trefethen97}. \cite{Pernice98} describe a similar approach using \texttt{BiCGSTAB}: the extension of \texttt{CG} to non-symmetric systems.
	
	
	
	\subsection{\texttt{GMRES}}
	
	The ``Generalized Minimum Residuals'' or \texttt{GMRES} method \citep{Saad86} solves a linear system $\vec{Ax} = \vec{b}$ by using the Arnoldi process to compute an orthonormal basis of the successive Krylov subspaces $\set{K}^{(n)}$ and then takes approximate solutions from those subspaces having least squares residuals. See \cite{Trefethen97} for a good introduction to Krylov methods in general, including the Arnoldi process and \texttt{GMRES}. In the $n\ith$ stage, \texttt{GMRES} ``factors'' $\vec{A}$ as $\vec{A}\vec{Q}^{(n)} = \vec{Q}^{(n+1)}\tilde{\vec{H}}^{(n)}$ where $\vec{Q}^{(n)} \in \R^{N \times n}$ is an orthonormal basis for $\set{K}^{(n)}$, $\vec{Q}^{(n+1)} \in \R^{N \times (n+1)}$ is an orthonormal basis for $\set{K}^{(n+1)} \supset \set{K}^{(n)}$, and $\tilde{\vec{H}}^{(n)} \in \R^{(n+1) \times n}$ is upper-Hessenberg. Any vector $\vec{x} \in \set{K}^{(n)} \subset \R^N$ can be written $\vec{x} = \vec{Q}^{(n)}\vec{y}$ for some $\vec{y} \in \R^n$ and thus the least-squares residual problem becomes
	\begin{equation*}
		\min_{ \vec{x} \in \set{K}^{(n)} } \norm{ \vec{As} - \vec{b} }_2
			= \min_{ \vec{y} \in \R^n } \norm{ \vec{A}\vec{Q}^{(n)}\vec{y} - \vec{b} }_2
			= \min_{ \vec{y} \in \R^n } \norm{ \tilde{\vec{H}}^{(n)}\vec{y} - (\vec{Q}^{(n+1)})^\top\vec{b} }_2. 
	\end{equation*}
	The orthonormal basis is typically chosen so that $(\vec{Q}^{(n+1)})^\top\vec{b} = \beta \vec{e}_1$ for some $\beta \in \R$, and hence the \texttt{GMRES} solution $\vec{x}^{(n)} = \vec{Q}^{(n)}\vec{y}$ where $\vec{y}$ solves $\min_{ \vec{q} \in \R^n } \norm{ \tilde{\vec{H}}^{(n)}\vec{y} - \beta\vec{e}_1 }_2$. This least squares problem can be solved using the QR factorization of $\tilde{\vec{H}}^{(n)}$. Furthermore this factorization can be efficiently updated in each iteration, instead of computed from scratch. Moreover the actual solution vector need not be formed until the residual is suitably small. 
	
	
	\subsubsection{Householder GMRES}
	
	We have implemented a variant of \texttt{GMRES} based on Householder transformations due to \cite{Walker88}; this is also the version implemented in \texttt{matlab}'s \texttt{gmres} code. We have verified that our implementation generates results matching \texttt{matlab}'s implementation. In this version of the \texttt{GMRES} process applied to the generic problem $\vec{Ax} = \vec{b}$, Householder reflectors $\vec{P}^{(n)} \in \R^{N \times N}$ are used to generate the orthonormal matrices
	\begin{equation*}
		\vec{Q}^{(n)} 
			= \vec{P}^{(1)} \dotsb \vec{P}^{(n)} \begin{bmatrix} \vec{I} \\ \vec{0} \end{bmatrix}
				\in \R^{N \times n}
			\quad\quad( \vec{I} \in \R^{n \times n}, \; \vec{0} \in \R^{(N - n) \times n})
	\end{equation*}
	satisfying
	\begin{equation*}
		\vec{A}\vec{Q}^{(n)} = \vec{P}^{(1)} \dotsb \vec{P}^{(n+1)} \vec{H}^{(n)} = \vec{Q}^{(n+1)} \tilde{\vec{H}}^{(n)}
	\end{equation*}
	where $\vec{H}^{(n)} \in \R^{N \times n}$ is
	\begin{equation*}
		\vec{H}^{(n)} = \begin{bmatrix} \tilde{\vec{H}}^{(n)} \\ \vec{0} \end{bmatrix}
	\end{equation*}
	for upper Hessenberg $\tilde{\vec{H}}^{(n)} \in \R^{(n+1)\times n}$ and $\vec{0}\in \R^{ ( N - n - 1 ) \times n}$. $\vec{P}^{(1)}$ is chosen to satisfy $\vec{P}^{(1)} \vec{b} = - \beta \vec{e}_1$ where $\beta = \sign( b_1 ) \norm{\vec{b}}_2$,  and hence $(\vec{Q}^{(n)})^\top\vec{b} = - \beta \vec{e}_1$. The $n\ith$ approximate solution $\vec{x}^{(n)}$ is taken to be $\vec{x}^{(n)} = \vec{Q}^{(n)} \vec{y}^{(n)}$ where $\vec{y}^{(n)} \in \R^n$ solves
	\begin{equation*}
		\min_{\vec{y} \in \R^n} \norm{ \tilde{\vec{H}}^{(n)} \vec{y} - \beta \vec{e}_1 }_2. 
	\end{equation*}
	Again these problems can be solved cheaply by updating QR factorizations with Givens rotations. Neither the solution vector nor the residual vector be formed until \texttt{GMRES} converges. An efficient implementation requires $\order(Jn)$ \texttt{flops} and a matrix multiply in the $n\ith$ iteration, so that taking $N$ iterations requires $\order(JN^2)$ of ``overhead'' in addition to the $\order(NJ^2)$ work required for the matrix multiplications (using the actual Jacobians). So long as $N < J$, using \texttt{GMRES} with the actual Jacobians is cheaper than solving for the actual Jacobian with QR. With small $N$, as we achieve using $\bsym{\eta}$ and $\bsym{\zeta}$, the savings is quite substantial. 
	
	We note the following formulae specific to the Newton system case. For $\vec{A} = (D\vec{F})(\vec{x})$ and $\vec{b} = - \vec{F}(\vec{x})$, $\beta = - \sign( F_1(\vec{x})  ) \norm{ \vec{F}(\vec{x}) }_2$ and $- \beta \vec{e}_1 = \vec{P}^{(1)} \vec{b} = - \vec{P}^{(1)} \vec{F}(\vec{x})$ so that
	\begin{equation*}
		\vec{P}^{(1)} \vec{F}(\vec{x}) = \beta \vec{e}_1 = - \sign( F_1(\vec{x})  ) \norm{ \vec{F}(\vec{x}) }_2 \vec{e}_1. 
	\end{equation*}
	Moreover, $\vec{P}^{(n)}\vec{e}_1 = \vec{e}_1$ for all $n > 1$ so that 
	\begin{equation*}
		(\vec{Q}^{(n)})^\top\vec{F}(\vec{x}) = - \sign( F_1(\vec{x})  ) \norm{ \vec{F}(\vec{x}) }_2 \vec{e}_1. 
	\end{equation*}
	
	
	\subsubsection{Preconditioning}
	\label{ECSUBSEC:GMRESPreconditioner}
	
	As is well known, preconditioning is key to the effectiveness of iterative linear solvers; see \cite{Golub96}. We have not found the linear systems in $\bsym{\eta}$-NM or $\bsym{\zeta}$-NM to need preconditioning. However we have found the preconditioned system 
	\begin{equation}
		\label{PreCondNewtonSystem}
		\bsym{\Lambda}(\vec{p})\inv(D\tilde{\nabla}\hat{\pi})(\vec{p})\vec{s}^{IN} 
			= - \bsym{\Lambda}(\vec{p})\inv(\tilde{\nabla}\hat{\pi})(\vec{p})
			= \vec{c} + \bsym{\zeta}(\vec{p}) - \vec{p}
	\end{equation}
	to be very necessary for rapid solution of the Newton system in CG-NM. This preconditioner is motivated by the following relationship of the Jacobian of $(\tilde{\nabla}\hat{\pi})$ to the Jacobian of $\bsym{\zeta}$ in equilibrium. 
	
	\begin{lemma}
		\label{DzetaDPIIdentity}
		$\vec{I} - (D\bsym{\zeta})(\vec{p}) = \bsym{\Lambda}(\vec{p})\inv(D\tilde{\nabla}\hat{\pi})(\vec{p})$ for any simultaneously stationary $\vec{p}$. 
	\end{lemma}
	\proof{}
		This follows from differentiating $(\tilde{\nabla}\hat{\pi})(\vec{p}) = \bsym{\Lambda}(\vec{p})( \vec{p} - \vec{c} - \bsym{\zeta}(\vec{p}) )$ via the product rule, recognizing that $\vec{p} - \vec{c} - \bsym{\zeta}(\vec{p}) = \vec{0}$ in equilibrium and $D[\vec{p} - \vec{c} - \bsym{\zeta}(\vec{p})] = \vec{I} - (D\bsym{\zeta})(\vec{p})$.
	\endproof
	
	In other words, Newton's methods applied to $\vec{F}_\pi(\vec{p})$ preconditioned as above ends up being essentially the same iteration as $\vec{F}_\zeta(\vec{p})$, close enough to equilibrium. 
	
	\texttt{GMRES}, if used successfully on this preconditioned system Eqn. (\ref{PreCondNewtonSystem}), will ensure that
	\begin{equation}
		\label{PreCondInexactNewtonCondition}
		\norm{ \bsym{\Lambda}(\vec{p})\inv(\tilde{\nabla}\hat{\pi})(\vec{p})
				+ \bsym{\Lambda}(\vec{p})\inv(D\tilde{\nabla}\hat{\pi})(\vec{p})\vec{s}^{IN} }
			\leq \delta^\prime \norm{ \bsym{\Lambda}(\vec{p})\inv(\tilde{\nabla}\hat{\pi})(\vec{p}) }
	\end{equation}
	for some $\delta^\prime$. This is distinct from the inexact Newton condition Eqn. (\ref{InexactNewtonCondition}). The following proposition gives modified tolerances for the preconditioned system to ensure satisfaction of the original system.
	\begin{proposition}
		Let $\delta > 0$ be given. If Eqn. (\ref{PreCondInexactNewtonCondition}) is satisfied with $\delta^\prime(\vec{p},\delta) \leq \delta$ given by
		\begin{equation}
			\label{PreCondTolerance}
			\delta^\prime(\vec{p},\delta)
				= \left( \frac{ \norm{(\tilde{\nabla}\hat{\pi})(\vec{p})}_2 }
						{ \max_j \left \{ \abs{\lambda_j(\vec{p})} \right \}
							\norm{ \bsym{\Lambda}(\vec{p})\inv(\tilde{\nabla}\hat{\pi})(\vec{p})}_2 }
					\right) \delta, 
		\end{equation}
		then Eqn. (\ref{InexactNewtonCondition}) is satisfied. 
	\end{proposition}
	This is a consequence of the following general result, which we state without proof.
	\begin{lemma}
		Let $\vec{b} \in \R^N$ and $\vec{A},\vec{M} \in \R^{N \times N}$ be nonsingular. Then
		\begin{equation}
			\label{RelNormBound}
			\frac{ \norm{ \vec{Ax} - \vec{b} } }{ \norm{ \vec{b} } }
				\leq \alpha \left( \frac{ \norm{ \vec{M}\inv\vec{Ax} - \vec{M}\inv\vec{b} } }{ \norm{ \vec{M}\inv\vec{b} } } \right)
		\end{equation}
		where $\alpha \in [1,\kappa(\vec{M})]$ is given by
		\begin{equation*}
			\alpha 
				= \frac{ \norm{ \vec{M} } \norm{ \vec{M}\inv\vec{b} } }{ \norm{ \vec{b} } }
				= \norm{ \vec{M} } \left \lvert \left \lvert \vec{M}\inv \left( \frac{ \vec{b} }{ \norm{ \vec{b} } } \right) \right \rvert \right \rvert . 
		\end{equation*}
	\end{lemma}
	This implies that
	\begin{equation*}
		\frac{ \norm{\vec{Ax} -\vec{b} } }{ \norm{\vec{b} } }
			\leq\delta
		\quad\quad\text{if}\quad\quad
		\frac{ \norm{ \vec{M}\inv\vec{Ax} - \vec{M}\inv\vec{b} } }{ \norm{ \vec{M}\inv\vec{b} } }
			\leq \frac{\delta}{\alpha}. 
	\end{equation*}
	Note that the preconditioned system must always be solved to a {\em stricter} tolerance than is desired for the un-preconditioned system using this bound. Additionally, computing $\alpha$ for a generic preconditioner $\vec{M}$ relies on the ability to compute $\norm{\vec{M}}$. 
	
	Eqn. (\ref{PreCondTolerance}) simply adopts the 2-norm and applies the formula \citep{Golub96} 
	\begin{equation*}
		\norm{\bsym{\Lambda}(\vec{p})}_2 
			= \sqrt{ \max_j \{ \abs{\lambda_j(\vec{p})}^2 \} }
			= \max_j \{ \abs{\lambda_j(\vec{p})} \}
	\end{equation*}
	
\begin{figure*}
		
			\begin{center}
				\includegraphics{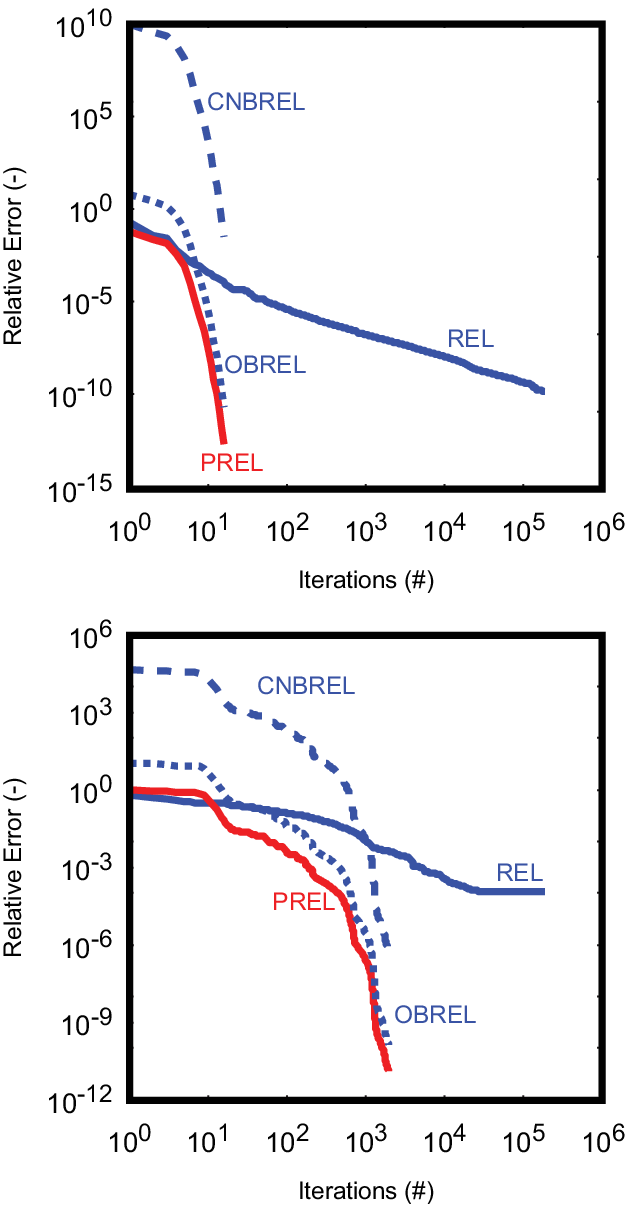}
				\caption[Relative error in computed solutions to the Newton system and its preconditioned form using \texttt{GMRES} under the \cite{Berry95} model.]{Relative error in computed solutions to the CG-NM Newton system and its preconditioned form using \texttt{GMRES} in the vehicle example from \cite{Morrow10} using the \cite{Berry95} model. On the top, prices are $\vec{p} = \vec{p}^* + 100 \bsym{\nu}$ where $\vec{p}^*$ are equilibrium prices and $\bsym{\nu} \in [-\vec{1},\vec{1}]$ is a sample from a uniformly distributed random vector. For this case $\kappa(\bsym{\Lambda}(\vec{p})) = 1.56 \times 10^{11}$ while the multiplier in Eqn. (\ref{RelNormBound}) is only 106.41. On the bottom, prices are $\vec{p} = 20,000\bsym{\nu} + 5,000$ where $\bsym{\nu}$ is a sample from a random vector uniformly distributed on $[\vec{0},\vec{1}]$. For this case $\kappa(\bsym{\Lambda}(\vec{p})) = 4.6 \times 10^{4}$ while the multiplier in Eqn. (\ref{RelNormBound}) is only 10.73. Abbreviations are as follows. REL: relative error in the Newton System; PREL: relative error in the pre-conditioned Newton System; OBREL: our bound, Eqn. (\ref{RelNormBound}), on the relative error in the Newton System as determined from the relative error in the preconditioned Newton system; CNBREL: condition number bound on the relative error in the Newton System as determined from the relative error in the preconditioned Newton system. }
				\label{FIG_BLP95_gmrestest}
			\end{center}
			
	\end{figure*}

	Eqn. (\ref{RelNormBound}) also implies that if Eqn. (\ref{PreCondInexactNewtonCondition}) holds with $\delta^\prime > 0$, then
	\begin{equation*}
		\frac{ \norm{ (\tilde{\nabla}\hat{\pi})(\vec{p}) + (D\tilde{\nabla}\hat{\pi})(\vec{p})\vec{s}^{IN} }_2 }
			{ \norm{ (\tilde{\nabla}\hat{\pi})(\vec{p}) }_2 }
			\leq \kappa_2(\bsym{\Lambda}(\vec{p})) \delta^\prime
	\end{equation*}
	where $\kappa_2(\bsym{\Lambda}(\vec{p})) = \norm{ \bsym{\Lambda}(\vec{p}) }_2\norm{ \bsym{\Lambda}(\vec{p})\inv }_2$ is the (2-norm) condition number of $\bsym{\Lambda}(\vec{p})$. This equation, while the more compact representation, can also be overly conservative as clearly illustrated in Fig. \ref{FIG_BLP95_gmrestest}. It is unlikely that $\kappa(\bsym{\Lambda}(\vec{p}))$ is a {\em tight} upper bound on the multiplier in Eqn. (\ref{PreCondTolerance}). In fact, the multiplier on $\delta$ depends only on the norm of $\bsym{\Lambda}(\vec{p})\inv\vec{x}$ at a single point on the surface of the unit sphere in $\R^J$ rather than $\norm{ \bsym{\Lambda}(\vec{p})\inv }_2$, the maximum norm of $\bsym{\Lambda}(\vec{p})\inv\vec{x}$ over this entire sphere. Our examples in Fig. \ref{FIG_BLP95_gmrestest} bear this out, having condition numbers many orders of magnitude larger than the multiplier in Eqn. (\ref{PreCondTolerance}). 
	
	The power of the preconditioning is that the preconditioned system Eqn. (\ref{PreCondInexactNewtonCondition}) appears to be solved to a relative error of $\delta^\prime(\vec{p},\delta)$ {\em much} faster than the original system can be solved to a relative error of $\delta$, even though $\delta^\prime(\vec{p},\delta) \leq \delta$. As can be seen in Fig. \ref{FIG_BLP95_gmrestest}, solving the preconditioned system to $\delta^\prime(\vec{p},\delta)$ can achieve a relative error in the original system below $\delta = 10^{-10}$ in roughly {\em four orders of magnitude fewer iterations} than solving the original system to this same relative error for prices near equilibrium. Away from equilibrium, \texttt{GMRES} may not be able to solve the original system to small relative errors like $10^{-6}$ at all. Thus using the original system would appear to slow, if not halt, an implementation of the inexact Newton's method. 
	
	\subsection{The \texttt{GMRES} Hookstep}
	
	Suitable modifications of each of the globalization strategies originally developed for ``exact'' Newton methods can be applied in the inexact context. \cite{Brown90} directly extend line search and a dogleg steps to \texttt{GMRES}-Newton methods. \cite{Eisenstat96} and \cite{Pernice98} apply a safeguarded backtracking line search to facilitate global convergence. More recently, \cite{Pawlowski06,Pawlowski08} have studied dogleg steps suitable for \texttt{GMRES}-Newton methods in some detail. Finally \cite{Viswanath07} has derived an elegant version of the hookstep method suitable for \texttt{GMRES}-Newton methods. In contrast with the hookstep approach for the ``exact'' Newton method with Jacobian $(D\vec{F})(\vec{p})$, Viswanath's approach requires computing the SVD only of a matrix whose size is determined by the number of iterations taken by \texttt{GMRES}. For reasonable applications of \texttt{GMRES}, this can be {\em far} less than the size of $(D\vec{F})(\vec{p})$ itself. For the examples in \cite{Morrow10}, the size difference is roughly two orders of magnitude: the \texttt{GMRES}-Newton hookstep worked with roughly $10 \times 10$ instead of $1,000 \times 1,000$ matrices. Thus, the \texttt{GMRES}-hookstep can accumulate a tremendous savings over an exact-Newton implementation of the hookstep method. Again, each of these approaches iterates until an acceptable step is found, and can, in principle, involve many additional evaluations of $\vec{F}$ or fail to find an acceptable step altogether. 
	
	Here we describe an implementation of the Levenberg-Marquardt method or ``hookstep'' \citep{Dennis96} suitable for \texttt{GMRES} as first suggested by \cite{Viswanath07}. See also \cite{Viswanath09a, Viswanath09b, Halcrow09}. First, we recall the basic structure of model trust region methods; see \cite[Chapter 6, Section 4]{Dennis96}. We then adopt this structure to the case of Krylov subspace methods, particularly \texttt{GMRES}. Again, see \cite{Morrow10b} for a more detailed discussion of this method. 
	
		
	\subsubsection{Model Trust Region Methods.}
	
	Trust region methods assume that for steps $\vec{s}$ satisfying $\norm{\vec{s}}_2 \leq \delta$, the function
	\begin{equation*}
		\hat{m}_{\vec{x}}(\vec{s})
			= \left( \frac{1}{2} \right) \norm{ \vec{F}(\vec{x}) }_2^2
				+ ( (D\vec{F})(\vec{x})^\top\vec{F}(\vec{x}) )^\top \vec{s}
				+ \left( \frac{1}{2} \right) \vec{s}^\top (D\vec{F})(\vec{x})^\top(D\vec{F})(\vec{x}) \vec{s}
	\end{equation*}
	is an accurate {\em local model} of $f(\vec{x}) = \norm{ \vec{F}(\vec{x}) }_2^2/2$ for suitably small steps. Note that $\hat{m}_{\vec{x}}$ is not the usual, quadratic model of $f$ derived from a Taylor series because $(D\vec{F})(\vec{x})^\top(D\vec{F})(\vec{x}) \neq (D \nabla f)(\vec{x})$ \cite[pg. 149]{Dennis96}. The idea is to solve 
	\begin{equation}
		\label{LocalModelProblem}
		\min_{ \norm{\vec{s}}_2 \leq \delta } \hat{m}_{\vec{x}}(\vec{s}). 
	\end{equation}
	The solution $\vec{s}_*$ is given as follows: take $\vec{s}_* = \vec{s}^N = - (D\vec{F})(\vec{x})\inv\vec{F}(\vec{x})$ if $\norm{ \vec{s}^N }_2 \leq \delta$; if $\norm{ \vec{s}^N }_2 > \delta$, take $\vec{s}_*  = \vec{s}(\mu_*)$ where
	\begin{equation*}
		\vec{s}(\mu) = - \big( (D\vec{F})(\vec{x})^\top(D\vec{F})(\vec{x}) + \mu \vec{I} \big)\inv(D\vec{F})(\vec{x})^\top\vec{F}(\vec{x})
	\end{equation*}
	and $\mu_* > 0$ is the unique $\mu > 0$ such that $\norm{ \vec{s}(\mu) }_2 = \delta$. These follow from the standard optimality conditions, or rather that the gradient $(\nabla\hat{m}_{\vec{x}})(\vec{s})$ must lie in the negative normal cone to $\bar{\mathbb{B}}_\delta(\vec{0}) = \{ \vec{y} \in \R^N : \norm{\vec{y}}_2 \leq \delta \}$ at $\vec{x}$ \citep{Clarke75}; see \citep[Lemma 6.4.1, pg. 131]{Dennis96}.
	
	Solving the problem above exactly generates the Levenberg-Marquardt method \citep{Levenberg44,Marquardt63} or ``hookstep.'' By computing the SVD of $(D\vec{F})(\vec{x})$ we can easily solve for $\vec{s}(\mu)$ when $\norm{\vec{s}^N}_2 > \delta$ \citep{Dennis96}; see \cite[Section 12.1, pgs. 580-583]{Golub96} for closely related results. Let $(D\vec{F})(\vec{x}) = \vec{U}\bsym{\Sigma}\vec{V}^\top$. We can then set $\vec{s}(\mu) = \vec{V}\vec{r}(\mu)$ where 
	\begin{equation*}
		\vec{r}(\mu) = - ( \bsym{\Sigma}^2 + \mu\vec{I} ) \inv \bsym{\Sigma}\vec{U}^\top \vec{F}(\vec{x}). 
	\end{equation*}
	A simple single-dimensional iteration can then be used to solve for the unique $\mu_*$ such that $\norm{\vec{s}(\mu_*)}_2 = \delta$. \cite{Morrow10b} derives two globally convergent methods for this task using Newton's method and a nonlinear local model \citep{Dennis96}. The difficulty here is computing the SVD of $(D\vec{F})(\vec{x})$, requiring $\order(J^3)$ \texttt{flops} \citep[Chapter 5, pg. 254]{Golub96}. 
	
	The step $\vec{s}_*$ computed by either approach is {\em acceptable} if it generates sufficient decrease in the squared 2-norm of $\vec{F}$. Specifically, fix $\rho \in (0,1)$, $\alpha > 1$, and $\beta_2 \leq \beta_1 < 1$. If
	\begin{equation*}
		\norm{ \vec{F}(\vec{x}) }_2^2 
				- \norm{ \vec{F}(\vec{x}+\vec{s}_*) }_2^2
			\geq \rho ( \norm{ (D\vec{F})(\vec{x}) }_2^2 - \norm{ \vec{F}(\vec{x}) + (D\vec{F})(\vec{x})\vec{s}_* }_2^2 )
	\end{equation*}
	then $\vec{p} \leftarrow \vec{p} + \vec{s}_*$ and a the step length bound is expanded to $[\delta,\alpha\delta]$ for the next iteration. Otherwise, $\delta$ is chosen from $[\beta_1\delta,\beta_2\delta]$ and the corresponding $\vec{s}_*$ is computed. While this process of specifying an acceptable $\vec{s}_*$ is iterative, much of the work required to build a trial step does not need to be repeated. Specifically the SVD required for the hookstep does not change (so long as it was computed in a previous iteration) while in the doglep step the Newton and Cauchy steps remain the same. However every time the step size bound is decreased $\vec{F}$ must be re-evaluated at the new trial step, with a computational burden equivalent to taking a fixed-point step. 
	
	
	\subsubsection{Model Trust Region Methods on a Subspace}
	
	A Krylov method for solving $(D\vec{F})(\vec{x})\vec{s}^N = - \vec{F}(\vec{x})$ builds approximate solutions in the successive Krylov subspaces $\set{K}^{(n)}$. This has the effect of further constraining the local model problem (\ref{LocalModelProblem}) to
	\begin{equation}
		\min_{ \vec{s} \in \set{K}^{(n)} , \; \norm{ \vec{s} }_2 \leq \delta } \hat{m}_{\vec{x}}(\vec{s}). 
	\end{equation}
	
	For any $\vec{Q} \in \R^{J \times n}$ with orthonormal columns (generated by \texttt{GMRES} or not) we can set $\hat{m}_{\vec{x},\vec{Q}}(\vec{y}) = \hat{m}_{\vec{x}}(\vec{Q}\vec{y})$ and restrict attention to the trust region problem $\min_{ \norm{ \vec{y} }_2 \leq \delta } \hat{m}_{\vec{x},\vec{Q}}(\vec{y})$. See \cite[pgs. 149-150]{Brown90}. The first-order conditions for this problem are equivalent to either
	\begin{itemize}
		\item[(i)] $(\nabla \hat{m}_{\vec{x},\vec{Q}})(\vec{y}) = \vec{0}$ and $\norm{\vec{y}}_2 \leq \delta$
		\item[(ii)] or $(\nabla \hat{m}_{\vec{x},\vec{Q}})(\vec{y}) + \mu \vec{y} = \vec{0}$ for $\norm{\vec{y}}_2 = \delta$ and some $\mu > 0$. 
	\end{itemize}
	By the definition of $\hat{m}_{\vec{x},\vec{Q}}$, (i) implies
	\begin{equation*}
		 \vec{Q}^\top(D\vec{F})(\vec{x})^\top(D\vec{F})(\vec{x})\vec{Q}\vec{y}
			+ \vec{Q}^\top(D\vec{F})(\vec{x})^\top\vec{F}(\vec{x})
				= \vec{0}
	\end{equation*}
	and (ii) implies
	\begin{equation*}
		\left( \vec{Q}^\top (D\vec{F})(\vec{x})^\top(D\vec{F})(\vec{x})\vec{Q}
					 + \mu\vec{I} \right) \vec{y}
						+ \vec{Q}^\top(D\vec{F})(\vec{x})^\top\vec{F}(\vec{x})
			= \vec{0}. 
	\end{equation*}
	Note that these are square problems that can be solved exactly. 
	
	
	\subsubsection{The \texttt{GMRES}-Newton Hookstep}
	
	Using \texttt{GMRES} {\em started at zero}, $(D\vec{F})(\vec{x})\vec{Q}^{(n)} = \vec{Q}^{(n+1)} \tilde{\vec{H}}^{(n)}$ and $( \vec{Q}^{(n+1)} )^\top\vec{F}(\vec{x}) = - \sign( F_1(\vec{x}) ) \norm{ \vec{F}(\vec{x}) }_2 \vec{e}_1$. Thus we consider the family of $n \times n$ linear systems
	\begin{align*}
		&(\vec{Q}^{(n)})^\top (D\vec{F})(\vec{x})^\top(D\vec{F})(\vec{x})\vec{Q}^{(n)}\vec{q} 
					 + \mu \vec{q}
						+ (\vec{Q}^{(n)})^\top(D\vec{F})(\vec{x})^\top\vec{F}(\vec{x}) \\
		&\quad\quad\quad\quad\quad\quad\quad\quad
			= \big( (\tilde{\vec{H}}^{(n)})^\top \tilde{\vec{H}}^{(n)} + \mu \vec{I} \big) \vec{q} 
						- \sign( F_1(\vec{x}) ) \norm{ \vec{F}(\vec{x}) }_2 ( \tilde{\vec{H}}^{(n)})^\top \vec{e}_1 = \vec{0}
	\end{align*}
	defined for all $\mu \geq 0$. 
	
	By computing the (``thin'') Singular Value Decomposition of $\tilde{\vec{H}}^{(n)}$, $\tilde{\vec{H}}^{(n)} = \tilde{\vec{U}}\bsym{\Sigma}\vec{V}^\top$ where $\tilde{\vec{U}} \in \R^{(n+1) \times n}$, $\vec{V} \in \R^{n \times n}$, and $\bsym{\Sigma} \in \R^{n \times n}$, we can easily solve each such problem. See \cite[Section 12.1, pgs. 580-583]{Golub96} for closely related results. Particularly, 
	\begin{align*}
		& ( ( \tilde{\vec{H}}^{(n)} )^\top\tilde{\vec{H}}^{(n)} + \mu \vec{I} ) \vec{q} - \sign( F_1(\vec{x}) )\norm{\vec{F}(\vec{x})}_2 ( \tilde{\vec{H}}^{(n)} )^\top \vec{e}_1 = \vec{0}
	\end{align*}
	is solved by $\vec{q}(\mu) = \vec{V}\vec{r}(\mu)$ where 
	\begin{equation*}
		\vec{r}(\mu) 
			= \sign( F_1(\vec{x}) ) \norm{\vec{F}(\vec{x})}_2 ( \bsym{\Sigma}^2 + \mu \vec{I} )\inv \bsym{\Sigma} \tilde{\vec{U}} ^\top \vec{e}_1. 
	\end{equation*}
	Because the diagonal elements of $\bsym{\Sigma}^2$ are positive, $\vec{r}(\mu)$ is well defined for all $\mu \geq 0$. Note also that we only need the first row of $\vec{U}$, but all of $\vec{V}$, to compute $\vec{q}(\mu)$. 
	
	In particular, $\vec{q}(0) = \sign( F_1(\vec{x}) ) \norm{\vec{F}(\vec{x})}_2 \vec{V} \bsym{\Sigma}\inv \vec{U}^\top \vec{e}_1$. Invoking the {\em full} SVD of $\tilde{\vec{H}}^{(n)}$, 
	\begin{align*}
		\tilde{\vec{H}}^{(n)}
			= \begin{bmatrix} \tilde{\vec{U}} & \vec{u}_{n+1} \end{bmatrix}
				\begin{bmatrix} \bsym{\Sigma} \\ \vec{0}^\top \end{bmatrix}
				\vec{V}^\top
	\end{align*}
	for some $\vec{u}_{n+1} \perp \mathrm{span} \{ \vec{u}_i \}_{i=1}^n$, we can write
	\begin{align*}
		\norm{ \tilde{\vec{H}}^{(n)}\vec{q} - \sign( F_1(\vec{x}) ) \norm{\vec{F}(\vec{x})}_2 \vec{e}_1 }_2
			&= \left\lvert\left\lvert \begin{bmatrix} \tilde{\bsym{\Sigma}}\vec{V}^\top\vec{q}  \\ 0 \end{bmatrix}
						- \sign( F_1(\vec{x}) ) \norm{\vec{F}(\vec{x})}_2 \begin{bmatrix} \tilde{\vec{U}}^\top \vec{e}_1 \\ u_{1,n+1} \end{bmatrix} \right\rvert\right\rvert_2. 
	\end{align*}
	We thus see that $\vec{q}(0)$ solves the $(n+1) \times n$ \texttt{GMRES} least squares problem
	\begin{equation*}
		\min_{\vec{q}} \norm{ \vec{H}^{(n+1,n)}\vec{q} - \sign( F_1(\vec{x}) ) \norm{\vec{F}(\vec{x})}_2 \vec{e}_1 }_2. 
	\end{equation*}
	with residual $\abs{u_{1,n+1}} \norm{ \vec{F}(\vec{x}) }_2$. $\abs{u_{1,n+1}}$ is unique: First, note that $\vec{u}_{n+1}$ is a unit vector in the span of a single vector, say $\vec{v}$, that is orthogonal to the span of the columns of $\tilde{\vec{U}}$. There are only two unit vectors in this span, specifically $\pm \vec{v} / \norm{\vec{v}}_2$, and thus $\vec{u}_{n+1} \in \{ \pm \vec{v} / \norm{\vec{v}}_2 \}$. Thus $\abs{u_{1,n+1}} \in \abs{ \pm v_1 / \norm{ \vec{v} }_2 } = \abs{v_1} / \norm{ \vec{v} }_2$.
	
	It is also easy to see that
	\begin{align*}
		\vec{F}(\vec{x})^\top(D\vec{F})(\vec{x})\vec{s}^{(n)}(\mu)
			&= \vec{F}(\vec{x})^\top(D\vec{F})(\vec{x})\vec{Q}^{(n)}\vec{q}^{(n)}(\mu) \\
			&= \left( \big( \vec{Q}^{(n+1)} \big)^\top\vec{F}(\vec{x}) \right)^\top \tilde{\vec{H}}^{(n)}\vec{q}^{(n)}(\mu) \\
			&= - \beta^2 \left( \bsym{\nu}_1^\top \vec{D}(\mu) \bsym{\nu}_1 \right) \\
			&= - \norm{ \vec{F}(\vec{x}) }_2^2  \left( \bsym{\nu}_1^\top \vec{D}(\mu) \bsym{\nu}_1 \right) < 0
	\end{align*}
	where $\bsym{\nu}_1$ is the first row of $\tilde{\vec{U}}$ and $\vec{D}(\mu) = \diag( d_1(\mu) , \dotsc , d_n(\mu) )$ for $d_i(\mu) = \sigma_i^2 / ( \sigma_i^2 + \mu )$. That is, {\em the Householder \texttt{GMRES}-Newton Hookstep always lies in a descent direction for the globalizing objective} $f(\vec{x}) = \norm{ \vec{F}(\vec{x}) }_2^2/2$. 
	

	\subsubsection{Directional Finite Differences}
	\label{ECSUBSEC:DirectionalFiniteDerivatives}
	
	Recall that one advantage to using an iterative solver like \texttt{GMRES} to solve the Newton system is that only products of the type $(D\vec{F})(\vec{p})\vec{s}$ will be required to solve the Newton system for $\vec{F}$ at $\vec{p}$ \citep{Brown90,Pernice98}. Such products can be approximated by a {\em single} additional evaluation of $\vec{F}$ in a ``directional'' finite difference \citep{Brown90,Pernice98}. For example, the first-order formula
	\begin{equation*}
		(D\vec{F})(\vec{x})\vec{s} \approx h\inv \big( \vec{F}(\vec{x} + h\vec{s}) - \vec{F}(\vec{x}) \big), 
	\end{equation*}
	requires only a single additional evaluation of $\vec{F}$ per (approximate) evaluation of $(D\vec{F})(\vec{x})\vec{s}$. Higher-order formulae requiring 2 and 4 additional evaluations of $\vec{F}$ are easy to derive; see \cite{Pernice98}. In their implementation of the \texttt{GMRES} method in the context of an inexact Newton method, \cite{Pernice98} only use higher order finite-differencing formulas at restarts. \cite{Brown90} provide a practical formula for computing an appropriate value of $h$. 
	
	Since directional finite derivatives must be repeated {\em at each step} of iterative linear solvers, each step of an iterative Newton system solver using directional finite differences could be at least as expensive as a $\bsym{\zeta}$-FPI step. That is, if an iterative solver should take 100 steps to compute an inexact Newton step having small enough residual to satisfy the inexact Newton condition, then we could have equivalently taken 100, 200, and 400 $\bsym{\zeta}$-FPI steps with the first, second, and fourth order formulae available in \cite{Pernice98}. In our examples, using \texttt{GMRES} regularly solves the $\bsym{\eta}$-NM and $\bsym{\zeta}$-NM Newton systems in approximately 10 steps. This implies that each $\bsym{\eta}$-NM and $\bsym{\zeta}$-NM step is roughly equivalent to $10$ $\bsym{\zeta}$-FPI steps. 
	
	In the Newton context, whether efficiency is ultimately gained by using directional finite differences instead of computing the Jacobian matrices and using standard matrix-vector products depends on the number of steps taken by the iterative linear solver. If \texttt{GMRES} takes $N \in \N$ iterations to find an inexact Newton step for $\vec{F}$, computing and using the Jacobian requires $\order((S+N)J^2)$ \texttt{flops} while using directional finite differences requires $\order(SN\sum_{f=1}^FJ_f^2)$ \texttt{flops}. 
	
\begin{figure}
		\includegraphics[height=7.5in]{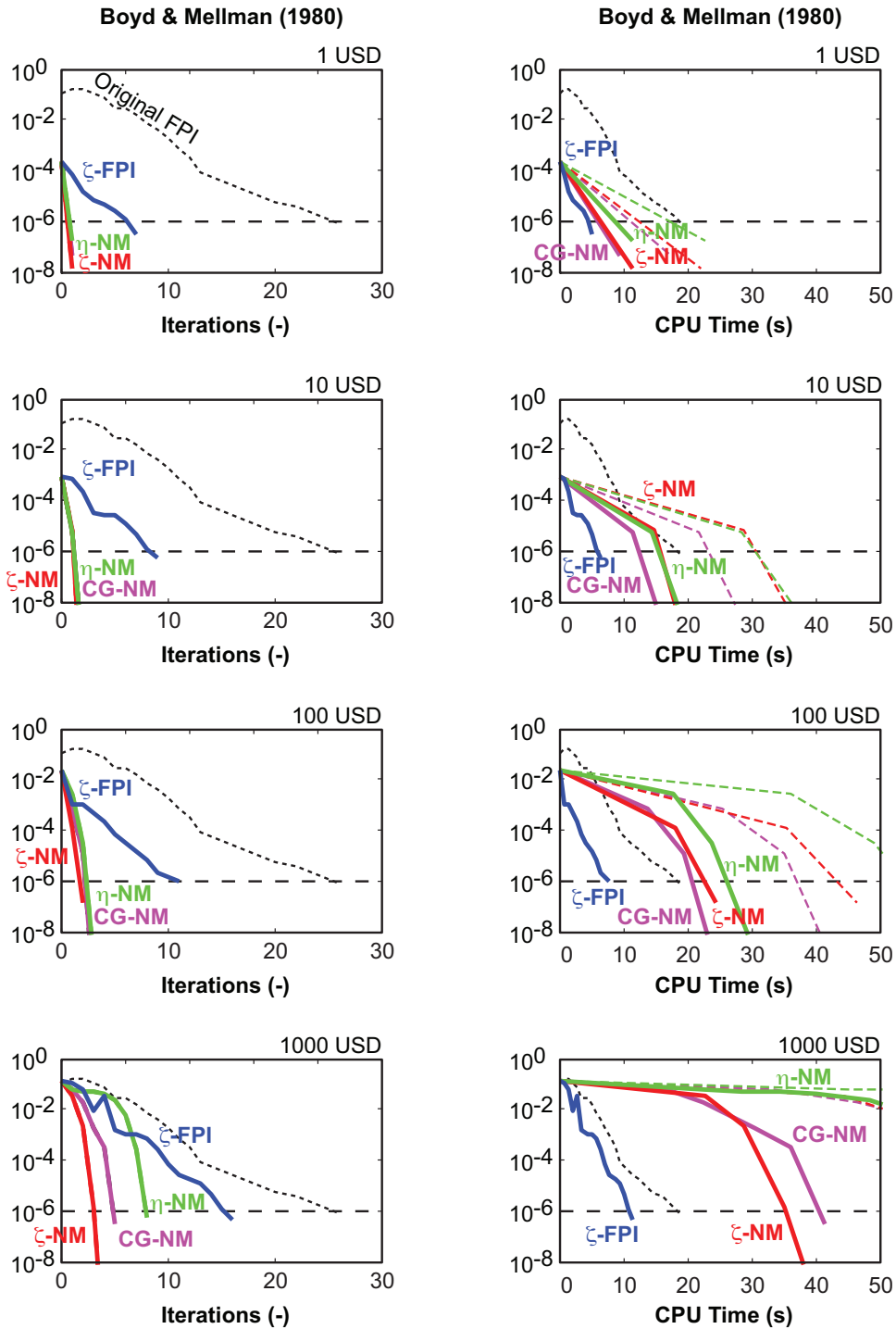}
		\caption{Typical convergence curves for perturbation trials under the \cite{Boyd80} model using both analytical and directional finite difference Jacobians. See also Fig. \ref{FIG_BM80_BestCasePert}. Convergence curves for analytical Jacobian are drawn with solid lines, whereas convergence curves for directional finite differences are drawn with dashed lines of the same color.}
		\label{FIG_BM80_PERTRIAL_C}
	\end{figure}

	
	We have observed that for $\bsym{\eta}$-NM and $\bsym{\zeta}$-NM, using the actual Jacobian takes roughly {\em half} the computation time than using directional finite differences, even though \texttt{GMRES} converges in very few iterations ($N \sim 10$). Fig. \ref{FIG_BM80_PERTRIAL_C} plots the sample trials for the \cite{Boyd80} model provided in \cite{Morrow10} using both analytical Jacobians and directional finite differences. First note that the $\bsym{\zeta}$-FPI regularly takes about 1 s per iteration. For $\kappa = 1$ USD, the single-step convergence of the \texttt{GMRES}-Newton Hookstep method translates into about 10 $\bsym{\zeta}$-FPI steps, or about $10$ s. Because \texttt{GMRES} itself requires some small overhead ($\order(Jn)$ in the $n\ith$ step), this is a somewhat reasonable estimate of the work required. Two \texttt{GMRES}-Newton steps are required with $\kappa = 10$ USD and we would expect about $20$ s, a somewhat less sound estimate of the time required. Three \texttt{GMRES}-Newton steps are required with $\kappa = 100$ USD, leading us to expect about $20$ s, a further less sound estimate of the time required. (These observations can be matched with an asymptotic analysis of the work required.) Note also that the $\bsym{\eta}$-NM has the greatest increase in time as a consequence of using the directional finite differences. This is a consequence of having to repeat block QR factorizations when evaluating $\bsym{\eta}$ at different points, while evaluating $(D\bsym{\eta})$ requires only a single factorization. 
	
\begin{figure}
		\includegraphics[height=7.5in]{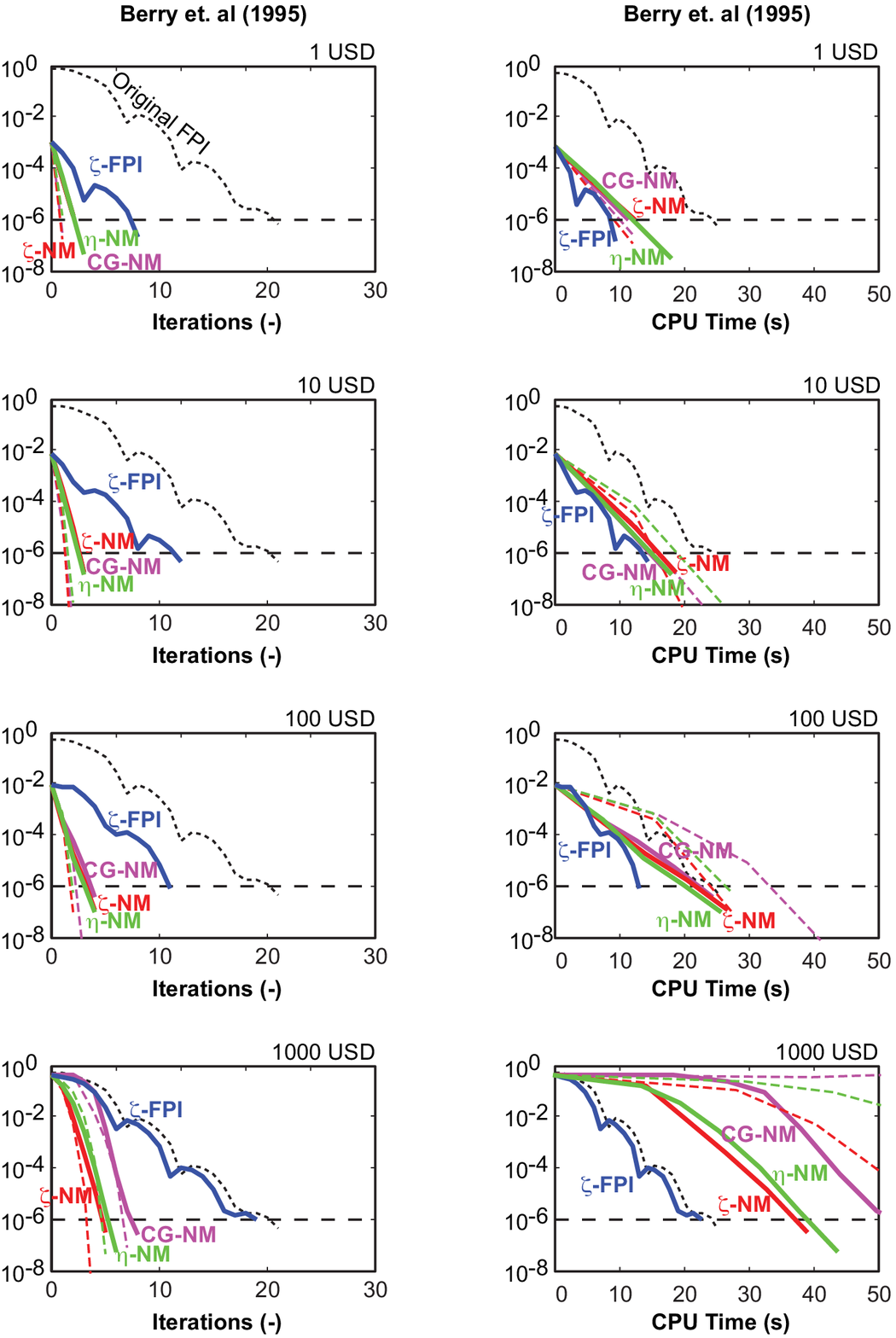}
		\caption{Typical convergence curves for perturbation trials under the \cite{Boyd80} model using both analytical and directional finite difference Jacobians. See also Fig. \ref{FIG_BLP95_BestCasePert}. Convergence curves for analytical Jacobian are drawn with solid lines, whereas convergence curves for directional finite differences are drawn with dashed lines of the same color.}
		\label{FIG_BLP95_PERTRIAL_C}
	\end{figure}

	Fig. \ref{FIG_BLP95_PERTRIAL_C} plots the sample trials for the \cite{Berry95} model provided in \cite{Morrow10} using both analytical Jacobians and directional finite differences. Interestingly, in this case use of the directional finite differences appears to generate a convergence {\em rate} improvement. Otherwise, the story remains much the same as that discussed above for the \cite{Boyd80} model. 
	
%% file: othermethods.tex

\section{Other Methods}
\label{ECSEC:OtherMethods}

	
	\subsection{Variational Methods} 
	\label{ECSUBSEC:VariationalMethods}
	
	Equilibrium problems are commonly formulated as variational inequalities or complementarity problems \citep{Harker90,Ferris97}. To be nontrivially distinct from nonlinear equations, such formulations require restricting the variables to a proper, convex subset of $\R^J$. When $\varsigma_* < \infty$ there is an appropriate variational formulation of the equilibrium pricing problem:
	\begin{equation}
		\label{VI}
		\text{find}\quad
		\vec{p} \in [0,\varsigma_*]^J
		\quad\text{such that}\quad
		(\tilde{\nabla}\hat{\pi})(\vec{p})^\top( \vec{p} - \vec{q} ) \geq 0
		\quad\text{for all}\quad
		\vec{q} \in [0,\varsigma_*]^J. 
	\end{equation}
	
	
	\subsubsection{The VI formulation is poorly posed} 
	
	Unfortunately, the Variational Inequality (\ref{VI}) is poorly posed when the derivatives of profit vanish as prices approach $\varsigma_* < \infty$. There are two specific issues with Eqn. (\ref{VI}) in this case. First, $\varsigma_*\vec{1} \in \set{P}^J$ is {\em always} a solution but {\em never} an equilibrium when profits vanish as all prices approach $\varsigma_*$; see Section \ref{ECSUBSEC:Profits} and Lemma \ref{InfinitySolvesVI}. Second, Eqn. (\ref{VI}) can be solved by {\em any} equilibrium of {\em any} differentiated product market model constructed with a {\em subset} of the products offered (Prop. \ref{SubEquilibriaSolveVI}). Equilibria of such ``sub-problems'' are {\em not} necessarily equilibria of the original problem, as demonstrated in Example \ref{VIExample} below. This issue with Eqn. (\ref{VI}) is, in fact, equivalent to the problem with CG-NM discussed in Section \ref{SUBSEC:CG-NM}. 
	
	These issues imply that variational methods can compute many ``spurious'' solutions. If an equilibrium problem and all its sub-problems have unique equilibria with all prices less than $\varsigma_*$, Eqn. (\ref{VI}) has $2^J$ solutions that might be recovered by a global method such as \texttt{PATH} \citep{Ralph94,Dirkse95}. However, only {\em one} of these solutions is an equilibrium of the original problem, by assumption. A simple example demonstrates this phenomenon. 
	
	\begin{example}
		\label{VIExample}
		Consider a monopoly with two products produced at the same unit cost $c$. Demand is given by a simple Logit model with product-specific utility functions $u_j(p_j) = \alpha \log( \varsigma - p_j ) + v_j$ for $j \in \{1,2\}$, where $\varsigma \in ( c , \infty )$, $v_1,v_2 \in \R$, and $\vartheta > -\infty$. The firm has unique profit-maximizing prices $(p_1^*,p_2^*)$. Furthermore $p_1^*,p_2^* < \varsigma$, and $(p_1^*,p_2^*)$ is the unique fixed-point of the map $\vec{c} + \bsym{\zeta}(\cdot)$ on all of $\set{P}^2$. 
		
		However the variational inequality formulation contains four distinct solutions, only one of which is profit-maximizing. These four solutions are $(p_1^*,p_2^*)$, $(\varsigma,\varsigma)$, $(q_1^*,\varsigma)$, and $(\varsigma,q_2^*)$, where $q_j^* < \varsigma$ for $j \in \{1,2\}$ are the unique profit-maximizing prices that exist should the firm offer only product $1$ or $2$. Only the first solution, $(p_1^*,p_2^*)$, is profit-maximizing.
	\end{example}
	
	\proof{} We complete the details of  Example \ref{VIExample}. 
	
		Consider a monopoly with two products produced at the same unit cost ($c = c_1 = c_2 > 0$), $\vartheta > -\infty$, and simple Logit model with utility
		\begin{equation*}
			u_1(p_1) = \alpha \log( \varsigma - p_1 ) + v_1
			\quad\text{and}\quad
			u_2(p_2) = \alpha \log( \varsigma - p_2 ) + v_2
		\end{equation*}
		for some fixed $\varsigma \in ( c , \infty )$, $\alpha > 1$, and arbitrary $v_1,v_2 \in \R$. Let $p_2 \leq \varsigma$, and observe that
		\begin{equation*}
			\lim_{ p_1 \uparrow \varsigma } \Big( p_1 - c - \zeta_1(p_1,p_2) \Big)
				= \varsigma - c - P_2(\varsigma,p_2)( p_2 - c )
				= (\varsigma - c) \left[ 1 - P_2(\varsigma,p_2)\left( \frac{ p_2 - c }{\varsigma - c} \right) \right]. 
		\end{equation*}
		Since $p_2 \leq \varsigma$ and $P_2(p_1,p_2) < 1$ for all $p_1,p_2$, we have $\lim_{ p_1 \uparrow \varsigma } ( p_1 - c - \zeta_1(p_1,p_2) ) > 0$. Thus $(D_1\hat{\pi})(p_1,p_2) < 0$ for all $p_1$ sufficiently close to $\varsigma$. A similar argument can be made for $(D_2\hat{\pi})(p_1,p_2)$. 
		
		Note also that this proves that $\varsigma + \epsilon > c + \zeta_1(\varsigma + \epsilon,p_2)$ for any $\epsilon \geq 0$ and $p_2$, where $\zeta_1$ is the extended map. A similar result holds for $\zeta_2$, instead of $\zeta_1$. Thus no $(p_1,p_2)$ outside of $(0,\varsigma)$ is fixed for the extended map $\vec{c} + \bsym{\zeta}(\vec{p})$.  
		
		We now prove that there exists a {\em unique} pair of profit-maximizing prices $\vec{p}^* = (p_1^*,p_2^*) \in (0,\varsigma)^2$. Since
		\begin{equation*}
			\lim_{ p_j \uparrow \varsigma } \Big( p_j - c - \zeta_j(p_1,p_2) \Big) < \infty
		\end{equation*}
		for $j \in \{1,2\}$, $\bsym{\zeta} = (\zeta_1,\zeta_2)$ is bounded and continuous on $\set{P}^2$. By Brower's fixed-point theorem, there exists a stationary point $\vec{p}^* = (p_1^*,p_2^*)$. Both prices must both be less than $\varsigma$, since profits decrease for all prices sufficiently close to $\varsigma$. We now show that these prices are also unique, borrowing a technique from \cite{Morrow08a}. 
		
		The first step is to prove that $(D\nabla\hat{\pi})(\vec{p}^*)$ is negative definite at any stationary $\vec{p}^*$. Note that $(D\nabla\hat{\pi})(\vec{p}^*) = \bsym{\Lambda}(\vec{p}^*)(\vec{I} - (D\bsym{\zeta})(\vec{p}^*))$; this relationship is valid for Mixed Logit models with multiple firms as well. Furthermore $\zeta_j(\vec{p}) = \hat{\pi}(\vec{p}) - (Dw_k)(p_k)\inv$ for {\em any} simple Logit model and any number of products. Hence
		\begin{equation*}
			(D_k\zeta_j)(\vec{p})
				= (D_k\hat{\pi})(\vec{p}) + \delta_{j,k} \left( \frac{(D^2w_k)(p_k)}{(Dw_k)(p_k)^2} \right), 
		\end{equation*}
		and $\vec{I} - (D\bsym{\zeta})(\vec{p}^*)$ is a diagonal matrix with elements
		\begin{equation*}
			1 - \frac{(D^2w_k)(p_k)}{(Dw_k)(p_k)^2}. 
		\end{equation*}
		In the case of this example,
		\begin{equation*}
			1 - \frac{(D^2w_1)(p_1)}{(Dw_1)(p_1)^2}
				= 1 - \frac{(D^2w_2)(p_2)}{(Dw_2)(p_2)^2}
				= 1 + \frac{1}{\alpha}
				> 0. 
		\end{equation*}
		Thus $(D\nabla\hat{\pi})(\vec{p}^*)$ is negative definite at any stationary point, and any stationary point maximizes profits. 
		
		The next step is to prove that the existence of only maximizers of profits proves that there is a unique pair of profit-maximizing prices. \cite{Morrow08a} accomplish this with an application of the Poincare-Hopf theorem \citep{Milnor65}, as follows. Consider $- \hat{\pi}(\vec{p})$. This function is minimized at any stationary $\vec{p}^* = (p_1,p_2)$, and thus the gradient vector field $-(\nabla\hat{\pi})(\vec{p})$ has index $1$ at any stationary point $\vec{p}^*$ \citep{Milnor65}. Note also that 
		\begin{align*}
			\mathrm{sign} \{ -(D_j\hat{\pi})(p_1,p_2) \}
				&= \mathrm{sign} \left\{ p_j - c - \hat{\pi}(p_1,p_2) - \frac{\varsigma - p_j}{\alpha} \right\}
		\end{align*}
		for $j \in \{1,2\}$. This equation shows that the gradient vector field $-(\nabla\hat{\pi})(\vec{p})$ points outward on the boundary of the compact, convex set $[c,\varsigma]^2$, as can be checked. Thus the Poincare-Hopf theorem states that the sum of the indices of the critical (stationary) points equals one, the Euler characteristic of $[c,\varsigma]^2$. Since the index of any critical (stationary) point of $-(\nabla\hat{\pi})(\vec{p})$ is one, there can only be one stationary point. 
		
		Using similar arguments, we see that the sub-problems formed by offering product 1 or product 2 alone also have unique profit-maximizing prices $q_1^*$ and $q_2^*$, respectively. Because $v_1$ and $v_2$ may be distinct, these prices need not be the same. 
		
		We have claimed that variational formulation of this problem has four solutions, only one of which is an equilibrium. Indeed, these four solutions are $(p_1^*,p_2^*)$, $(q_1^*,\varsigma)$, $(\varsigma,q_2^*)$, and $( \varsigma , \varsigma )$ but, as shown above, only $(p_1^*,p_2^*)$ is an equilibrium. While this follows from Props. \ref{InfinitySolvesVI} and \ref{SubEquilibriaSolveVI} above, we prove it directly here. Of course, $(p_1^*,p_2^*)$ is a solution since $(\nabla\hat{\pi})(p_1^*,p_2^*) = (0,0)$. Since
		\begin{equation*}
			\lim_{p_j \uparrow \varsigma} \lambda_j(p_1,p_2)
				= \lim_{p_j \uparrow \varsigma}
					\left[ \left( \frac{\varsigma - p_j}{\alpha} \right) P_j(p_1,p_2) \right]
				= 0
		\end{equation*}
		for $j \in \{1,2\}$, $\lim_{p_j \uparrow \varsigma} (D_j\hat{\pi})(p_1,p_2) = 0$ (i.e., Assumption \ref{LVA} holds). Thus $(\nabla\hat{\pi})(\varsigma,\varsigma) = (0,0)$, and the variational inequality is satisfied at $(\varsigma,\varsigma)$. Furthermore, 
		\begin{equation*}
			(D_1\hat{\pi})(\varsigma,p_2)(\varsigma - q_1)
					+ (D_2\hat{\pi})(\varsigma,p_2)(p_2 - q_2)
				= (D_2\hat{\pi})(\varsigma,p_2)(p_2 - q_2)
		\end{equation*}
		and thus $(\varsigma,q_2^*)$ is also a solution to the variational inequality. Similarly, $(q_1^*, \varsigma)$ is also a solution. This completes the proof. 
	\endproof
		
	Example \ref{VIExample} is easily generalized to include $J > 2$ products and a variational inequality with $2^J$ solutions. One of these solutions is the unique vector of profit-maximizing prices for the original problem, one is $\varsigma\vec{1} \in \set{P}^J$ and is not profit-maximizing for any sub-problem, and the rest are profit-maximizing for some sub-problem but not profit-maximizing for the original problem. 
	
	This property of variational formulations is especially problematic since computations of equilibrium prices must often be performed using models with $\varsigma_* < \infty$. Such models may be derived from simulation-based approximations to Mixed Logit models with reservation prices that are finite $\mu$-a.e., as in the \cite{Berry95} model of Example \ref{EX:BLPExample}. 
	
	Fortunately methods based on the $\bsym{\zeta}$ map resolve {\em only} equilibria of the original problem. In Section \ref{ECSUBSEC:ResolutionOfEquilibria} we consider the important class of simulation-based approximations to Mixed Logit models like those from Example \ref{EX:BLPExample} and prove that fixed-points of $\vec{c} + \bsym{\zeta}(\cdot)$ {\em cannot} be equilibria of a sub-problem that is not an equilibria of the original model. This is essentially a consequence of Eqn. (\ref{Zeta-CGRelation}), which connects the sign of $(D_k\hat{\pi}_f)(\vec{p})$ directly to the sign of $p_k - c_k - \zeta_k(\vec{p})$. 

	Similar results may apply to the markup equation. However because Eqn. (\ref{Eta-CGRelation}) involves $(\tilde{D}\vec{P})(\vec{p})^\top$ instead of simply the diagonal matrix $\bsym{\Lambda}(\vec{p})$, the relationship between the sign of $p_k - c_k - \eta_k(\vec{p})$ and the sign of $(D_k\hat{\pi}_f)(\vec{p})$ is not clear.

	
	\subsubsection{General Results.} 
	
	We now prove the results stated above concerning a variational formulation of the price equilibrium problem when $\varsigma_* < \infty$. 
	
	\begin{proposition}
		\label{InfinitySolvesVI}
		Suppose $\varsigma_* < \infty$ and Assumptions \ref{MixedLogitUtilityAssump}-\ref{LVA} hold. Then the variational inequality (\ref{VI}) always contains $\varsigma_*\vec{1} \in \set{P}^J$ as a solution. 
	\end{proposition}
	\proof{}
		Since $(\tilde{\nabla}\hat{\pi})(\varsigma_*\vec{1}) = \vec{0}$, Eqn. (\ref{VI}) is trivially satisfied.
	\endproof
	
	The following proposition states that this variational formulation is poorly posed in the sense that it contains solutions to all sub-problems. 
	\begin{proposition}
		\label{SubEquilibriaSolveVI}
		Let $\varsigma_* < \infty$ and Assumptions \ref{MixedLogitUtilityAssump}-\ref{LVA} hold. Consider a proper subset $\set{J}^\prime \subset \N(J)$ of $J^\prime = \abs{\set{J}^\prime}$ product indices, and any solution $\vec{p}_{\set{J}^\prime}^* = \{ p_j^* : j \in \set{J}^\prime )$ to the sub-variational inequality 
		\begin{equation*}
			\sum_{j \in \set{J}^\prime} 
				(D_j\hat{\pi}_{f(j)})(\vec{p}_{\set{J}^\prime}^*) 
				(p_j^* - q_j )
				\geq 0
			\quad\text{for all}\quad
			\vec{q}_{\set{J}^\prime}
				= \{ q_j : j \in \set{J}^\prime \}
				\subset [0,\varsigma_*]^{J^\prime}. 
		\end{equation*}
		If we define $\vec{p} \in [0,\varsigma_*]^J$ by $p_j = p_j^*$ for all $j \in \set{J}^\prime$ and $p_k = \varsigma_*$ for all $k \notin \set{J}^\prime$ then $\vec{p}$ solves the full variational inequality (\ref{VI}). 
	\end{proposition}
	\proof{}
		Because
		\begin{equation*}
			(D_j\hat{\pi}_{f(j)})(\vec{p})
				= \left\{ \begin{aligned}
					&(D_j\hat{\pi}_{f(j)})(\vec{p}_{\set{J}^\prime}^*)
						&&\quad\text{if } j \in \set{J}^\prime \\
					&\quad\quad 0
						&&\quad\text{if } j \notin \set{J}^\prime \\
				\end{aligned} \right .
		\end{equation*}
		we have
		\begin{equation*}
			\sum_{j = 1}^J (D_j\hat{\pi}_{f(j)})(\vec{p}) (p_j - q_j )
				= \sum_{j \in \set{J}^\prime} 
					(D_j\hat{\pi}_{f(j)})(\vec{p}_{\set{J}^\prime}^*) 
					(p_j^* - q_j )
				\geq 0
		\end{equation*}
		for all $\vec{q} \in [0,\varsigma_*]^J$.
	\endproof
	

	
	\subsubsection{The Resolution of Equilibria with $\bsym{\zeta}$}
	\label{ECSUBSEC:ResolutionOfEquilibria}
	
	We have shown that variational formulations of the equilibrium problem nest equilibria of all sub-problems, which may {\em not} be equilibria of the original problem as Example \ref{VIExample}. In this section we show that methods based on the $\bsym{\zeta}$ map need not have this unfortunate shortcoming. This result strongly distinguishes nonlinear system methods based on the $\bsym{\zeta}$ map from variational approaches. 
	
	We motivate this result with an example. 
	\begin{example}
		\label{EX:ResolutionExample}
		Consider a finite-sample approximation to the \cite{Berry95} model of Example \ref{EX:BLPExample}. That is, choose $S \in \N$ and draw $\{ \bsym{\theta}_s \}_{s=1}^S$ where $\bsym{\theta}_s = ( \phi_s , \bsym{\beta}_s , \beta_{0,s} )$. These samples could be drawn via standard sampling from $\mu$ or from another technique like importance or quasi-random sampling. In any case, suppose that the $\phi$'s drawn are distinct with probability one: $\phi_s \neq \phi_r$ for all $s,r \in \N(S)$ with probability one. Without loss of generality we take $\phi_1 < \phi_2 < \dotsb < \phi_S$, and note that $\varsigma_* = \phi_S < \infty$. 
		If $\vec{p} = \vec{c} + \bsym{\zeta}(\vec{p})$ and $p_k > \varsigma_*$, then firm $f(k)$'s profits increase with the price of the $k\ith$ product in some neighborhood of $\varsigma_*$.
	\end{example}
	
	Thus if we compute some fixed-point $\vec{p} = \vec{c} + \bsym{\zeta}(\vec{p})$ with $p_k > \varsigma_*$, we know that excluding product $k$ is profit-optimal for firm $f(k)$. As shown in Example \ref{VIExample}, this is not the case with the VI formulation. 
	
	\proof{}
		We will first define $\bsym{\zeta}$ on all of $\set{P}^J$, and then consider fixed-points $\vec{p} = \vec{c} + \bsym{\zeta}(\vec{p})$ with $p_k \geq \phi_S = \varsigma_*$.  
		
		To extend $\bsym{\zeta}$, we define
		\begin{equation*}
			\zeta_k(p_1,\dotsc,p_k,\dotsc,p_J)
				= \zeta_k(p_1,\dotsc,\varsigma_*,\dotsc,p_J)
				= \lim_{ q \to \varsigma_* } \zeta_k(p_1,\dotsc,q,\dotsc,p_J). 
		\end{equation*}
		when $p_k \geq \varsigma_*$. Note that for all $k$ and all $\vec{p} \in (0,\varsigma_*)^J$ we can write
		\begin{equation*}
			\zeta_k(\vec{p})
				= \sum_{ s : \phi_s > p_k }
					\left( 
						\sum_{j\in\set{J}_{f(k)}} P_j^L(\bsym{\theta}_s,\vec{p})(p_j-c_j)
							+ \frac{\phi_s - p_k}{\alpha} 
					\right)
					\left( \frac{ P_k^L(\bsym{\theta}_s,\vec{p}) / (\phi_s - p_k) }
						{ \sum_{ r : \phi_r > p_k } P_k^L(\bsym{\theta}_r,\vec{p})/(\phi_r - p_k) }
					\right)
		\end{equation*}
		We first define $\lim_{p_k \uparrow \phi_S} \zeta_k(\vec{p})$, we first note that for all $p_k \in ( \phi_{S-1} , \phi_S ) \neq \{ \emptyset \}$, we have
		\begin{equation*}
			\zeta_k(\vec{p})
				= \sum_{ j \in \set{J}_{f(k)} } P_j^L(\bsym{\theta}_S,\vec{p}) ( p_j - c_j )
							+ \frac{\phi_S - p_k}{\alpha}
		\end{equation*}
		since $p_k > \phi_s$ for all $s \in \{1,\dotsc,S-1\}$. Thus
		\begin{equation*}
			\lim_{p_k \uparrow \phi_S} \zeta_k(\vec{p})
				= \sum_{ j \in \set{J}_{f(k)} \setminus \{ k \} } 
					\left[ \lim_{p_k \uparrow \phi_S} P_j^L(\bsym{\theta}_S,\vec{p}) \right] 
					( p_j - c_j ). 
		\end{equation*}
		In other words, as $p_k$ approaches $\phi_S = \varsigma_*$, $\zeta_k$ approaches the profits firm $f(k)$ accrues from selling all products {\em other} than $p_k$ to the sampled individual with the highest income. This establishes that the extended $\bsym{\zeta}$ is well-defined and continuous. 
		
		Now suppose $p_k = c_k + \zeta_k(\vec{p})$, where $p_k > \phi_S = \varsigma_*$. Thus
		\begin{equation*}
			0 = p_k - c_k - \zeta_k(\vec{p})
				> \phi_S - c_k - \zeta_k(\vec{p}) 
				= \lim_{q_k \uparrow \phi_S} \Big( q_k - c_k - \zeta_k(p_1,\dotsc,q_k,\dotsc,p_J) \Big), 
		\end{equation*}
		and there must exist some $\delta > 0$ such that 
		\begin{equation*}
			q_k - c_k - \zeta_k(p_1,\dotsc,q_k,\dotsc,p_J) < 0
		\end{equation*}
		for all $q_k \in (\varsigma_* - \delta,\varsigma_*)$. Hence $(D_k\hat{\pi}_{f(k)})(p_1,\dotsc,q_k,\dotsc,p_J) > 0$
		\begin{equation*}
			(D_k\hat{\pi}_{f(k)})(p_1,\dotsc,q_k,\dotsc,p_J) 
				= \lambda_k(p_1,\dotsc,q_k,\dotsc,p_J)
					\big( q_k - c_k - \zeta_k(p_1,\dotsc,q_k,\dotsc,p_J) \big) > 0.
		\end{equation*}
		In other words, if $\vec{p} = \vec{c} + \bsym{\zeta}(\vec{p})$ and $p_k > \varsigma_*$, then firm $f(k)$'s profits increase with the price of the $k\ith$ product in some neighborhood of $\varsigma_*$.
	\endproof
		
%
	
	Fortunately this example is fairly general. In the following proposition we prove that all finite-sample simulators generate $\bsym{\zeta}$ maps that do not have equilibria of sub-problems as fixed points {\em unless} they are, in fact, equilibria of the original problem. Three assumptions are added: utilities must be twice continuously differentiable in prices, $\varsigma(\bsym{\theta})$ is finite $\mu$-a.e. as in the \cite{Berry95} model, and the sampled values $\varsigma(\bsym{\theta}_s)$ must be distinct with probability one.
	
	\begin{proposition}
		Consider a Mixed Logit model satisfying Assumptions \ref{MixedLogitUtilityAssump}, \ref{LeibnizRuleCondition}, and \ref{LVA} with $w_j (\bsym{\theta},\cdot) : (0,\varsigma(\bsym{\theta})) \to \R$ twice continuously differentiable in price and $\varsigma : \set{T} \to \set{P}$ finite $\mu$-a.e.. 
		
		Generate a finite-sample simulator to this Mixed Logit model with $\{ \bsym{\theta}_s \}_{s=1}^S$ for some $S \in \N$. Let $\varsigma_s = \varsigma(\bsym{\theta}_s)$, and assume that $\varsigma_s \neq \varsigma_r$ with probability one for any $s \neq r$. Subsequently, order the samples so that $\varsigma_1 < \dotsb < \varsigma_S = \varsigma_*$. 
		
		Suppose that $\vec{p} \in \set{P}^J$ satisfies $\vec{p} = \vec{c} + \bsym{\zeta}(\vec{p})$ where $\bsym{\zeta}$ is the extended map as in Example \ref{EX:ResolutionExample}. If $p_k \geq \varsigma_S$, then excluding product $k$ is profit-optimal for firm $f = f(k)$; particularly, there exists $\delta > 0$ such that
		\begin{equation*}
			(D_k\hat{\pi}_{f(k)})(p_1,\dotsc,p_k,\dotsc,p_J) > 0
		\end{equation*}
		for all $p_k \in (\varsigma_S - \delta,\varsigma_S)$. 
	\end{proposition}
	\proof{}
		The case $p_k > \varsigma_S$ is handled exactly as in Example \ref{EX:ResolutionExample}. We must only consider the case where
		\begin{equation*}
			0 = \varsigma_S -  c_k - \lim_{p_k \uparrow \varsigma_S} \zeta_k(\vec{p}) 
				= \lim_{p_k \uparrow \varsigma_S} \Big[ p_k - c_k - \zeta_k(\vec{p}) \Big]. 
		\end{equation*} 
		Our approach is to show that $D_k [ p_k - c_k - \zeta_k(\vec{p})] > 0$ for all $p_k$ near enough to $\varsigma_S$, and thus
		\begin{align*}
			p_k - c_k - \zeta_k(\vec{p})
				&= p_k - c_k - \zeta_k(\vec{p})
					- \Big[ \varsigma_S - c_k - \lim_{p_k \uparrow \varsigma_S} \zeta_k(\vec{p}) \Big] \\
				&\quad\quad\quad\quad
					= - \int_{p_k}^{\varsigma_S} D_k [ p_k - c_k - \zeta_k(\vec{p})] dp_k
				< 0
		\end{align*}
		(with a slight abuse of notation in the integral). More specifically, we prove that $\lim_{p_k \uparrow \varsigma_S} D_k [ p_k - c_k - \zeta_k(\vec{p})] > 0$, which implies that $D_k [ p_k - c_k - \zeta_k(\vec{p})] > 0$ for all $p_k$ near enough to $\varsigma_S$. Because $p_k - c_k - \zeta_k(\vec{p}) < 0$ for $p_k$ near enough to $\varsigma_S$, 
		\begin{equation*}
			(D_k\hat{\pi}_{f(k)})(p_1,\dotsc,q_k,\dotsc,p_J) 
				= \lambda_k(\vec{p}) \Big( p_k - c_k - \zeta_k(\vec{p}) \Big)
				> 0. 
		\end{equation*}
		
		As in Example \ref{EX:ResolutionExample}, note that for all $p_k \in ( \varsigma_{S-1} , \varsigma_S )$ we have
		\begin{equation*}
			\zeta_k(\vec{p})
				= \sum_{ j \in \set{J}_{f(k)} } P_j^L(\bsym{\theta}_S,\vec{p}) ( p_j - c_j )
							- \frac{1}{(Dw_k)(\bsym{\theta}_S,p_k)}
		\end{equation*}
		since $p_k > \varsigma_s$ for all $s \in \{1,\dotsc,S-1\}$. From this equation we derive
		\begin{equation*}
			(D_k\zeta_k)(\vec{p})
				= \sum_{ j \in \set{J}_{f(k)} } (D_kP_j^L)(\bsym{\theta}_S,\vec{p}) ( p_j - c_j )
					+ P_k^L(\bsym{\theta}_S,\vec{p})
						+ \frac{(D^2w_k)(\bsym{\theta}_S,p_k)}{(Dw_k)(\bsym{\theta}_S,p_k)^2}
		\end{equation*}
		and thus
		\begin{align*}
			&D_k \Big[ p_k - c_k - \zeta_k(\vec{p}) \Big] \\
			&\quad\quad
				= 1 - (Dw_k)(\bsym{\theta}_S,p_k)
						P_k^L(\bsym{\theta}_S,\vec{p})
							\sum_{ j \in \set{J}_{f(k)} } P_j^L(\bsym{\theta}_S,\vec{p}) ( p_j - c_j )
						- P_k^L(\bsym{\theta}_S,\vec{p})
						- \frac{(D^2w_k)(\bsym{\theta}_S,p_k)}{(Dw_k)(\bsym{\theta}_S,p_k)^2}
		\end{align*}
		Now $\lim_{p_k \uparrow \varsigma_S} P_k^L(\bsym{\theta}_S,\vec{p}) = 0$, we have assumed that 
		\begin{equation*}
			\lim_{p_k \uparrow \varsigma_S} 
				\big[ (Dw_k)(\bsym{\theta}_S,p_k)P_k^L(\bsym{\theta}_S,\vec{p}) \big]
				= 0
		\end{equation*}
		(Assumption \ref{LVA}), and 
		\begin{equation*}
			\lim_{p_k \uparrow \varsigma_S} 
				\left[ \sum_{ j \in \set{J}_{f(k)} } P_j^L(\bsym{\theta}_S,\vec{p}) ( p_j - c_j ) \right]
			= \sum_{ j \in \set{J}_{f(k)} \setminus \{k\} } \lim_{p_k \uparrow \varsigma_S} 
				\left[ P_j^L(\bsym{\theta}_S,\vec{p}) \right] ( p_j - c_j ) 
				< \infty,
		\end{equation*}
		we have
		\begin{align*}
			\lim_{p_k \uparrow \varsigma_S} D_k \Big[ p_k - c_k - \zeta_k(\vec{p}) \Big] 
				= 1 - \lim_{p_k \uparrow \varsigma_S}
						\left[ \frac{(D^2w_k)(\bsym{\theta}_S,p_k)}
								{(Dw_k)(\bsym{\theta}_S,p_k)^2} \right]
		\end{align*}
		So long as 
		\begin{align*}
			\lim_{p_k \uparrow \varsigma_S}
						\left[ \frac{(D^2w_k)(\bsym{\theta}_S,p_k)}
								{(Dw_k)(\bsym{\theta}_S,p_k)^2} \right]
				< 1
		\end{align*}
		we have $\lim_{p_k \uparrow \varsigma_S} D_k \Big[ p_k - c_k - \zeta_k(\vec{p}) \Big] > 0$. This must be true, as Claim \ref{MustBeSubQuadratic} below demonstrates. This completes the proof.
		\endproof
		
		\begin{claim}
			\label{MustBeSubQuadratic}
			Let $w : (0,\varsigma) \to \R$ be twice continuously differentiable, with $(Dw)(p) < 0$ for all $p \in (0,\varsigma)$ and $(Dw)(p) \downarrow -\infty$ as $p \uparrow \varsigma$. Then 
			\begin{equation*}
				\lim_{p \uparrow \varsigma} \left[ \frac{(D^2w)(p)}{(Dw)(p)^2} \right] < 1. 
			\end{equation*}
		\end{claim}
		\proof{Proof}
			We prove this by contradiction. Note that
			\begin{equation*}
				D \left[ \frac{1}{\abs{(Dw)(p)}} \right]
					= \frac{(D^2w)(p)}{(Dw)(p)^2}. 
			\end{equation*}
			Now if 
			\begin{align*}
				\lim_{p \uparrow \varsigma_S} \left[ \frac{(D^2w)(p)}{(Dw)(p)^2} \right]
					\geq 1, 
			\end{align*}
			there must exist some $\bar{p} \in (0,\varsigma)$ such that
			\begin{align*}
				\frac{(D^2w)(p)}{(Dw)(p)^2} > 0 \quad\text{for all}\quad p \in [\bar{p},\varsigma). 
			\end{align*}
			But then
			\begin{equation*}
				0 \leq \int_p^\varsigma \frac{(D^2w)(q)}{(Dw)(q)^2} dq	
					= \int_p^\varsigma D \left[ \frac{1}{\abs{(Dw)(q)}} \right] dq
					= \lim_{q \uparrow \varsigma} \left[ \frac{1}{\abs{(Dw)(q)}} \right]
						- \frac{1}{\abs{(Dw)(p)}}
					= - \frac{1}{\abs{(Dw)(p)}} < 0, 
			\end{equation*}
			a contradiction.
		\endproof
	
	
	\subsection{Tatonnement} 
	
	Some authors iterate best responses $-$ i.e. tatonnement $-$ to compute equilibria. See, for example, \cite{Choi90, CBO03, Michalek04, Austin05, Bento05, Hu07}. For this process Newton's method, or another algorithm of (unconstrained) optimization, will be required. Tatonnement should be an efficient way to compute ``equilibrium'' if all firm's profit-maximizing prices are independent of their competitor's decisions, but wasteful if some firm's optimal pricing depends heavily on their competitors' prices. Furthermore no convergence guarantees exist for tatonnement while there are at least theoretical guarantees that Newton's method, properly constructed, will converge to simultaneously stationary prices. 
	
	
	\subsection{Least-Squares Minimization and the Gauss-Newton Method} 
	
	In principle one could also use optimization methods to {\em explicitly} minimize $f(\vec{p}) = \norm{\vec{F}(\vec{p})}_2^2/2$ for any of our choices of $\vec{F}$. In fact, line search and trust-region strategies for global convergence {\em implicitly} minimize this function \citep{Dennis96}. Computations of equilibrium prices benefit from leaving this implicit, as explicit minimization via Newton's method requires third-order derivatives of $\vec{F}$, increasing both differentiability requirements and computational burden. The Gauss-Newton method \citep{Ortega70} is obtained by neglecting the influence of the third-order derivatives of $\vec{F}$. This defines the Gauss-Newton step as a solution to the (symmetric) normal equation $(D\vec{F})(\vec{p})^\top(D\vec{F})(\vec{p})\vec{s} = - (D\vec{F})(\vec{p})^\top\vec{F}(\vec{p})$; note that the same problem arises should one wish to use the Conjugate Gradient method to solve the Newton system. So long as $(D\vec{F})(\vec{p})$ is nonsingular the standard Newton steps will be recovered from the Gauss-Newton method. However they are explicitly formulated as solutions to linear systems that are more poorly conditioned \citep{Golub96,Trefethen97} and thus we should at least expect to accumulate more error in the process of solving for the same steps. The burden of computing these steps also increases because of the requirement to multiply by the transpose of the Jacobian of $\vec{F}$.